
\magnification=1100
\overfullrule0pt

\input prepictex
\input pictex
\input postpictex

\input amssym.def


\def\qed{\hbox{\hskip 1pt\vrule width4pt height 6pt depth1.5pt \hskip 1pt}}
\def\mapleftright#1{\smash{
   \mathop{\longleftrightarrow}\limits^{#1}}}

\def\CC{{\Bbb C}}
\def\RR{{\Bbb R}}
\def\ZZ{{\Bbb Z}}
\def\cA{{\cal A}}

\def\cF{{\cal F}}
\def\fh{{\frak h}}
\def\gp{{\goth p}}




\font\smallcaps=cmcsc10
\font\titlefont=cmr10 scaled \magstep1

\font\sectionfont=cmbx10
\font\tinyrm=cmr10 at 8pt


\newcount\sectno
\newcount\subsectno
\newcount\resultno

\def\section #1. #2\par{
\sectno=#1
\resultno=0
\bigskip\noindent{\sectionfont #1.  #2}~\medbreak}

\def\subsection #1\par{ \global\advance\resultno by 1
\bigskip\noindent{\bf  (\the\sectno.\the\resultno) #1}\ \ }


\def\prop{ \global\advance\resultno by 1
\bigskip\noindent{\bf Proposition \the\sectno.\the\resultno. }\sl}
\def\lemma{ \global\advance\resultno by 1
\bigskip\noindent{\bf Lemma \the\sectno.\the\resultno. }
\sl}
\def\remark{ \global\advance\resultno by 1
\bigskip\noindent{\bf Remark \the\sectno.\the\resultno. }}
\def\example{ \global\advance\resultno by 1
\bigskip\noindent{\bf Example \the\sectno.\the\resultno. }}
\def\cor{ \global\advance\resultno by 1
\bigskip\noindent{\bf Corollary \the\sectno.\the\resultno. }\sl}
\def\thm{ \global\advance\resultno by 1
\bigskip\noindent{\bf Theorem \the\sectno.\the\resultno. }\sl}
\def\defn{ \global\advance\resultno by 1
\bigskip\noindent{\it Definition \the\sectno.\the\resultno. }\slrm}
\def\endthm{\rm\bigskip}

\def\endlemma{\rm\bigskip}
\def\endexample{\rm\bigskip}
\def\endprop{\rm\bigskip}

\def\pf{\rm\bigskip\noindent{\it Proof. }}
\def\endpf{\qed\hfil\bigskip}


\def\formula{\global\advance\resultno by 1
\eqno{(\the\sectno.\the\resultno)}}
\def\formulano{\global\advance\resultno by 1 (\the\sectno.\the\resultno)}
\def\tableno{\global\advance\resultno by 1
\the\sectno.\the\resultno. }
\def\lformula{\global\advance\resultno by 1
\leqno(\the\sectno.\the\resultno)}

\def\monthname {\ifcase\month\or January\or February\or March\or April\or
May\or June\or
July\or August\or September\or October\or November\or December\fi}

\newcount\mins  \newcount\hours  \hours=\time \mins=\time
\def\now{\divide\hours by60 \multiply\hours by60 \advance\mins by-\hours
     \divide\hours by60         
     \ifnum\hours>12 \advance\hours by-12
       \number\hours:\ifnum\mins<10 0\fi\number\mins\ P.M.\else
       \number\hours:\ifnum\mins<10 0\fi\number\mins\ A.M.\fi}


\nopagenumbers
\def\runningtitle{\smallcaps affine hecke algebras}
\headline={\ifnum\pageno>1\eoheadline\else\firstheadline\fi}
\def\names{\smallcaps arun ram}
\def\firstheadline{\noindent Preliminary Draft \hfill  \today}
\def\firstheadline{}
\def\eoheadline{\ifodd\pageno\oddheadline\else\evenheadline\fi}
\def\oddheadline{\tenrm\hfil\runningtitle\hfil\folio}
\def\evenheadline{\tenrm \folio\hfil{\names}\hfil}

\vphantom{$ $}  
\vskip.75truein

\centerline{\titlefont Affine Hecke algebras and generalized standard Young
tableaux}
\bigskip
\centerline{\rm Arun Ram${}^\ast$}
\centerline{Department of Mathematics}
\centerline{University of Wisconsin}
\centerline{Madison, WI 53706 USA}
\centerline{{\tt ram@math.wisc.edu}}

\bigskip
\centerline{\sl  Dedicated to Robert Steinberg}

\footnote{}{\tinyrm ${}^\ast$ Research supported in part by the National
Science Foundation (DMS-0097977) and the National Security Agency (MDA904-01-1-0032).}

\bigskip

\noindent{\bf Abstract.}
This paper introduces calibrated representations 
for affine Hecke algebras and classifies and constructs all finite
dimensional irreducible calibrated representations.  The primary technique
is to provide indexing sets for controlling the weight space structure
of finite dimensional modules for the affine Hecke algebra.
Using these indexing sets we show that 
(1) irreducible calibrated representations are indexed by skew local regions,
(2)  the dimension of an irreducible calibrated representation is 
the number of chambers in the local region,
(3)  each irreducible calibrated representation is constructed explicitly
by formulas which describe the action of the generators
of the affine Hecke algebra on a specific basis in the representation space.  
The indexing sets for weight spaces are generalizations of standard 
Young tableaux
and the construction of the irreducible calibrated affine Hecke algebra modules
is a generalization of A. Young's seminormal construction of the irreducible
representations of the symmetric group.  In this sense Young's
construction has been generalized to arbitrary Lie type.

\section 0. Introduction

The classical representation theory of the symmetric
group, as developed by G. Frobenius and 
A. Young [Yg1,2], has the following features:
\smallskip\noindent
\item{(a)}  The irreducible representations $S^\lambda$
of the symmetric group $S_n$ 
are indexed by partitions $\lambda$ with $n$ boxes,
\smallskip\noindent
\item{(b)}  The dimension of $S^\lambda$ is the number of
standard tableaux of shape $\lambda$,
\smallskip\noindent
\item{(c)} The $S_n$-module has an elegant
explicit construction: $S^\lambda$ is the span of a basis
$\{v_T\}$ parametrized by standard tableaux $T$
and the action of each generator of $S_n$ is given by a
simple formula,
$$
s_i v_T = {1\over c(T(i))-c(T(i+1))}v_T
+ \left(1+ {1\over c(T(i))-c(T(i+1))}\right)v_{s_iT}.$$

In this paper we prove analogous results for representations
of affine Hecke algebras.
\smallskip\noindent
\item{(A)}  The irreducible calibrated representations 
$\tilde H^{(t,J)}$ of the affine Hecke algebra $\tilde H$
are indexed by skew local regions $(t,J)$,
\smallskip\noindent
\item{(B)} The dimension of $\tilde H^{(t,J)}$ is the number
of chambers in the local region $(t,J)$,
\smallskip\noindent
\item{(C)} The $\tilde H$-module $\tilde H^{(t,J)}$ has an
elegant explicit construction:
$\tilde H^{(t,J)}$ is the span of a basis 
$\{v_w\ |\ w\in {\cal F}^{(t,J)}\}$
parametrized by chambers in the local region
and the action of each generator of $\tilde H$ is
given by a simple formula,
$$\eqalign{
X^\lambda v_w &= q^{\langle\lambda,w\gamma\rangle}v_w, \cr
T_i v_w &= 
{q-q^{-1}\over 1-t(X^{w^{-1}\alpha_i}) }
+
\left(q^{-1}+
{q-q^{-1}\over 1-t(X^{w^{-1}\alpha_i}) }
\right)v_{s_iw}. \cr
}$$

\noindent
In fact, the classical theory of standard Young tableaux and
partitions is a special case of our theory of chambers 
and local regions;  
this is proved in Sections 5 and 6 of this paper.
Section 1 serves to fix notations and fundamental
data in the form which will need it.  
The bulk of this material can be found in [Bou, Ch. IV-VI]
and Steinberg's Yale Lecture Notes [Sb1].
Two known results which are included in Section 1 
are: (a) the determination of the
center of the affine Hecke algebra, and 
(b) the Pittie-Steinberg theorem, which provides a nice
basis for the affine Hecke algebra over its center.
In each case we have given an elementary proof, which,
hopefully, illustrates the beautiful simplicity of these
powerful results.  Section 2 treats the notion of weight
spaces for affine Hecke algebra representations and 
shows how certain combinatorially defined
indexing sets ${\cal F}^{(t,J)}$ give explicit information
about the weight space structure of affine Hecke algebra modules.
Section 3 classifies and constructs all irreducible calibrated
affine Hecke algebra modules (for any $q$ such that 
$q^2\ne \pm 1$, {\it including roots of unity}).
Section 4 gives the main results about the 
structure of the labeling sets ${\cal F}^{(t,J)}$ and
defines a conjugation involution on them.  Sections 5 and
6 show that the classical theory of standard Young tableaux is
very special case of the analysis of the combinatorial structure
of the sets ${\cal F}^{(t,J)}$.
Section 7 works out the generalized standard Young
tableaux in the type A, root of unity case.  The resulting
objects are $\ell$-periodic standard Young tableaux.
Section 8 describes how the generalized standard Young
tableaux look in the type C, non root of unity case.  In
this case the objects are negative rotationally symmetric standard
Young tableaux.  It should not be difficult to work out similar 
explicit tableaux in terms of fillings of boxes in the other classical
types.

\smallskip
Let us put these results into perspective.
\smallskip\noindent
{\sl (1) $p$-adic groups and affine Hecke algebras.}
\smallskip\noindent
The affine Hecke algebra was introduced by Iwahori and Matsumoto [IM]
as a tool for studying the representations of a $p$-adic Lie group.
In some sense, all irreducible principal series representations of
the $p$-adic group can be determined by classifying the representations
of the corresponding affine Hecke algebra.  
Kazhdan and Lusztig [KL] (see also [CG])
gave a geometric classification of all irreducible representations of
the affine Hecke algebra. This classification is a $q$-analogue of 
Springer's construction of the irreducible representations of the Weyl 
group on the cohomology
of unipotent varieties.  In the $q$-case, K-theory takes the place
of cohomology and the irreducible representations of the affine
Hecke algebra are constructed
as quotients of the K-theory of the Steinberg varieties.
It is difficult to obtain combinatorial
information from this geometric construction.  So the combinatorial
approach in this paper gives new information.

\medskip\noindent
{\sl (2) The theory of Young tableaux.}
\smallskip\noindent
The word ``Young tableau'' 
is commonly used for three very different
objects in representation theory:
\smallskip\noindent
\itemitem{(1a)} {\it partitions} with $n$ boxes, which index representations
of the symmetric group $S_n$,
\smallskip\noindent
\itemitem{(1b)} {\it partitions with $\le n$ rows}, which index the
polynomial representations of $GL_n(\CC)$,
\smallskip\noindent
\itemitem{(2)} {\it standard tableaux}, which label the basis elements
of an irreducible representation of $S_n$,
\smallskip\noindent
\itemitem{(3)} {\it column strict tableaux}, which label the basis elements
of an irreducible polynomial representation of $GL_n(\CC)$.
\smallskip\noindent
The partitions in (1b) were generalized to all Lie types
by H. Weyl in 1926, who showed that finite dimensional irreducible
representations of compact Lie groups are indexed by the
dominant integral weights.  There was much important work
generalizing the column strict tableaux in (3) to
other Lie types, for a survey of this work see [Su].  The
problem of generalizing the column strict tableaux in (3) 
to {\it all} Lie types was finally solved by the path model 
of Littelmann [Li1-2].  This
paper provides a generalization of the partitions of (1a)
and the standard tableaux of (2) which are valid for 
{\it all Lie types}.
For important earlier work in this direction see
[Mac, I App. B], Hoefsmit [Ho], and Ariki and Koike [AK].

This paper is a revised, expanded, and updated version of the
preprints [R2] and [R3].  The original preprints will not
be published since the results there are contained in 
and expanded in this paper.  Those preprints will 
remain available at {\tt http://www.math.wisc.edu/\~{\ }ram/preprints.html}.

\medskip\noindent
{\it Acknowledgements}

\medskip
During this work I have benefited from
conversations with many people, including,
but not limited to,
G. Benkart, 
H. Barcelo,
P. Deligne, 
S. Fomin, 
T. Halverson, 
F. Knop, 
R. Macpherson, 
R. Simion,
L. Solomon, 
J. Stembridge, 
M. Vazirani, 
D.-N. Verma,
and N. Wallach.  I sincerely thank everyone who has let
me tell them my story. Every one of these sessions was helpful to 
me in solidifying my understanding.
I thank A. Kleshchev for thrilling energetic
conversations which pushed me to work the examples out 
carefully for type A root of unity case and I thank
J. Olsson for his wonderful gift to me of A. Young's
collected papers [Yg1].

\bigskip

\section 1. The affine Hecke algebra

\medskip
Though we shall never really use the data $(G,B,T)$ it is
conceptually useful to note that
there is an affine Hecke algebra associated to each triple 
$(G\supseteq B\supseteq T)$ where 
\smallskip\noindent
\itemitem{} $G$ is a connected reductive complex algebraic group,
\smallskip\noindent
\itemitem{} $B$ is a Borel subgroup,
\smallskip\noindent
\itemitem{} $T$ is a maximal torus.
\smallskip\noindent
An example of this data is when
$G=GL_n(\CC)$, $B$ is
the subgroup of upper triangular invertible
matrices, and $T$ is the subgroup of invertible diagonal matrices.

The reason that we can avoid the 
data $(G\supseteq B\supseteq T)$ is that it is equivalent to a different
data $(W,C,L)$ where
\smallskip\noindent
\itemitem{} $W$ is a finite real reflection group with
reflection representation $\fh_\RR^*$,
\smallskip\noindent
\itemitem{} $C$ is a fixed fundamental chamber for the $W$-action,
\smallskip\noindent
\itemitem{} $L$ is a $W$-invariant lattice in $\fh_\RR^*$.
\smallskip\noindent
This will be our basic data.
In the example where $G=GL_n(\CC)$ and $B$ and $T$ are the
upper triangular and diagonal matrices, respectively,
$$
W = S_n,
\quad
\fh_\RR^* = \RR^n = \sum_{i=1}^n \RR\varepsilon_i,
\quad
C = \Big\{\mu = \sum_{i=1}^n \mu_i\varepsilon_i
\ \Big\vert\ 
\mu_1\le\cdots\le \mu_n\Big\},
\quad
L = \sum_{i=1}^n \ZZ\varepsilon_i,
\formula$$
where $W=S_n$ is the symmetric group, acting on
$\fh_{\RR}^*=\RR^n$ by permuting
the orthonormal basis $\varepsilon_1,\ldots, \varepsilon_n$.
This example will be treated in depth in Sections 5, 6 and 7.  
We shall show that the labeling sets ${\cal F}^{(t,J)}$
for weight spaces of affine Hecke algebra representations 
that are introduced in (2.18) and Corollary 2.19 and used for the classification
in Theorem 3.6 are generalizations of standard Young tableaux.

The components $W$ and $L$ in the data $(W,C,L)$ are 
obtained from $(G\supseteq B\supseteq T)$ by 
$$W = N(T)/T, \qquad X={\rm Hom}(T,\CC^*)= \{X^\lambda\ |\ \lambda\in L\},$$
where $N(T)$ is the normalizer of $T$ in $G$ and
${\rm Hom}(T,\CC^*)$ is the set of algebraic group homomorphisms
from $T$ to $\CC^*$.  The notation is designed so that the multiplication
in the group $X$ is 
$$X^\lambda X^\mu = X^{\lambda+\mu} = X^\mu X^\lambda,
\qquad\hbox{for $\mu,\lambda\in L$,}\formula$$
see [CSM, III \S 8].  
The reflection (or defining) representation of the group $W$ 
is given by its action on 
$\fh_\RR^* = \RR\otimes_\ZZ L \cong \RR^n$
and with respect to a $W$-invariant inner product $\langle,\rangle$
on $\fh_\RR^*$ the group $W$ is generated by reflections $s_\alpha$ in the 
hyperplanes 
$$H_\alpha = \{ x\in \fh_{\RR}^*\ |\ \langle x,\alpha\rangle=0\},
\qquad \alpha\in R^+.\formula$$
See the picture which appears just before Theorem 1.17.
The {\it chambers} are the connected components of 
$\fh_{\RR}^*-\left(\bigcup_{\alpha\in R^+} H_\alpha\right)$
and these are the fundamental regions for the action of $W$ on $\fh_\RR^*$.
Fixing a choice of a fundamental chamber $C$ corresponds to the choice
of the set $R^+$ of positive roots, 
which corresponds to the choice of $B$ in $G$.

In our formulation we may view the set $R^+$ as a labeling set for the
reflecting hyperplanes $H_\alpha$ in $\fh_\RR^*$ so that
$$C = \{ x\in \fh_{\RR}^*\ |\ \langle x,\alpha\rangle>0
\hbox{\ for all $\alpha\in R^+$}\}.\formula$$
For a root $\alpha\in R$, the
{\it positive side} of the hyperplane $H_\alpha$ is the side towards
$C$, i.e.\ $\{\lambda\in \fh_{\RR}^*\ |\ \langle\lambda,\alpha\rangle>0\}$,
and the {\it negative side} of $H_\alpha$ is the side away from $C$.  

For $w\in W$, the {\it inversion set} of $W$ is
$$R(w) = \{ \alpha\in R^+\ |\ w\alpha\in R^-\},\formula$$
where $R^- = -R^+$.  There is a bijection
$$\matrix{
W &\longleftrightarrow 
&\{\hbox{fundamental chambers for $W$ acting on $\fh_{\RR}^*$}\} \cr
w &\longmapsto &w^{-1}C \cr
}
\formula$$
and the chamber $w^{-1}C$ is the unique chamber which is on the 
positive side of $H_\alpha$ for $\alpha\not\in R(w)$ and on the
negative side of $H_\alpha$ for $\alpha\in R(w)$.

The {\it simple roots} $\alpha_1,\ldots, \alpha_n$ in $R^+$ 
index the walls $H_{\alpha_i}$ of the fundamental chamber $C$ 
and the corresponding reflections $s_1, \ldots s_n$ generate
$W$.  In fact, $W$ can be presented by 
generators $s_1,s_2,\ldots, s_n$ and relations
$$\matrix{
s_i^2&=&1, &&\hbox{for $1\le i\le n$,} \cr
\underbrace{s_is_js_i\cdots}_{m_{ij} {\rm \ \ factors}}
&=& \underbrace{s_js_is_j\cdots}_{m_{ij} {\rm \ \ factors}}\;,
&\qquad &\hbox{for $i\ne j$,}  \cr
}
\formula$$
where the (acute) angle $\pi/m_{ij}$ between the 
hyperplanes $H_{\alpha_i}$ and $H_{\alpha_j}$ determines
the value $m_{ij}$.

Fix $q\in \CC^*$ with $q^2\ne \pm 1$.
The {\it Iwahori-Hecke algebra}
$H$ associated to $(W,C)$
is the associative algebra over $\CC$ defined by 
generators $T_1,T_2,\ldots, T_n$ and relations
$$\matrix{
T_i^2&=& (q-q^{-1})T_i+1, &&\hbox{for $1\le i\le n$,} \cr
\underbrace{T_iT_jT_i\cdots}_{m_{ij} {\rm \ \ factors}}
&=& \underbrace{T_jT_iT_j\cdots}_{m_{ij} {\rm \ \ factors}}\;,
&\qquad &\hbox{for $i\ne j$,}  \cr
}
\formula$$
where $m_{ij}$ are the same as in the presentation of $W$.
For $w\in W$ define $T_w=T_{i_1}\cdots T_{i_p}$ where
$s_{i_1}\cdots s_{i_p} = w$ is a reduced
expression for $w$.  By [Bou, Ch. IV \S 2 Ex. 23], the element
$T_w$ does not depend on the choice of the reduced expression.
The algebra $H$ has dimension $|W|$ and the set $\{T_w\}_{w\in W}$
is a basis of $H$.

The {\it affine Hecke algebra} $\tilde H$ associated to $(W,C,L)$
algebra given by
$$\tilde H  
= \CC\hbox{-span} \{ T_w X^\lambda \ |\ w\in W, X^\lambda\in X\}
\formula$$ 
where the multiplication of the $T_w$ is as in the Iwahori-Hecke algebra
$H$, the multiplication of the $X^\lambda$ is as in (1.2) and
we impose the relation
$$X^\lambda T_i = T_i X^{s_i\lambda} + (q-q^{-1})
{X^\lambda-X^{s_i\lambda}\over 1-X^{-\alpha_i}}, 
\qquad\hbox{for $1\le i\le n$ and $X^\lambda\in X$.}
\formula$$
This formulation of the definition of $\tilde H$ is due to Lusztig [Lu]
following work of Bernstein and Zelevinsky.  
The elements $T_wX^\lambda$, $w\in W$, $X^\lambda\in X$, form a basis
of $\tilde H$.  

The group algebra of $X$,
$$\CC[X] = \CC\hbox{-span}\{ X^\lambda\ |\ \lambda\in L\},
\formula$$
is a subalgebra of $\tilde H$ with a $W$-action obtained by linearly
extending the $W$-action on $X$,
$$wX^\lambda = X^{w\lambda}
\qquad\hbox{for $w\in W$, $X^\lambda\in X$.}
\formula$$

\thm (Bernstein, Zelevinsky, Lusztig [Lu, 8.1])
The center of $\tilde H$ is  
$\CC[X]^W = \{ f\in \CC[X] \ |\ wf=f \hbox{\ for all $w\in W$} \}$.
\pf
Assume
$$
z=\sum_{\lambda\in L,w\in W} c_{\lambda,w}X^\lambda T_w\in Z(\tilde H).
$$
Let $m\in W$ be maximal in Bruhat order subject to 
$c_{\gamma,m}\neq 0$ for some $\gamma\in L$.
If $m\ne 1$ there exists a dominant $\mu\in L$ such that
$c_{\gamma+\mu-m\mu,m}=0$ (otherwise  $c_{\gamma+\mu-m\mu,m}\neq 0$ for every
dominant $\mu\in L$, which is impossible since $z$ is a finite linear
combination of $X^\lambda T_w$). Since $z\in Z(\tilde H)$
we have
$$
z  =   X^{-\mu}zX^\mu   =
\sum_{\lambda\in L,w\in W} c_{\lambda,w} X^{\lambda-\mu} T_w X^\mu.
$$
Repeated use of the relation (1.10) yields
$$
T_wX^\mu=\sum_{\nu\in L,v\in W} d_{\nu,v}X^\nu T_v
$$
where $d_{\nu,v}$ are constants such that
$d_{w\mu,w}=1$, $d_{\nu,w}=0$ for
$\nu\ne w\mu$, and $d_{\nu,v}=0$ unless $v\le w$.
So
$$
z  = \sum_{\lambda\in L,w\in W} c_{\lambda,w}X^\lambda T_w
= \sum_{\lambda\in L,w\in W}\sum_{\nu\in L,v\in W}
c_{\lambda,w}d_{\nu,v} X^{\lambda-\mu+\nu} T_v
$$
and comparing the coefficients of $X^\gamma T_m$ gives
$
c_{\gamma,m}=c_{\gamma+\mu-m\mu,m} d_{m\mu,m}.
$
Since $c_{\gamma+\mu-m\mu,m}=0$ it follows that
$c_{\gamma,m}=0$, which is a contradiction.  Hence
$z=\sum_{\lambda\in L} c_\lambda X^\lambda\in \CC[X]$.

The relation (1.10) gives
$$
zT_i=T_iz=(s_iz)T_i+(q-q^{-1})z'
$$
where $z'\in\CC[X]$. Comparing coefficients of
$X^\lambda$  on both sides yields $z' = 0$. Hence
$zT_i=(s_iz)T_i$, and therefore $z=s_iz$ for $1\leq i\leq n$.  So
$z\in \CC[X]^W$.
\endpf

It is often convenient to assume that $W$ acts irreducibly
on $\fh_\RR^*$ and that the lattice $L$ is the weight lattice
$$P = \{ x\in \fh_{\RR}^*\ |\ \langle x,\alpha\rangle\in\ZZ
\hbox{\ for all $\alpha\in R^+$}\}
= \sum_{i=1}^n \ZZ\omega_i ,
\formula$$
where the {\it fundamental weights} are the elements 
$\omega_1,\ldots, \omega_n$
of $\RR^n$ given by
$$\langle \omega_i, \alpha_j^\vee\rangle = \delta_{ij},
\qquad\hbox{where}\quad
\alpha_i^\vee ={2\alpha_i\over \langle \alpha_i,\alpha_i\rangle}
\formula
$$
and $\delta_{ij}$ is the Kronecker delta.  Many facts
are easier to state in this case and the general case
can always be reduced to this one. We will make some further
remarks on this reduction at the end of this section.

Consider the connected regions of the negative
{\it Shi arrangement}
$\cA^-$ ([Sh1-3], [St1-2], [AL]), i.e. the arrangement of 
(affine) hyperplanes given by
$$\cA^- = \{ H_\alpha, H_{\alpha-\delta} \ |\ \alpha\in R^+\}
\quad\hbox{where}\quad
\eqalign{
H_\alpha &= \{ x\in \RR^n \ |\ \langle x,\alpha\rangle=0 \}, \cr
H_{\alpha-\delta} &= \{ x\in \RR^n \ |\ \langle x,\alpha\rangle=-1 \}, \cr
}\formula$$
Each chamber $w^{-1}C$, $w\in W$, contains a unique region 
of ${\cal A}^-$ which is a cone, and the vertex of this cone is 
the point $\lambda_w$ which appears in the following theorem.
$$
\matrix{
\beginpicture
\setcoordinatesystem units <1cm,1cm>         
\setplotarea x from -4 to 4, y from -4 to 4    
\put{$H_{\alpha_1}$}[b] at 0 3.1
\put{$H_{\alpha_2}$}[l] at 3.1 3.1
\put{$H_{\alpha_1+\alpha_2}$}[r] at -3.1 3.1
\put{$H_{\alpha_1+2\alpha_2}$}[l] at 3.1 0
\put{$H_{\alpha_1+2\alpha_2-\delta}$}[br] at -3.1 -0.9
\put{$H_{\alpha_1+\alpha_2-\delta}$}[tl] at 1.6 -3.5
\put{$H_{\alpha_2-\delta}$}[tr] at -1.6 -3.5
\put{$H_{\alpha_1-\delta}$}[t] at -0.9 -3.5
\put{$C$} at 1.5 3
\put{$s_1C$} at -1.7 3
\put{$s_2C$} at  2.7 1.5
\put{$s_1s_2C$} at  -3 1.7
\put{$s_2s_1C$} at   3  -1.7 
\put{$s_1s_2s_1C$} at  1.9 -3
\put{$s_2s_1s_2C$} at  -3 -1.7
\put{$\bullet$} at 0 0
\put{$\bullet$} at 2 0
\put{$\bullet$} at 0 -2
\put{$\bullet$} at -2 0
\put{$\bullet$} at -1 1
\put{$\bullet$} at -1 -1
\put{$\bullet$} at 1 -1
\put{$\bullet$} at -1 -3
\put{$\lambda_1$}[bl] at 0.2 0.05
\put{$\lambda_{s_1}$}[l] at -0.9 1.1
\put{$\lambda_{s_2}$}[br] at 1.9  0.1  
\put{$\lambda_{s_1s_2}$}[bl] at -1.9 0.1 
\put{$\lambda_{s_2s_1}$}[t] at 1.1 -1.1 
\put{$\lambda_{s_1s_2s_1}$}[br] at -1.1 -0.9 
\put{$\lambda_{s_2s_1s_2}$}[l] at  0.1 -2 
\put{$\lambda_{s_1s_2s_1s_2}$}[l] at -0.9 -3 
\plot -3 -3   3 3 /
\plot  3 -3  -3 3 /
\plot  0  3   0 -3 /
\plot  3  0  -3  0 /
\setdots
\putrule from -1 3.5 to -1 -3.5
\putrule from 3.5 -1 to -3.5 -1
\plot 3.5 1.5  -1.5 -3.5 /
\plot -3.5 1.5  1.5 -3.5 /
\endpicture
\cr
\cr
\hbox{The arrangement $\cA^-$}
\cr
}$$

\thm [Sb3]
Suppose that $W$ acts irreducibly on $\fh_\RR^*$ and that
$X=\{X^\lambda\ |\ \lambda\in P\}$ where $P$ is the weight lattice.
The algebra $\CC[X]$ is a free $\CC[X]^W$-module with 
$$\hbox{basis\quad $\{ X^{\lambda_w}\ |\ w\in W\}$},
\qquad\hbox{where}\qquad
\lambda_w = w^{-1}\left(\sum_{s_iw<w} \omega_i\right).$$
\pf
The proof is accomplished by establishing three facts:
\smallskip\noindent
\item{(a)}  Let $f_y$, $y\in W$, be a family of elements of
$\ZZ[X]$.  Then
$\det(zf_y)$ is divisible by 
$\displaystyle{ \prod_{\alpha\in R^+} (X^{\alpha}-1)^{|W|/2}. }$
\item{(b)}  $\displaystyle{
\det\big( zX^{\lambda_y} \big)_{z,y\in W} 
= \prod_{\alpha>0} (1-X^{\alpha})^{|W|/2}. }$
\smallskip\noindent
\item{(c)}  If $f\in \ZZ[X]$ then there is a unique solution
to the equation
$$\sum_{w\in W} a_w X^{\lambda_w} = f,
\qquad\hbox{with}\qquad a_w\in \ZZ[X]^W.$$
\smallskip\noindent
(a)  For each $\alpha\in R^+$
subtract row $zf_y$ from row $s_\alpha z f_y$.
Then this row is divisible by $(1-X^{-\alpha})$.
Since there are $|W|/2$ pairs of rows $(zf_y,s_\alpha zf_y)$
the whole determinant is divisible by $(1-X^{-\alpha})^{|W|/2}$.
For $\alpha,\beta\in R^+$ the factors $(1-X^{-\alpha})$ and 
$(1-X^{-\beta})$ are coprime, and so $\det(zf_y)$ is 
divisible by $\displaystyle{
\prod_{\alpha\in R^+} (1-X^{-\alpha})^{|W|/2}.}$
This product and the product in the statement of (a) differ 
by the unit $(X^{2\rho})^{|W|/2}$ in $\ZZ[X]$.
\smallskip\noindent
(b)  By (a), $\det(z X^{\lambda_y})$ is divisible by 
$\displaystyle{
\prod_{\alpha\in R^+} (X^{\alpha}-1)^{|W|/2}. }$
The top coefficient of $\det(zX^{\lambda_y})$ is equal to
$$\prod_{z\in W} zX^{\lambda_z}
=\prod_{z\in W} \prod_{i\atop s_iz<z} X^{\omega_i}
=\prod_{i=1}^n X^{(|W|/2)\omega_i}
=(X^{\rho})^{|W|/2},$$
and the top coefficient of 
$\displaystyle{
\prod_{\alpha\in R^+}(X^{\alpha}-1)^{|W|/2} }$
is $ (X^{2\rho})^{|W|/2}$.
\smallskip\noindent
(c)  Assume that $a_y\in \ZZ[X]^W$ are solutions of the equation
$\sum_{y\in W} X^{\lambda_y} a_y = f$.
Act on this equation by the elements of $W$
to obtain the system of $|W|$ equations
$$\sum_{y\in W} (z X^{\lambda_y})a_y = z f,
\qquad z\in W.$$
By (a) the matrix $(zX^{\lambda_y})_{z,y\in W}$ is invertible
and so this system has a unique solution with $a_y\in \ZZ[X]^W$.
In fact, the $a_y$ can be obtained by Cramer's rule.
Cramer's rule provides an expression for $a_y$ as 
a quotient of two determinants.
By (a) and (b) the denominator divides the numerator
to give an element of $\ZZ[X]$.
Since each determinant is an alternating function, the quotient
is an element of $\ZZ[X]^W$.
\endpf

\smallskip\noindent
{\bf Remark.}  In [Sb2] Steinberg proves this type of result in full
generality without the assumptions that $W$ acts irreducibly
on $\fh_{\RR}^*$ and $L=P$.  Note also that the proof given above is 
sketchy, particularly in the aspect that the top coefficient 
of the determinant is what we have claimed it is.  See
[Sb2] for a proper treatment of this point.

\subsection Deducing the $\tilde H_L$ representation theory from $\tilde H_P$.

It is often easier to work with the representation theory
of $\tilde H$ in the case when $L=P$.  It is important to
be able to convert from this case to the case of a general
lattice $L$.  If $W$ acts irreducibly on $\fh_{\RR}^*$ then
the lattice $L$ satisfies
$$Q\subseteq L\subseteq P,
\qquad
\hbox{where}\qquad
P = \sum_{i=1} \ZZ\omega_i
\quad\hbox{and}\quad
Q = \sum_{i=1} \ZZ\alpha_i
$$
are the weight lattice and the {\it root lattice} respectively.
The group $\Omega = P/Q$ is a finite group 
(either cyclic or isomorphic to $\ZZ/2\ZZ\times \ZZ/2\ZZ$). It
corresponds to the center of the corresponding complex algebraic
group.  Let us denote the corresponding affine Hecke algebras by
$$\tilde H_Q \subseteq \tilde H_L \subseteq \tilde H_P,$$
according which lattice is used to make the group $X$.

\thm [RR]
Then there is an action of the finite group $P/L$ on $\tilde H_P$,
by ring automorphisms, such that
$$\tilde H_L = (\tilde H_P)^{P/L} = \{ h\in \tilde H_P\ |\ 
gh=h \hbox{\ for all $g\in P/L$}\},$$
is the subalgebra of fixed points
under the action
of the group $P/L$.
\endthm

This theorem is exactly what is needed to apply a (not very well known) 
version of Clifford theory to completely classify the 
representations of $\tilde H_L$ in terms of the representations 
of $\tilde H_P$, see [RR].


\section 2. $\tilde H$-modules

\subsection{Weights.}

In view of the results in (1.18) we shall 
(for the remainder of this paper, except sections 5, 6 and 7 
where we use the data in (1.1)) assume
that $L=P$ in the definition of the group $X$ and $\tilde H$, 
see (1.2), (1.9) and (1.14).
The Weyl group acts on 
$$T={\rm Hom}(X,\CC^*)
=\{ \hbox{group homomorphisms $t\colon X \to \CC^*$}\},
\qquad\hbox{by}\quad
(wt)(X^\lambda) = t(X^{w^{-1}\lambda}).$$
Let $M$ be a finite dimensional $\tilde H$-module and let $t\in T$.
The {\it $t$-weight space} 
and the {\it generalized $t$-weight space} of $M$ are 
$$\eqalign{
M_t &= \{ m\in M \ |\ X^\lambda m = t(X^\lambda)m
\hbox{\ for all $X^\lambda\in X$}\},
\qquad\hbox{and}  \cr
\cr
M_t^{\rm gen} &= \{ m\in M\ |\
\hbox{for each $X^\lambda\in X$,
$(X^\lambda-t(X^\lambda))^k m =0$ for some $k\in \ZZ_{>0}$}
\},  \cr
}
$$
respectively.
Then
$$
M=\bigoplus_{t\in T} M_t^{\rm gen}.
\formula$$
is a decomposition of $M$ into Jordan blocks
for the action of $\CC[X]$,  
and we say that $t$ is a {\it weight} of $M$ if $M_t^{\rm gen}\ne 0$.
Note that $M_t^{\rm gen}\ne 0$ if and only if $M_t\ne 0$.
A finite dimensional $\tilde H$-module
$$\hbox{$M$ is {\it calibrated} if $M_t^{\rm gen}=M_t$,
for all $t\in T$.}$$

\medskip\noindent
{\bf Remark.}  The term tame is sometimes used in place
of the term calibrated particularly in the context of
representations of of Yangians, see [NT].  The word
calibrated is preferable since tame also has many other meanings
in different parts of mathematics.

\medskip
Let $M$ be a simple $\tilde H$-module.  As an $X(T)$-module,
$M$ contains a simple submodule and this submodule must be
one-dimensional since all irreducible representations of a commutative
algebra are one-dimensional.  Thus, 
a simple module always has $M_t\ne 0$ for some $t\in T$.

\subsection{Central characters.}

The Pittie-Steinberg theorem, Theorem 1.17, shows that, as vector spaces,
$$\tilde H = H\otimes \CC[X] 
= H\otimes \CC[X]^W\otimes {\cal K},
\qquad\hbox{where}\qquad
{\cal K} = \CC\hbox{-span}\{X^{\lambda_w}\ |\ w\in W\},
$$
and $H$ is the Iwahori-Hecke algebra defined in (1.8).
Thus $\tilde H$ is a free module
over $Z(\tilde H)=\CC[X]^W$ of rank $~\dim(H)\cdot\dim({\cal K})=|W|^2$.
By Dixmier's version of Schur's lemma (see [Wa, Lemma 0.5.1]),
$Z(\tilde H)$ acts on a simple $\tilde H$-module by
scalars and so it follows that
every simple $\tilde H$-module is finite dimensional
of dimension $\le |W|^2$.  Theorem 2.12(d) below will show that,
in fact, the dimension of a simple module is $\le |W|$.

Let $M$ be a simple $\tilde H$-module.
The {\it central character} of $M$ is an element $t\in T$ such that 
$$pm = t(p)m,
\qquad\hbox{for all $m\in $M, $p\in \CC[X]^W=Z(\tilde H)$.}$$
The element $t$ is only determined up to the action of
$W$ since $t(p)=wt(p)$ for all $w\in W$.  Because of this,
any element of the orbit $Wt$ is referred to as the {\it central
character} of $M$.

Because $P=L$ in the construction of $X$, a theorem of 
Steinberg [Sb2, 3.15, 4.2, 5.3] tells us that 
the stabilizer $W_t$ of a point $t\in T$ under the action of $W$ is
the reflection group
$$W_t =  \langle s_\alpha \ | \ \alpha\in Z(t) \rangle,
\qquad\hbox{where}\qquad
Z(t) = \{ \alpha\in R^+\ |\ t(X^{\alpha})=1\}.$$
Thus the orbit $Wt$ can be viewed in several different 
ways via the bijections
$$Wt \longleftrightarrow
W/W_t \longleftrightarrow
\{ w\in W\ |\ R(w)\cap Z(t)=\emptyset\}
\longleftrightarrow
\left\{ \matrix{ \hbox{chambers on the positive}\cr
\hbox{side of $H_\alpha$ for $\alpha\in Z(t)$} \cr}
\right\},
\formula$$
where the last bijection is the restriction of the map
in (1.6).  If the root system $Z(t)$ is generated by the
simple roots $\alpha_i$ that it contains then $W_t$ is 
a parabolic subgroup of $W$ and $\{w\in W\ |\ R(w)\cap Z(t)\}$
is the set of minimal length coset representatives of the
cosets in $W/W_t$.

\subsection{Principal series modules.}

For $t\in T$ let $\CC v_t$ be the one-dimensional
$\CC[X]$-module given by
$$X^\lambda v_t =t(X^\lambda)v_t,
\qquad\hbox{for $X^\lambda\in X$.}
$$  
The {\it principal series representation} $M(t)$ is the
$\tilde H$-module defined by
$$
M(t) = \tilde H \otimes_{\CC[X]} \CC v_t
= {\rm Ind}_{\CC[X]}^{\tilde H}(\CC v_t).
\formula$$
The module $M(t)$ has basis $\{ T_w\otimes v_t \ |\ w\in W  \}$ 
with $H$ acting by left multiplication. 

If $w\in W$ and $X^\lambda\in X$ then the defining relation (1.10) for
$\tilde H$ implies that
$$X^\lambda (T_w\otimes v_t)= t(X^{w\lambda}) (T_w\otimes v_t)
+\sum_{u<w} a_u (T_u\otimes v_t),
\formula$$
where the sum is over $u<w$ in the Bruhat-Chevalley order and $a_u\in \CC$. 
Let $W_t = {\rm Stab}(t)$
be the stabilizer of $t$ under the $W$-action.
It follows from (2.7)
that the eigenvalues of $X$ on $M(t)$ are of the form 
$wt$, $w\in W$,
and by counting the multiplicity of each eigenvalue we have
$$M(t) = \bigoplus_{wt\in Wt} M(t)_{wt}^{\rm gen}
\qquad\hbox{where}\quad\dim(M(t)_{wt}^{\rm gen}) = |W_t|,
\quad\hbox{for all $w\in W$.}
\formula$$
In particular, if $t$ is regular (i.e. when $W_t$ is trivial), there
is a unique basis $\{v_{wt}\ |\ w\in W\}$ of $M(t)$ determined
by
$$\eqalign{
X^\lambda v_{wt} &= (wt)(X^\lambda) v_{wt},
\qquad\hbox{for all $w\in W$ and $\lambda\in P$,} \cr
v_{wt} &= T_w\otimes v_t +\sum_{u<w} a_{wu}(t)(T_u\otimes v_t),
\qquad\hbox{where $a_{wu}(t)\in \CC$.} \cr
}\formula$$

Let $t\in T$.  The {\it spherical vector} in $M(t)$ is
$${\bf 1}_t = \sum_{w\in W} q^{\ell(w)}T_w\otimes v_t.
\formula$$
Up to multiplication by constants this is the unique
vector in $M(t)$ such that 
$T_w{\bf 1}_t=q^{\ell(w)}{\bf 1}_t$ for all $w\in W$.  
The following is due to Kato, [Kt1, Proposition 1.20 and Lemma 2.3].

\prop Let $t\in T$ and let $W_t$ be the stabilizer of $t$ under
the $W$-action. 
\smallskip\noindent
\item{(a)}  If $W_t=\{1\}$ and $v_{wt}$,
$w\in W$ is the basis of $M(t)$ defined in (2.9) then
$${\bf 1}_t = \sum_{z\in W} t(c_z),
\qquad\hbox{where}\qquad
c_z = \prod_{\alpha\in R(w_0z)} {q-q^{-1}X^\alpha\over 1-X^\alpha}.$$
\item{(b)} The spherical vector ${\bf 1}_t$ generates
$M(t)$\quad if and only if \quad
$\displaystyle{
t\left(\prod_{\alpha\in R^+} (q^{-1}-qX^\alpha)\right) \ne 0.}$
\item{(c)} The module $M(t)$ is irreducible
if and only if ${\bf 1}_{wt}$ generates $M(wt)$ for all $w\in W$. 
\pf
The proof is accomplished in exactly the same way as done for 
the graded Hecke algebra in [KR, Proposition 2.8].  The
only changes which need to be made to [KR] are 
\smallskip\noindent
\item{(1)}
Use $\displaystyle{
T_i \left(\sum_{w\in W} q^{\ell(w)} T_w\right)
= q\left(\sum_{w\in W} q^{\ell(w)} T_w\right)}$
and
$\displaystyle{
{\bf 1}_t = 
\left(\sum_{w\in W} q^{\ell(w)} T_w\right)v_t}$ and the
$\tau$-operators defined in Proposition 2.14 for the proof of (a).
(We have included this result in this section since it is really 
a result about the structure of principal series modules.
Though the proof uses the $\tau$-operators, which we will define
in the next section, there is no logical gap here.)
\smallskip\noindent
\item{(2)}  For the proof of (b) use the Steinberg basis 
$\{ X^{\lambda_y}\ |\ y\in W\}$
and the determinant $\det(X^{z^{-1}\lambda_y})$
from Theorem 1.17(b) in place of the basis $\{b_y\ |\ w\in W\}$ and the 
determinant used in [KR].
\endpf

Part (b) of the following theorem is due to Rogawski
[Rg, Proposition 2.3] and part (c) is due to Kato [Kt1, Theorem 2.1].  
Parts (a) and (d) are classical.

\thm Let $t\in T$ and $w\in W$ and define 
$P(t) = \{ \alpha\in R^+\ |\ t(X^{\alpha})=q^{\pm2}\}$.
\item{(a)}  If $W_t=\{1\}$ then $M(t)$ is calibrated.
\item{(b)}  $M(t)$ and $M(wt)$ have the same composition factors.
\item{(c)}  $M(t)$ is irreducible if and only if $P(t)=\emptyset$.
\item{(d)}  If $M$ is a simple $\tilde H$-module with $M_t\ne 0$
then $M$ is a quotient of $M(t)$.
\pf
(a) follows from (2.8) and the definition of calibrated.
Part (b) accomplished exactly as in [KR, Proposition 2.8] and
(c) is a direct consequence of the definition of $P(t)$ and
Proposition 2.11.
\smallskip\noindent
(d)  Let $m_t$ be a nonzero vector in $M_t$.  
If $v_t$ is as in the construction of $M(t)$ in (2.6)
then, as $\CC[X]$-modules, $\CC m_t\cong \CC v_t$. Thus, since
induction is the adjoint functor to restriction there is
a unique $\tilde H$-module homomorphism given by
$$
\matrix{
\phi \colon &M(t) &\longrightarrow &M \cr
&v_t &\longmapsto &m_t \cr
}
$$
This map is surjective since $M$ is irreducible and so
$M$ is a quotient of $M(t)$.
\endpf

\subsection{The $\tau$ operators.}

The following proposition defines maps 
$\tau_i\colon M_t^{\rm gen}\to M_{s_it}^{\rm gen}$ 
on generalized weight spaces of finite dimensional $\tilde H$-modules
$M$.  These are ``local operators'' and are only defined on 
weight spaces $M_t^{\rm gen}$ such that $t(X^{\alpha_i})\ne 1$.
In general, $\tau_i$ does not extend to an operator on all of $M$.

\prop
Fix $i$, let $t\in T$ be such that $t(X^{\alpha_i})\ne 1$ and 
let $M$ be a finite dimensional $\tilde H$-module.  Define
$$
\matrix{
\tau_i \colon 
&M_t^{\rm gen} &\longrightarrow &M_{s_it}^{\rm gen} \cr
\cr
& m & \longmapsto &
\displaystyle{
\left(T_i - {q-q^{-1}\over 1-X^{-\alpha_i} }\right) m.} \cr
}
$$
\item{(a)}
The map $\tau_i\colon M_t^{\rm gen} \longrightarrow M_{s_it}^{\rm gen}$
is well defined.
\smallskip
\item{(b)}  As operators on $M_t^{\rm gen}$, \quad
$\displaystyle{ X^\lambda \tau_i = \tau_i X^{s_i\lambda}}$,
for all $X^\lambda\in X$.
\smallskip
\item{(c)}  As operators on $M_t^{\rm gen}$, \quad
$\displaystyle{
\tau_i\tau_i=
{(q-q^{-1}X^{\alpha_i})(q-q^{-1}X^{-\alpha_i})\over 
(1-X^{\alpha_i})(1-X^{-\alpha_i}) }.
}
$
\smallskip
\item{(d)}  Both maps $\tau_i\colon M_t^{\rm gen}\to M_{s_it}^{\rm gen}$
and $\tau_i\colon M_{s_it}^{\rm gen}\to M_t^{\rm gen}$
are invertible if and only if $t(X^{\alpha_i})\ne q^{\pm 2}$.
\smallskip
\item{(e)}  Let $1\le i\ne j\le n$ and let $m_{ij}$ be as in (1.7).
Then   
$$\underbrace{\tau_i\tau_j\tau_i\cdots}_{m_{ij} {\rm\ factors}}
= \underbrace{\tau_j\tau_i\tau_i\cdots}_{m_{ij} {\rm\ factors}},
$$
whenever both sides are well defined operators on $M_t^{\rm gen}$.
\pf
(a) 
The element $X^{\alpha_i}$ acts on $M_t^{\rm gen}$ by
$t(X^{\alpha_i})$ times a unipotent transformation.
As an operator on $M_t^{\rm gen}$, 
$1-X^{-\alpha_i}$ is invertible since it has determinant
$(1-t(X^{-\alpha_i}) )^d$ where $d=\dim(M_t^{\rm gen})$.
Since this determinant is nonzero $(q-q^{-1})/(1-X^{-\alpha_i})
=(q-q^{-1})(1-X^{-\alpha_i})^{-1}$ 
is a well defined operator on $M_t^{\rm gen}$.  
Thus the definition of $\tau_i$ makes sense.

Since $(q-q^{-1})/(1-X^{-\alpha_i})$ is not an
element of $\tilde H$ or $\CC[X]$ it should be viewed only
as an operator on $M_t^{\rm gen}$ in calculations.
With this in mind it is straightforward to use the defining
relation (1.10) to check that 
$$
\eqalign{
X^\lambda \tau_i m
&=X^\lambda\left(T_i - {q-q^{-1}\over 1-X^{-\alpha_i} } \right)m
= \left(
T_i - {q-q^{-1}\over 1-X^{-\alpha_i} } \right)
X^{s_i\lambda}m 
= \tau_i X^{s_i\lambda}m, \qquad\hbox{and} \cr
\tau_i\tau_i m
&=
\left(T_i - {q-q^{-1}\over 1-X^{-\alpha_i} }\right)
\left(T_i - {q-q^{-1}\over 1-X^{-\alpha_i} }\right) m 
=
{ (q-q^{-1}X^{\alpha_i})(q-q^{-1}X^{-\alpha_i})\over 
(1-X^{\alpha_i})(1-X^{-\alpha_i}) } 
 m, \cr
}
$$
for all $m\in M_t^{\rm gen}$ and $X^\lambda\in X$.
This proves (a), (b) and (c).
\smallskip\noindent
(d)  The operator $X^{\alpha_i}$ acts on 
$M_t^{\rm gen}$ as $t(X^{\alpha_i})$ times a unipotent transformation.
Similarly for $X^{-\alpha_i}$.  Thus, as an operator on $M_t^{\rm gen}$
$\det((q-q^{-1}X^{\alpha_i})(q-q^{-1}X^{-\alpha_i}))=0$
if and only if $t(X^{\alpha_i})=q^{\pm 2}$.
Thus part (c) implies that $\tau_i\tau_i$, and each factor
in this composition,
is invertible if and only if $t(X^{\alpha_i})\ne q^{\pm 2}$.
\smallskip\noindent
(e) 
Let $t\in T$ be regular.
By part (a), the definition of the $\tau_i$, and the
uniqueness in (2.9), the basis $\{v_{wt}\}_{w\in W}$
of $M(t)$ in (2.9) is given by
$$v_{wt} = \tau_w v_t,
\formula
$$
where $\tau_w=\tau_{i_1}\cdots \tau_{i_p}$ for a reduced word
$w=s_{i_1}\cdots s_{i_p}$ of $w$.  
Use the 
defining relation (1.10) for $\tilde H$ to expand the product of
$\tau_i$ and compute
$$\eqalign{
v_{w_0t} = \underbrace{\cdots \tau_i\tau_j\tau_i}_{m_{ij} {\rm\ factors}}v_{t} 
&= \underbrace{\cdots T_iT_jT_i}_{m_{ij} {\rm\ factors}} v_{t} +
\sum_{w<w_0} T_w P_w v_{t} 
=T_{w_0} \otimes v_{t} +
\sum_{w<w_0} t(P_w) T_w\otimes v_{t}  \cr
\phantom{v_{w_0t}} 
= \underbrace{\cdots \tau_j\tau_i\tau_j}_{m_{ij} {\rm\ factors}}v_{t} 
&= \underbrace{\cdots T_jT_iT_j}_{m_{ij} {\rm\ factors}} v_{t} +
\sum_{w<w_0} T_w Q_w v_{t} 
=T_{w_0} \otimes v_{t} +
\sum_{w<w_0} t(Q_w) T_w\otimes v_{t}  \cr
}$$
where $P_w$ and $Q_w$ are rational functions in the $X^\lambda$.
By the uniqueness in (2.9), $t(P_w)=a_{w_0 w}(t)=t(Q_w)$ for all $w\in W$, $w\ne w_0$.
Since the values of $P_w$ and $Q_w$ coincide on all generic
points $t\in T$ it follows that 
$$P_w=Q_w
\qquad\hbox{for all $w\in W$, $w\ne w_0$.}
\formula$$
Thus,
$$
\underbrace{\cdots \tau_i\tau_j\tau_i}_{m_{ij} {\rm\ factors}}
= T_{w_0} + \sum_{w<w_0} T_w P_w  
= T_{w_0} + \sum_{w<w_0} T_w Q_w 
= \underbrace{\cdots \tau_j\tau_i\tau_j}_{m_{ij} {\rm\ factors}}, 
$$
whenever both sides are well defined operators on $M_t^{\rm gen}$.
\endpf

Let $t\in T$ and recall that 
$$
Z(t) = \{ \alpha\in R^+\ |\ t(X^{\alpha})=1\}
\qquad\hbox{and}\qquad
P(t) = \{ \alpha\in R^+\ |\ t(X^{\alpha})=q^{\pm2}\}.
\formula$$
If $J\subseteq P(t)$ define
$${\cal F}^{(t,J)} = 
\{w\in W\ |\ R(w)\cap Z(t) = \emptyset,\ \ 
R(w)\cap P(t)=J\}.\formula$$
We say that the pair $(t,J)$ is a {\it local region} if
${\cal F}^{(t,J)}\ne \emptyset$.
Under the bijection (2.4) the set ${\cal F}^{(t,J)}$ maps to
the set of chambers whose union is
the set of points $x\in \fh_{\RR}^*$ which are
\smallskip\noindent
\itemitem{(a)} on the positive side of the hyperplanes
$H_\alpha$ for $\alpha\in Z(t)$,
\itemitem{(b)} on the positive side of the hyperplanes
$H_\alpha$ for $\alpha\in P(t)\backslash J$,
\itemitem{(c)} on the negative side of the hyperplanes
$H_\alpha$ for $\alpha\in J$.
\smallskip\noindent
See the picture in Example 4.11(d).
In this way the local region $(t,J)$ really does correspond 
to a region in $\fh_{\RR}^*$.  This is a connected convex region
in $\fh_{\RR}^*$ since it is cut out by half spaces in 
$\fh_{\RR}^*\cong \RR^n$.  The elements $w\in {\cal F}^{(t,J)}$
index the {\it chambers $w^{-1}C$ in the local region}
and, as $J$ runs over the subsets of $P(t)$,
 the sets ${\cal F}^{(t,J)}$ form a partition of the
set $\{w\in W\ |\ R(w)\cap Z(t)=\emptyset\}$ (which, by (2.4),
indexes the cosets in $W/W_t$).

\cor Let $M$ be a finite dimensional $\tilde H$-module.
Let $t\in T$ and let $J\subseteq P(t)$.  Then
$$\dim(M_{wt}^{\rm gen}) = \dim(M_{w't}^{\rm gen}),
\qquad\hbox{for $w,w'\in {\cal F}^{(t,J)}$.}
$$
\pf
Suppose $w, s_iw\in {\cal F}^{(t,J)}$.  We may assume that
$s_iw>w$.  Then $\alpha=w^{-1}\alpha_i>0$, $\alpha\not\in R(w)$
and $\alpha\in R(s_iw)$.  Now, $R(w)\cap Z(t) = R(s_iw)\cap Z(t)$
implies $t(X^{\alpha})\ne 1$, and
$R(w)\cap P(t)$ implies $t(X^{\alpha})\ne q^{\pm 2}$.  Since
$wt(X^{\alpha_i})=t(X^{w^{-1}\alpha_i})=t(X^\alpha)\ne 1$
and $wt(X^{\alpha_i})\ne q^{\pm2}$ and thus, by Proposition
2.14(d), the map $\tau_i\colon M_{wt}^{\rm gen}\to M_{s_iwt}^{\rm gen}$
is well defined and invertible.  It remains to note that 
if $w,w'\in {\cal F}^{(t,J)}$, then $w'=s_{i_1}\cdots s_{i_\ell}w$
where $s_{i_k}\cdots s_{i_\ell}w\in {\cal F}^{(t,J)}$
for all $1\le k\le \ell$.  This follows from the fact that
${\cal F}^{(t,J)}$ corresponds to a connected convex region in 
$\fh_{\RR}^*$.
\endpf

\bigskip
\section 3.  Classification of calibrated representations

For simple roots $\alpha_i$ and $\alpha_j$ in $R$ and let
$R_{ij}$ be the rank two root subsystem of $R$ generated
by $\alpha_i$ and $\alpha_j$.  
A weight $t\in T$ is {\it calibratable} if, for every pair
$i,j$, $i\ne j$, $t$ is a weight of a calibrated representation 
of the rank two affine Hecke (sub)algebra generated by
$T_i, T_j$ and $\CC[X]$.  
A local region 
$$\hbox{$(t,J)$ is {\it skew}
\quad if\quad
$wt$ is calibratable
for all $w\in {\cal F}^{(t,J)}$.}
$$
The classification of irreducible representations of rank two 
affine Hecke algebras
given in [R1]  can be used to state this condition
combinatorially.  Specifically,
a weight $t\in T$ is {\it calibratable} if
\smallskip
\itemitem{(a)} For all simple roots $\alpha_i$,
$1\le i\le n$,\ \ $t(X^{\alpha_i})\ne 1$, and
\smallskip
\itemitem{(b)} For all pairs of simple roots $\alpha_i$ and $\alpha_j$
such that $\{\alpha\in R_{ij}\ |\ t(X^\alpha)=1\}\ne \emptyset$,
\itemitem{} the set $\{\alpha\in R_{ij}\ |\ t(X^\alpha)=q^{\pm2}\}$
contains more than two elements.
\smallskip\noindent
Condition (a) says that $t$ is regular for all rank 1 subsystems
of $R$ generated by simple roots.  This condition guarantees that 
the weight is ``calibratable'' (i.e. appears as a weight of
some calibrated representation) for all rank 1 affine Hecke subalgebras
of $\tilde H$.  Condition (b) is an ``almost regular'' condition
on $t$ with respect to rank 2 subsystems generated by simple roots.

\medskip\noindent
{\bf Remark.}  The conversion between the definition of
calibratable weight and the combinatorial condition given in (a) and (b)
is as follows.  Consider a rank two affine Hecke algebra $\tilde H$.
\item{(A)}  By Theorem 2.12(a) and Theorem 2.12(d), local regions 
$(t,J)$ with $t$ regular satisfy (a) and (b)
and always contribute calibrated representations of $\tilde H$.
\item{(B)}  Using the notation of [R1],
the local regions $(t,J)$ with $t$ nonregular and 
which satisfy both conditions (a) and (b) are 
\itemitem{} Type $A_2$: none,
\itemitem{} Type $C_2$: $(t_b,\{\alpha_1\})$ and 
$(t_b,\{\alpha_1,\alpha_1+\alpha_2\})$ (for each of these
$P(t)$ contains $3$ elements),
\itemitem{} Type $G_2$:  $(t_e,J)$ with $J\ne \emptyset$ and $J\ne P(t_e)$
(for each of these $P(t_e)$ contains $4$ elements).
\smallskip\noindent
From (A) and (B) it follows that the local regions which satisfy
(a) and (b) do contribute calibrated weights.  The following 
shows that the other local regions don't contribute calibratable weights.
\item{(C)}  By Lemma 3.1(a) local regions $(t,J)$ with a weight
$\xi=wt$, $w\in {\cal F}^{(t,J)}$ such that $\xi(X^{\alpha_i})=1$
don't satisfy (a) and, by inspection of the tables in [R1], they never
contribute a calibrated representation.
\item{(D)}  Using the notation of [R1],
the local regions which satisfy condition (a) but not
condition (b) are 
\itemitem{} Type $A_2$: $(t_c, \{\alpha_2\})$ and $(t_d,\{\alpha_1\})$
\itemitem{} Type $B_2$: $(t_d,\{\alpha_2\})$
\itemitem{} Type $G_2$: $(t_i, \{\alpha_2\})$, $(t_f,\{\alpha_1\})$,
\smallskip\noindent
(note that to satisfy (b) $Z(t)$ must be nonempty). 
From the tables in [R1] we see that none of these local regions
supprts a calibrated representation.

\medskip\noindent
{\bf Remark.}  The paper [R1] does not treat roots of unity.
However, it is interesting to note that, 
{\it provided} $q^2\ne\pm1$,  the methods of [R1]
go through without change to classify all representations
of rank two affine Hecke algebras even when $q^2$ is a root of unity.
This classification can be used (as in the previous remark)
to show that (a) and (b) above still characterize calibratable 
weights when $q^2$ is a root of unity such that $q^2\ne \pm 1$.
The key point is that Lemma 1.19 of [R1] still holds.  
If $q^2=-1$ then Lemma 1.19
of [R1] breaks down at the next to last line of the proof
in the statement
``$\ldots$ forces $\phi(wt(T_j)$ to have Jordan blocks of size $1$
$\ldots$''.  When $q^2={-1}$ it is possible that $\phi(wt(T_j))$
has a Jordan block of size 2.  If $q^2=1$ then one can change the
definition of the $\tau$-operators and use similar methods
to produce a complete analysis of simple $\tilde H$-modules, but
we shall not do this here, choosing instead to 
exclude the case $q^2=1$ for simplicity of exposition.

\medskip
\medskip
The following lemma provides fundamental results
about the structure of irreducible calibrated $\tilde H$-modules.
We omit the proof since it is accomplished 
in exactly the same way as in [KR, Lemma 4.1 and Lemma 4.2].

\lemma Let $M$ be an irreducible calibrated module. 
Then, for all $t\in T$ such that $M_t\ne 0$,
\smallskip
\item{(a)} If $t\in T$ such that $M_t\ne 0$ then
$t(X^{\alpha_i})\ne 1$ for all $1\le i\le n$.
\smallskip
\item{(b)} If $t\in T$ such that $M_t\ne 0$ then $\dim(M_t)=1$.
\smallskip\noindent
\item{(c)} If $t\in T$ such that both
$M_t$ and $M_{s_it}$ are both nonzero then the
map $\tau_i\colon M_t\to M_{s_it}$ is a bijection.
\endlemma

This lemma together with the classification of irreducible
modules for rank two affine Hecke algebras gives the following 
fundamental structural result for irreducible calibrated $\tilde H$-modules.
The proof is essentially the same as the proof of Proposition 4.3 in [KR].  
We repeat the proof here for continuity.

\thm  If $M$ is an irreducible calibrated $\tilde H$-module
with central character $t\in T$ then there is a unique
skew local region $(t,J)$ such that 
$$\dim(M_{wt}) = \cases{ 1, &for all $w\in {\cal F}^{(t,J)}$, \cr
0 &otherwise. \hfill \cr
}$$
\pf
By Lemma 3.1b all nonzero generalized weight spaces of $M$
have dimension $1$ and by Lemma 3.1c all $\tau$-operators
between these weight spaces are bijections.  This already guarantees
that there is a unique local region $(t,J)$ which satisfies
the condition.  It only remains to show that this local region
is skew.  

Let $\tilde H_{ij}$ be the subalgebra generated by $T_i$, $T_j$
and $\CC[X]$.  Since $M$ is calibrated as an $\tilde H$-module
it is calibrated as a $\tilde H_{ij}$-module and so all factors
of a composition series of $M$ as an $\tilde H_{ij}$-module are 
calibrated.  Thus the weights of $M$ are calibratable.
So $(t,J)$ is a skew local region.
\endpf

The following Proposition shows that the weight space structure
of calibrated representations, as determined in Theorem 3.2,
essentially forces the $\tilde H$-action on a weight
basis.
The proof is quite similar to the proof of Proposition 4.4 in [KR].  
However, we include the details since there is a technicality here; to 
make the conclusion in (3.4) we use the fact that the group $X$ 
corresponds to the weight lattice $L=P$.

\prop  Let $M$ be a calibrated
$\tilde H$-module and assume that for all $t\in T$ such that $M_t\ne 0$,
$$\hbox{(A1)\quad $t(X^\alpha_i)\ne 1$ for all $1\le i\le n$,
\qquad and
\qquad (A2)\quad $\dim(M_t)=1$.}
$$ 
For each $b\in T$ such that $M_b\ne 0$ let $v_b$ be a nonzero vector in
$M_b$.  The vectors $\{v_b\}$ form a basis of $M$.
Let $(T_i)_{cb}\in \CC$ and $b(X^\lambda)\in \CC$ be given by
$$T_iv_b = \sum_c (T_i)_{cb}v_c  
\qquad\hbox{and}\qquad 
X^\lambda v_b = b(X^\lambda) v_b.$$
Then
\smallskip
\itemitem{(a)} $\displaystyle{
(T_i)_{bb} = {q-q^{-1}\over 1-b(X^{-\alpha_i}) }
},$
for all $v_b$ in the basis, 
\smallskip
\itemitem{(b)} If $(T_i)_{cb}\ne 0$ then $c=s_ib$,
\smallskip
\itemitem{(c)} $(T_i)_{b,s_ib}(T_i)_{s_ib,b} 
= (q^{-1}+(T_i)_{bb})(q^{-1}+(T_i)_{s_ib,s_ib}).$
\pf
The defining equation for $\tilde H$,
$$
X^\lambda T_i - T_i X^{s_i\lambda} =
(q-q^{-1}) { X^\lambda-X^{s_i\lambda}\over 1-X^{-\alpha_i} },
$$
forces
$$
\sum_c \left( c(X^\lambda)(T_i)_{cb} - (T_i)_{cb}b(X^{s_i\lambda})
\right)v_c
=
(q-q^{-1}) { b(X^\lambda)-b(X^{s_i\lambda})\over 1-b(X^{-\alpha_i}) }v_b
$$
Comparing coefficients gives
$$\eqalign{
c(X^\lambda)(T_i)_{cb}-(T_i)_{cb} b(X^{s_i\lambda}) 
&= 0,
\qquad\hbox{if $b\ne c$, and } \cr
\cr
b(X^\lambda)(T_i)_{bb}-(T_i)_{bb}b(X^{s_i\lambda})
&=(q-q^{-1}){ b(X^\lambda)-b(X^{s_i\lambda})\over 1-b(X^{-\alpha_i}) }. \cr
}
$$
These relations give:
$$
\hbox{If \quad $(T_i)_{cb}\ne 0$ \quad then \quad $b(X^{s_i\lambda}) =
c(X^\lambda)$
 \quad for all $X^\lambda\in X$, and}
$$
$$
(T_i)_{bb} = {q-q^{-1}\over 1-b(X^{-\alpha_i}) }
\quad \hbox{ if $b(X^{-\alpha_i})\ne 1$ and $b(X^\lambda)\ne b(X^{s_i\lambda})$ for
some $X^\lambda\in X$.}
$$
By assumption (A1), $b(X^{\alpha_i})\ne 1$ for all $i$.
For each fundamental weight $\omega_i$, $X^{\omega_i}\in X$ and
$b(X^{s_i\omega_i})=b(X^{\omega_i-\alpha_i})\ne b(X^{\omega_i})$ since
$b(X^{\alpha_i})\ne 1$.  Thus we conclude that 
$$
T_i v_b = (T_i)_{bb} v_b + (T_i)_{s_ib,b}v_{s_ib},
\qquad\hbox{with}\qquad
(T_i)_{bb} = {q-q^{-1}\over 1-b(X^{-\alpha_i}) }.
\formula$$
This completes the proof of (a) and (b).
By the definition of $\tilde H$ the vector
$$
T_i^2 v_b = ((T_i)_{bb}^2 +(T_i)_{b,s_ib}(T_i)_{s_ib,b})v_b
+((T_i)_{bb}+(T_i)_{s_ib,s_ib})(T_i)_{s_ib,b}v_{s_ib}
$$
must equal
$$
((q-q^{-1})T_i+1)v_b 
=((q-q^{-1})(T_i)_{bb} +1)v_b + (q-q^{-1})(T_i)_{s_ib,b}v_{s_ib}.
$$
Using the formula for $(T_i)_{bb}$ and $(T_i)_{s_ib,s_ib}$ we find
$(T_i)_{bb}+(T_i)_{s_ib,s_ib} = (q-q^{-1})$.  
So, by comparing coefficients of $v_b$, we obtain the equation
$$
(T_i)_{b,s_ib}(T_i)_{s_ib,b} = (q-(T_i)_{bb})((T_i)_{bb}+q^{-1})
= (q^{-1}+(T_i)_{bb})(q^{-1}+(T_i)_{s_ib,s_ib}).\qquad\hbox{\qed}
$$

\thm  Let $(t,J)$ be a skew local region and let $\cF^{(t,J)}$ index
the chambers in the local region $(t,J)$.
Define 
$$\tilde H^{(t,J)} =
\CC\hbox{-span} \{ v_w \ |\ w\in {\cF}^{(t,J)} \},$$
so that the symbols $v_w$ are a labeled basis of the
vector space $\tilde H^{(t,J)}$.
Then the following formulas make $\tilde H^{(t,J)}$
into an irreducible $\tilde H$-module:  For each $w\in {\cF}^{(t,J)}$,
$$
\matrix{
\hfill X^\lambda v_w &=& (wt)(X^\lambda) v_w, \hfill
&&\hbox{for $X^\lambda\in X$, and}\hfil \cr
\cr
\hfill T_i v_w 
&=& (T_i)_{ww} v_w + (q^{-1}+(T_i)_{ww}) v_{s_iw}, \hfill
&\qquad &\hbox{for $1\le i\le n$,}\hfil \cr
}
$$
where $\displaystyle{
(T_i)_{ww} = {q-q^{-1}\over 1 - (wt)(X^{-\alpha_i})}\;, }
$
and we set $v_{s_iw} = 0$ if $s_iw\not\in {\cF}^{(t,J)}$.
\pf
Since $(t,J)$ is a skew local region 
$(wt)(X^{-\alpha_i})\ne 1$ for all $w\in
\cF^{(t,J)}$ and all simple roots $\alpha_i$.  This implies that the 
coefficient $(T_i)_{ww}$ is well defined for all $i$ and $w\in \cF^{(t,J)}$.

By construction, the nonzero weight spaces of $\tilde
H^{(t,J)}$ are 
$(\tilde H^{(t,J)})_{wt}^{\rm gen}=(\tilde H^{(t,J)})_{wt}$ 
where $w\in \cF^{(t,J)}$.  
Since $\dim(\tilde H^{(t,J)})=1$ for $u\in {\cal F}^{(t,J)}$,
any proper submodule $N$ of $\tilde H^{(t,J)}$ must have
$N_{wt}\ne 0$ and $N_{w't}=0$ for some $w\ne w'$ with
$w,w'\in {\cal F}^{(t,J)}$.  This is a contradiction to
Corollary 2.19.  So
$\tilde H^{(t,J)}$ is irreducible if it is an $\tilde H$-module.

It remains to show that the defining relations for $\tilde H$ are
satisfied.  This is accomplished as in the proof of [KR, Theorem 4.5].
The only relation which is tricky to check is the braid relation.
This can be verified as in [KR] or it can be checked by case by case 
arguments (as in [R2]).
\endpf

We summarize the results of this section with the following corollary
of Theorem 3.2 and the construction in Theorem 3.5.

\thm  Let $M$ be an irreducible calibrated $\tilde H$-module.  Let
$t\in T$ be (a fixed choice of) the central character of $M$
and let $J=R(w)\cap P(t)$ for any $w\in W$ such that $M_{wt}\ne 0$.
Then $(t,J)$ is a skew local region and $M\cong \tilde H^{(t,J)}$
where $\tilde H^{(t,J)}$ is the module defined in Theorem 3.5.
\endthm

\section 4. The structure of local regions

\smallskip
Recall that the Weyl group acts on 
$$T={\rm Hom}(X,\CC^*)=\{ \hbox{group homomorphisms $t\colon X \to \CC^*$}\},
\qquad\hbox{by}\quad
(wt)(X^\lambda) = t(X^{w^{-1}\lambda}).$$
Any element $t\in T$ is determined
by the values $t(X^{\omega_1}), t(X^{\omega_2}), \ldots, t(X^{\omega_n})$.
For $t\in T$ define the {\it polar decomposition} 
$$t=t_rt_c,\qquad\hbox{ $t_r,t_c\in T$ such that $t_r(X^\lambda)\in \RR_{>0}$,
and $|t_c(X^\lambda)|=1$,}$$
for all $X^\lambda\in X$.  There is a unique $\gamma\in \RR^n$ and a unique 
$\nu\in \RR^n/P$ such that
$$t_r(X^\lambda)=e^{\langle \gamma,\lambda\rangle}
\qquad\hbox{and}\qquad
t_c(X^\lambda)=e^{2\pi i\langle \nu,\lambda\rangle},
\qquad\hbox{for all $\lambda\in P$.}\formula$$
In this way we identify the sets $T_r=\{ t\in T\ |\ t=t_r\}$
and $T_c=\{ t\in T\ |\ t=t_c\}$ with $\fh_{\RR}^*$
and $\fh_{\RR}^*/P$, respectively.

For this paragraph (our goal here is (4.3) below)
assume that $q$ is not a root of unity (we will treat
the type A, root of unity case in detail in \S 7).
The representation theory of $\tilde H$ is ``the
same'' for any $q$ which is not a root of unity i.e., provided
$q$ is not a root of unity, the classification
and construction of simple $\tilde H$-modules can be stated
uniformly in terms of the parameter $q$.
Suppose $t\in T$ is such that $t=t_r$ and $\gamma\in \fh_{\RR}^*$ is
such that
$$t=e^{\gamma},
\qquad\hbox{in the sense that}\qquad
t(X^\lambda) = e^{\langle \gamma,\lambda\rangle},
\qquad\hbox{for all $X^\lambda\in X$}.$$
For the purposes of representation theory (as in Theorem 3.5) 
$t$ indexes a central character and so we should assume that $\gamma$
is chosen nicely in its $W$-orbit.   When
$$q=e\qquad\hbox{and}\qquad
\hbox{$\gamma$ is dominant,
i.e. $\langle\gamma,\alpha\rangle\ge0$ for all $\alpha\in R^+$,}
\formula$$
 then
$$
Z(t)=Z(\gamma),\qquad
P(t)=P(\gamma),
\qquad\hbox{and}\qquad
{\cal F}^{(t,J)} = {\cal F}^{(\gamma,J)}
\quad\hbox{for a subset $J\subseteq P(t)$, }$$
where 
$$Z(\gamma) = \{ \alpha\in R^+\ |\ \langle \gamma,\alpha\rangle=0\},
\qquad\hbox{and}\qquad
P(\gamma) = \{ \alpha\in R^+\ |\ \langle \gamma,\alpha\rangle=1\}
$$
$$
{\cal F}^{(\gamma,J)} = \{ w\in W\ |\ R(w)\cap Z(\gamma)=\emptyset,
\ \ R(w)\cap P(\gamma)=J\}.\formula$$
In this case the combinatorics of local regions is a new chapter
in the combinatorics of the Shi arrangement defined in (1.16).  
Other aspects of the combinatorics of the Shi arrangement
can be found in [Sh1-3], [St1-2], [AL], [ST], and there 
are several additional places in the literature [Sh3], 
[Xi,1.11, 2.6], [KOP], [Ks] 
which indicate that there that there is a deep (and not
yet completely understood) connection between the
structure and representation theory of the
affine Hecke algebra and the combinatorics of the Shi arrangement.

\subsection{Intervals in Bruhat order.}

Using the formulation in (4.3),
Theorem 4.6 will give a complete
description of the structure of ${\cal F}^{(\gamma,J)}$
as a subset of the Weyl group when $q$ is not 
a root of unity.  We will treat the type A, root of
unity cases in \S 7.

The {\it weak Bruhat order} is the partial order on $W$ given by
$$\hbox{$v\le w$\qquad if\qquad $R(v)\subseteq R(w)$},
\formula$$
where $R(w)$ denotes the inversion set of $w\in W$ as defined in (1.5).
This definition of the weak Bruhat order
is not the usual definition 
but is equivalent to the usual one by [Bj, Prop.\ 2].  
A set of positive roots $K$ is {\it closed} if  
$\alpha,\beta\in K$, $\alpha+\beta\in R^+$ implies
that $\alpha+\beta\in K$.  The {\it closure} $\overline{K}$ of a subset
$K\subseteq R^+$ is the smallest closed subset of $R^+$ containing $K$.
A set of positive roots $K\subseteq R^+$ is the inversion set
of some permutation $w\in W$ if and only if 
$K$ is closed and $K^c = R^+\backslash K$ is closed (see
[Bj, Prop.\ 2] or [KR, Theorem 5.1]).  

The following theorem is proved in [KR, \S 5].  The proof of
part (b) of the Theorem relies crucially on a
theorem of J. Losonczy [Lo].

\thm  Let $\gamma\in \fh_\RR^*$ be dominant 
(i.e. $\langle \gamma,\alpha\rangle\ge 0$
for all $\alpha\in R^+$) and let $J\subseteq P(\gamma)$.  
Let ${\cal F}^{(\gamma,J)}$ be as given in (4.3).
\smallskip\noindent
\item{(a)}
Then ${\cal F}^{(\gamma,J)}$ is nonempty if and only if
$J$ satisfies the condition
$$\hbox{if $\beta\in J$, $\alpha\in Z(\gamma)$
and $\beta-\alpha\in R^+$ 
\quad then\quad $\beta-\alpha\in J$.}$$
\item{(b)}
The sub-root system 
$R_{[\gamma]}=\{\alpha\in R\ |\ \langle \gamma,\alpha\rangle\in \ZZ\}$,
has Weyl group 
$$W_{[\gamma]} = \langle s_\alpha\ |\ \alpha\in R_{[\gamma]}\rangle
\qquad\hbox{and if}\qquad
W^{[\gamma]} = \{\sigma\in W\ |\ R(\sigma)\cap R_{[\gamma]}=\emptyset\}$$
then
$${\cal F}^{(\gamma,J)}
= W^{[\gamma]}\cdot[\tau_{\rm max},\tau_{\rm min}],$$
where $\tau_{\rm max},\tau_{\rm min}\in W_{[\gamma]}$
are determined by 
$$R(\tau_{\rm max})\cap R_{[\gamma]}=\overline{J}
\qquad\hbox{and}\qquad
R(\tau_{\rm min})\cap R_{[\gamma]} =
\overline{(P(\gamma)\backslash J)\cup Z(\gamma)}^c,$$ 
the complement is taken in the set of positive roots of $R_{[\gamma]}$,
and $[\tau_{\rm min},\tau_{\rm max}]$ denotes the interval between
$w_{\rm min}$ and $w_{\rm max}$ in the weak Bruhat order in $W_{[\gamma]}$.
\endthm

%

\subsection{Conjugation.}

Assume that $\gamma$ is dominant (i.e. $\langle \gamma,\alpha\rangle\ge 0$
for all $\alpha\in R^+$) and $J\subseteq P(\gamma)$.
Let ${\cal F}^{(\gamma,J)}$ be as given in (4.3).
The {\it conjugate} of $(\gamma,J)$ and of $w\in {\cal F}^{(\gamma,J)}$
are defined by
$$
(\gamma,J)'=(-u\gamma,-u(P(\gamma)\setminus J))
\qquad\hbox{and}\qquad
\matrix{
\cF^{(t,J)} &\mapleftright{1-1} &\cF^{(t,J)'} \cr
w &\longleftrightarrow &w'=wu^{-1}\,,\cr
}
\formula$$
where $u$ is the minimal length coset representative of
$w_0W_\gamma\in W/W_\gamma$ and $w_0$ is the longest element of $W$. 
In (6.7) we shall show that these maps
are generalizations of the classical conjugation 
operation on partitions.

\thm
The conjugation maps defined in (4.8) 
are well defined involutions.
\pf
\item{(a)}  Since $\gamma$ is dominant, $-u\gamma=-w_0\gamma$ is
dominant and thus $\langle -u\gamma, -u\alpha\rangle=1$
only if $-u\alpha>0$.  Thus
the equation
$\langle -u\gamma,-u\alpha\rangle =1 
\Longleftrightarrow \langle \gamma,\alpha\rangle = 1$
gives that 
$P(-u\gamma)=-uP(\gamma)$.
\smallskip
\item{(b)}  Let $v\in W_\gamma$ such that $w_0=uv$.  
(By [Bou, IV \S 1 Ex.\ 3], $v$ is unique.)
Then $R^+\supseteq -w_0Z(\gamma)=-uvZ(\gamma)=uZ(\gamma)$, and it follows 
that
$$Z(-u\gamma) 
=R^+ \cap \{ \alpha\in R\ |\ \langle u\gamma,\alpha\rangle=0\}
=R^+\cap (uZ(\gamma)\cup -uZ(\gamma))=uZ(\gamma).$$

\item{(c)}  Let $R^-=-R^+$ be the set of negative roots in $R$.
Let $v\in W_\gamma$ such that $w_0=uv$.  Then
$v$ is the longest element of $W_\gamma$ and $R(v)=Z(\gamma)$. Thus,
since $w_0R^-=R^+$,
$$
R(u)=
\{ \alpha\in R \ |\ \alpha\in R^+, w_0v\alpha\in R^-\} 
=\{ \alpha\in R \ |\ \alpha\in R^+, v\alpha\in R^+\},
= R^+\backslash R(v) = R^+\backslash Z(\gamma). 
$$
\item{(d)}
The weight $-u\gamma = -uv\gamma = -w_0\gamma$ is dominant
and  $-u(P(\gamma)\backslash J)\subseteq P(-u\gamma)$ since
$-uP(\gamma)=P(-u\gamma)$. This shows that $(\gamma,J)'$ is
well defined.
\smallskip
\item{(e)}
Write $w_0=uv$ where 
where $v$ is the longest element of $W_\gamma$.  Similarly,
write $w_0=u'v'$ where $u'$ is the minimal
length coset representative of $w_0W_{w_0\gamma}$ and $v'$ is
the longest element in $W_{w_0\gamma}$. 
Conjugation by $w_0$ is an involution on $W$ which takes simple reflections to
simple reflections and $W_{w_0\gamma}=w_0W_\gamma w_0$.  It follows that
$v'=w_0v w_0$.  This gives
$$u'u=(w_0v')(w_0v)=w_0w_0vw_0w_0v =1,$$
and so the second map in (4.8) is an involution.
\smallskip
\item{(f)}  Using (e) and (a),
$$
-u'(P(-u\gamma)\backslash(-u(P(\gamma)\backslash J)))
=-u'(-uP(\gamma)\backslash(-u(P(\gamma)\backslash J)))
=P(\gamma)\backslash(P(\gamma)\backslash J)=J,
$$
and so the first map in (4.8) is an involution.
\smallskip
\item{(g)}
Let $w\in \cF^{(\gamma,J)}$ and let $w'=wu^{-1}$.  
Since $R(w)\cap Z(\gamma)=\emptyset$, 
$$
\eqalign{
u^{-1}R(wu^{-1})\cap Z(\gamma)
&= \{\beta\in R\ |\ u\beta\in R(wu^{-1}), \beta\in Z(\gamma) \} \cr
&=\{\beta\in R\ |\ u\beta\in R^+,wu^{-1}u\beta\in  R^-,\beta\in Z(\gamma) \}\cr
&= \{\beta\in R\ |\ \beta\in u^{-1}R^+,w\beta\in  R^-,\beta\in Z(\gamma) \}\cr
&= \{\beta\in R\ |\ \beta\in u^{-1}R^+,\beta\in R(w), \beta\in Z(\gamma) \}, 
\quad\hbox{since $Z(\gamma)\subseteq R^+$,} \cr
&= \{\beta\in R\ |\ \beta\in u^{-1}R^+, \beta\in R(w)\cap Z(\gamma) \}\cr
&=\emptyset,\cr
}
$$
and thus, by (b),
$$R(w')\cap Z(-u\gamma)
=R(wu^{-1})\cap uZ(\gamma)
=u\left(u^{-1}R(wu^{-1})\cap Z(\gamma)\right)
=\emptyset.$$
Since $R(w)\cap P(\gamma)=J$,  
$$
\eqalign{
-u^{-1}R(wu^{-1})\cap P(\gamma)
&= \{\beta\in R\ |\ -u\beta\in R(wu^{-1}), \beta\in P(\gamma) \} \cr
&= \{\beta\in R\ |\ -u\beta\in R^+, -wu^{-1}u\beta\in  R^-,
\beta\in P(\gamma) \} \cr
&= \{\beta\in R\ |\ u\beta\in R^-, w\beta\in  R^+, \beta\in P(\gamma) \} \cr
&= \{\beta\in R\ |\ \beta\in R(u), \beta\in  R^+\backslash R(w),
\beta\in P(\gamma) \}, \quad\hbox{since $P(\gamma)\subseteq R^+$} \cr
&= \{\beta\in R\ |\ \beta\in R^+\backslash Z(\gamma),
\beta\in  R^+\backslash R(w), \beta\in P(\gamma) \} \cr
&=\{\beta\in R\ |\ \beta\in R^+\backslash Z(\gamma),
\beta\in  P(\gamma)\backslash J \},
\quad\hbox{since $R(w)\cap P(\gamma)=J$,} \cr
&= P(\gamma)\backslash J,
\quad\hbox{since $Z(\gamma)$ and $P(\gamma)$ are disjoint.} \cr
}
$$
Thus, by (a),
$$R(w')\cap P(-u\gamma)
=R(wu^{-1})\cap -u P(\gamma)
=-u\left( -u^{-1}R(wu^{-1})\cap P(\gamma) \right)
=-u\left(P(\gamma)\backslash J\right),$$
and so the second map in (4.8) is well defined.
\endpf

\remark
In type $A$, the conjugation involution coincides
with the duality operation for representations of 
$\gp$-adic $GL(n)$ defined by
Zelevinsky [Ze].  
Zelevinsky's involution has been studied
further in [MW], [KZ], [LTV] and extended to general Lie type 
by Kato [Kt2] and Aubert [Au].
For $\tilde H$-modules in type A, this is the involution on modules induced
by the Iwahori-Matsumoto involution of $\tilde H$ and
is detected on the level of characters: it sends an 
irreducible $\tilde H$-module
$L$ to the unique irreducible $L^*$ with
$\dim((L^*)_t^{\rm gen}) = \dim(L_{t^{-1}}^{\rm gen})$
for each $t\in T$.
I would like to thank J. Brundan for clarifying this remark and
making it precise.

\medskip\noindent
\subsection{Examples.}

\smallskip\noindent
\item{(a)}
If $\gamma$ is dominant and is generic (as an element of $C$) then 
$Z(\gamma)=P(\gamma)=\emptyset$ and
$\cF^{(\gamma,\emptyset)}=W$.
\smallskip\noindent
\item{(b)}
Let $\rho$ be defined by 
$\langle \rho,\alpha_i\rangle=1$, for all $1\le i\le n$.
Then
$$Z(\rho)=\emptyset,
\quad
P(\rho)=\{\alpha_1,\ldots,\alpha_n\},
\qquad\hbox{and}\qquad
\cF^{(\rho,J)}~=~\{ w\in W\ |\ D(w)=J\},$$
where
$D(w) = \{ \alpha_i \ |\ ws_i<w\}$
is the {\it right descent set} of $w\in W$.  The sets
${\cal F}^{(\gamma,J)}$ which arise here are fundamental to
the theory of descent algebras [So], [GR], [Re].
\smallskip\noindent
\item{(c)}  This example is a generalization of (b).
Suppose that $(\gamma,J)$ is a local region such that $\gamma$ is
regular and integral
(i.e. $\langle \gamma,\alpha\rangle\in \ZZ_{>0}$ for all $\alpha\in R^+$).
Then
$$Z(\gamma)=\emptyset,\quad
P(\gamma)\subseteq \{\alpha_1,\ldots,\alpha_n\},
\qquad\hbox{and}\qquad
\cF^{(\gamma,J)}~=~\{ w\in W\ |\ D(w)\cap P(\gamma)=J\}.$$
\smallskip\noindent
\item{(d)}  
Let $R$ be the root system of type $C_2$ with simple roots
$\alpha_1=\varepsilon_1$
$\alpha_2=\varepsilon_2-\varepsilon_1,$
where $\{\varepsilon_1,\varepsilon_2\}$ is an orthonormal basis
of $\fh_{\RR}^*=\RR^2$.  
The positive roots are 
$R^+ = \{ \alpha_1,\alpha_2,\alpha_1+\alpha_2,\alpha_1+2\alpha_2\}$.
Let $\gamma\in \RR^2$ be given by
$\langle \gamma,\alpha_1\rangle = 0$
and
$\langle \gamma,\alpha_2\rangle =1$.
Then $\gamma$ is dominant (i.e. in $\overline{C}$) and integral and
$$Z(\gamma)=\{\alpha_1\}
\qquad\hbox{and}\qquad
P(\gamma)=\{\alpha_2,\alpha_1+\alpha_2\}.$$
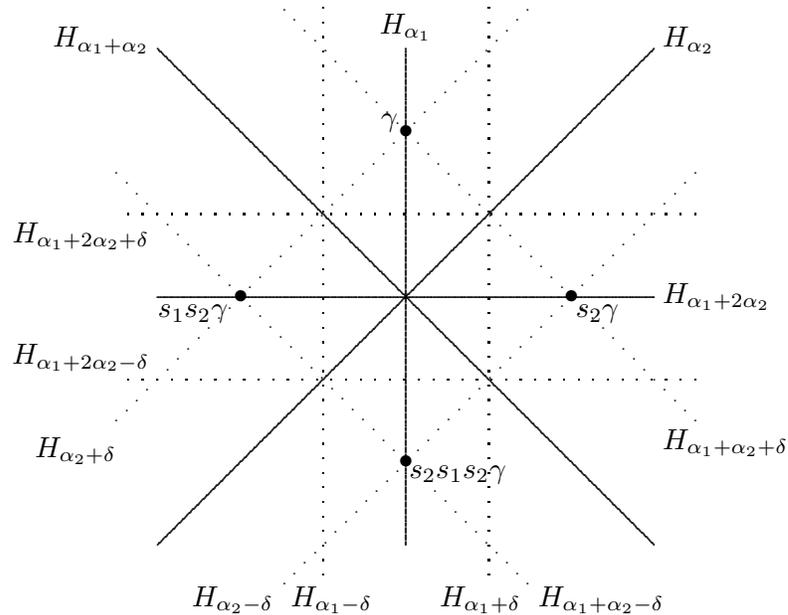
$$
\matrix{
\beginpicture
\setcoordinatesystem units <1cm,1cm>         
\setplotarea x from -4 to 4, y from -4 to 4    
\put{$H_{\alpha_1}$}[b] at 0 3.1
\put{$H_{\alpha_2}$}[l] at 3.1 3.1
\put{$H_{\alpha_1+\alpha_2}$}[r] at -3.1 3.1
\put{$H_{\alpha_1+2\alpha_2}$}[l] at 3.1 0
\put{$H_{\alpha_1+2\alpha_2+\delta}$}[tr] at -3.1 0.9
\put{$H_{\alpha_1+2\alpha_2-\delta}$}[br] at -3.1 -0.9
\put{$H_{\alpha_1+\alpha_2+\delta}$}[tl] at 3.1 -1.6
\put{$H_{\alpha_1+\alpha_2-\delta}$}[tl] at 1.6 -3.5
\put{$H_{\alpha_2+\delta}$}[tr] at -3.5 -1.7
\put{$H_{\alpha_2-\delta}$}[tr] at -1.6 -3.5
\put{$H_{\alpha_1-\delta}$}[t] at -0.9 -3.5
\put{$H_{\alpha_1+\delta}$}[t] at  0.9 -3.5
\put{$\bullet$} at 0 2
\put{$\gamma$}[r] at -0.1 2.1
\put{$\bullet$} at 2 0
\put{$s_2\gamma$}[tl] at 2.05 -0.1
\put{$\bullet$} at -2 0 
\put{$s_1s_2\gamma$}[tl] at -3 -0.1
\put{$\bullet$} at 0 -2
\put{$s_2s_1s_2\gamma$}[tl] at .04 -2.04 
\plot -3 -3   3 3 /
\plot  3 -3  -3 3 /
\plot  0  3   0 -3 /
\plot  3  0  -3  0 /
\setdots
\putrule from 1 3.5 to 1 -3.5
\putrule from -1 3.5 to -1 -3.5
\putrule from 3.5 1 to -3.5  1
\putrule from 3.5 -1 to -3.5 -1
\plot 3.5 1.5  -1.5 -3.5 /
\plot -3.5 -1.5  1.5 3.5 /
\plot -3.5 1.5  1.5 -3.5 /
\plot 3.5 -1.5  -1.5 3.5 /
\endpicture
\cr
\cr
}$$
The following picture displays the local regions ${\cal F}^{(\gamma,J})$
as regions in $\fh_{\RR}^*$, see the remarks after (2.18).
$$
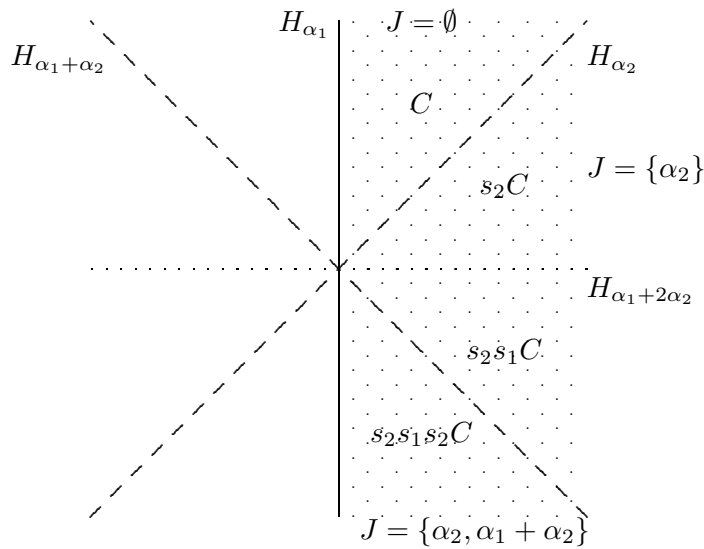
\beginpicture
\setcoordinatesystem units <1cm,1cm>         
\setplotarea x from -3 to 3, y from -3 to 3    
\put{$H_{\alpha_1}$}[r] at -0.1 2.95
\put{$H_{\alpha_2}$}[tl] at 3 2.7
\put{$H_{\alpha_1+\alpha_2}$}[tr] at -2.8 2.7
\put{$H_{\alpha_1+2\alpha_2}$}[tl] at 3 -0.1
\put{$J=\emptyset$} at 1 3
\put{$J=\{\alpha_2\}$}[l] at 3 1.2
\put{$J=\{\alpha_2,\alpha_1+\alpha_2\}$}[tl] at 0.25 -3.0
\put{$C$} at 1 2
\put{$s_2C$} at 2 1
\put{$s_2s_1C$} at 2 -1
\put{$s_2s_1s_2C$} at  1 -2
\putrule from 0 3 to 0 -3
\setdashes
\plot -3 -3   3 3 /
\plot  3 -3  -3 3 /
\setdots
\putrule from 3 0 to -3  0
\vshade 0 -3 3    3 -3 3 /
\endpicture
$$
The solid line is the hyperplane corresponding to the root in 
$Z(\gamma)$ and the dashed lines are the hyperplanes 
corresponding to the roots in $P(\gamma)$.
\smallskip\noindent
\item{(e)}
Let $R$ be the root system of type $C_2$ as
in (d).  Let $\gamma\in \RR^2$ be defined by
$$\langle \gamma, \alpha_1\rangle = 0,
\qquad
\langle \gamma,\alpha_2\rangle = \hbox{${1\over 2}$}.
\qquad
\hbox{Then}\quad
Z(\gamma)=\{\alpha_1\}, \quad P(\gamma)=\{\alpha_1+2\alpha_2\}.$$
If $J=P(\gamma)$ then the unique minimal element 
$w_{\rm min}$ of $\cF^{(\gamma,J)}$ has 
$R(w_{\rm min})=\{\alpha_2,\alpha_1+2\alpha_2\}\ne \overline{J}=J.$
%
%
%
%

\section 5. The connection to standard Young tableaux

In this section we shall show that the combinatorics
of local regions is a generalization of the
combinatorics of standard Young tableaux.
Let us first make some general definitions, which we will
show later provide generalizations of standard objects 
in the Young tableaux theory.  This section is a (purely
combinatorial) study of the local regions in the form
which appears in (4.3), and therefore corresponds to
the representation theory of affine Hecke algebras when
$q$ is not a root of unity.

\subsection{Definitions.}  

Let $\gamma\in \fh_{\RR}^*$ be dominant 
and let 
$$Z(\gamma)= \{ \alpha\in R^+\ |\ \langle \gamma,\alpha\rangle=0\},
\qquad
P(\gamma)= \{ \alpha\in R^+\ |\ \langle \gamma,\alpha\rangle=1\},
$$
$$
{\cal F}^{(\gamma,J)} = \{ w\in W\ |\ R(w)\cap Z(\gamma)=\emptyset,
\ \ R(w)\cap P(\gamma)=J\},$$
as in (4.3).  
\smallskip\noindent
\item{(a)}
A {\it local region} is a pair $(\gamma,J)$ such that 
${\cal F}^{(\gamma,J)}$ is nonempty.
\smallskip\noindent
\item{(b)}
A {\it ribbon} is a local region $(\gamma,J)$ such that 
$\gamma$ is regular, i.e. $\langle \gamma,\alpha\rangle\ne 0$
for all $\alpha \in R$.  
\smallskip\noindent
\item{(c)}
An element $\gamma\in \overline{C}$ is {\it calibratable} if
\smallskip
\itemitem{(1)} For all simple roots $\alpha_i$,
$1\le i\le n$,\ \ $\langle \gamma,\alpha_i\rangle \ne 0$, and
\smallskip
\itemitem{(2)} For all pairs of simple roots $\alpha_i$ and $\alpha_j$
such that 
$\{\alpha\in R_{ij}\ |\ \langle\gamma,\alpha\rangle=0\}\ne \emptyset$,
\itemitem{} the set 
$\{\alpha\in R_{ij}\ |\ \langle\gamma,\alpha\rangle=1\}$
contains more than two elements.
\smallskip\noindent
\item{(d)}
A {\it skew local region} is a local region $(\gamma,J)$
such that $w\gamma$ is calibratable
for all $w\in {\cal F}^{(\gamma,J)}$.  All ribbons are skew.
\smallskip\noindent
\item{(e)}
A {\it column (resp. row) reading tableau} is a
minimal (resp. maximal) element of $\cF^{(\gamma,J)}$ in the weak 
Bruhat order.   
\smallskip\noindent
\item{(f)}
If $\alpha\in R$ the {\it $\alpha$-axial distance} 
for $w\in {\cal F}^{(\gamma,J)}$ is the value
$d_\alpha(w) = \langle w\gamma, \alpha\rangle.$

{\it Remarks.}
(1) Theorem 4.6(b) shows
that, up to a shift, the set ${\cal F}^{(\gamma,J)}$ has
a unique maximal and a unique minimal element 
and is an interval in the weak Bruhat order.  This is the fundamental
importance of the notions of the row reading and the column reading
tableaux. Theorem 6.9 in Section 6 will show how Theorem 4.6(b)
is a generalization of a Young tableaux result of Bj\"orner and Wachs [BW, 
Theorem 7.2]. 
\smallskip\noindent
(2)  The definition of skew local regions is forced by
the representation theory of the affine Hecke algebra (see 
Theorem 3.6, the classification of irreducible calibrated
representations).  In Proposition 6.4 below we shall show
that the skew local regions and the ribbons
are generalizations of the
skew shapes and border strips which
are used in the theory of symmetric functions
[Mac, I \S 5 and I \S 3 Ex. 11]
\smallskip\noindent
(3) The axial distances control the denominators
which appear in the construction of irreducible
representations of the affine Hecke algebra in 
Theorem 3.5.  In (6.1) we shall see how they are analogues
of the axial distances used by A. Young [Yg2] in his constructions
of the irreducible representations of the symmetric group.

\smallskip\noindent
To summarize, a brief dictionary between local regions combinatorics and 
the Young tableaux combinatorics:
$$\matrix{
\hfill \hbox{skew local regions} 
&\qquad\longleftrightarrow\qquad
&\hbox{skew shapes $\lambda/\mu$} \hfill \cr
\hfill \hbox{ribbons} 
&\qquad\longleftrightarrow\qquad
&\hbox{border strips} \hfill \cr
\hfill \hbox{local regions}
&\qquad\longleftrightarrow\qquad
&\hbox{general configurations of boxes} \hfill \cr
\hfill {\cal F}^{(\gamma,J)}
&\qquad\longleftrightarrow\qquad
&\hbox{the set of standard tableaux ${\cal F}^{\lambda/\mu}$}\hfill  \cr
}
$$
The remainder of this section and the next section explain in greater
detail the conversions indicated in this dictionary.

\subsection{The root system.}

Let $\{\varepsilon_1,\ldots,\varepsilon_n\}$ be an
orthonormal basis  of $\fh_{\RR}^*=\RR^n$ so that each sequence
$\gamma=(\gamma_1,\ldots,\gamma_n)\in \RR^n$
is identified with the vector $\gamma=\sum_i \gamma_i\varepsilon_i$.
The root system of type $A_{n-1}$ is given by the sets
$$
R = \left\{\pm(\varepsilon_j-\varepsilon_i)
\,|\, 1\leq i,j\leq n\right\}
\quad\hbox{and}\quad
R^+=\left\{\varepsilon_j-\varepsilon_i\,|\, 1\leq i<j\leq n\right\}.
$$
The Weyl group is $W=S_n$, the symmetric group, acting by 
permutations of the $\varepsilon_i$.

\subsection{Partitions, skew shapes, and standard tableaux.}

A partition $\lambda$ is a collection of $n$
boxes in a  corner.  We shall conform to the conventions in 
[Mac] and assume
that gravity goes up and to the left. 

$$
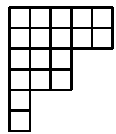
\beginpicture
\setcoordinatesystem units <0.25cm,0.25cm>         
\setplotarea x from 0 to 4, y from 0 to 3    
\linethickness=0.5pt                          
\putrule from 0 6 to 5 6          %
\putrule from 0 5 to 5 5          
\putrule from 0 4 to 5 4          %
\putrule from 0 3 to 3 3          %
\putrule from 0 2 to 3 2          %
\putrule from 0 1 to 1 1          %
\putrule from 0 0 to 1 0          %

\putrule from 0 0 to 0 6        %
\putrule from 1 0 to 1 6        %
\putrule from 2 2 to 2 6        %
\putrule from 3 2 to 3 6        
\putrule from 4 4 to 4 6        %
\putrule from 5 4 to 5 6        %
\endpicture
$$
Any partition $\lambda$ can be identified with
the sequence $\lambda=(\lambda_1\ge \lambda_2\ge \ldots )$
where $\lambda_i$ is the number of boxes in row $i$ of $\lambda$.
The rows and columns are numbered in the same way as for matrices.
We shall always use the word {\it diagonal} to mean a major diagonal.
In the example above $\lambda=(553311)$ and the diagonals of $\lambda$
(from southwest to northeast) contain 1,1,1,2,3,3,2,2,2, and 1 box
respectively.

If $\lambda$ and $\mu$ are partitions such that $\mu_i\le \lambda_i$
for all $i$ write $\mu\subseteq \lambda$.  
The {\it skew shape} $\lambda/\mu$
consists of all boxes of $\lambda$ which are not in $\mu$.
Let $\lambda/\mu$ be a skew shape with $n$ boxes.
Number the boxes of each skew shape
$\lambda/\mu$ along diagonals from southwest to northeast and
$$\hbox{write ${\rm box}_i$ to indicate the box numbered $i$.}$$
See Example 5.8 below.
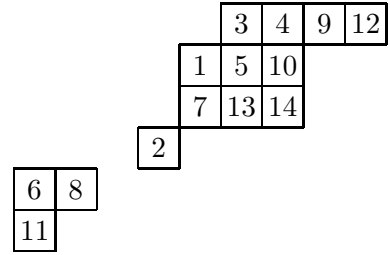
A {\it standard tableau of shape} $\lambda/\mu$ is a filling  
of the boxes in the skew shape
$\lambda/\mu$ with the numbers $1,\ldots,n$ such that 
the numbers increase from left to right in each row and from top to
bottom down each column. 
Let $\cF^{\lambda/\mu}$ be the set of standard tableaux of shape
$\lambda/\mu$.
Given a standard tableau $p$ of shape $\lambda/\mu$ define
the {\it word} of $p$ to be the permutation
$$
w_p=\left(\matrix{
1&\cdots & n \cr
p({\rm box}_1)&\ldots &p({\rm box}_n)
}\right)
\formula
$$ 
where $p({\rm box}_i)$ is the entry in ${\rm box}_i$ of the standard
tableau.

\subsection{Placed skew shapes.}

Let $\lambda/\mu$ be a skew shape with $n$ boxes.  Imagine placing
$\lambda/\mu$ on a piece of infinite graph paper where the diagonals of the 
graph paper are indexed consecutively (with elements of $\ZZ$) from southwest
to northeast.  
$$
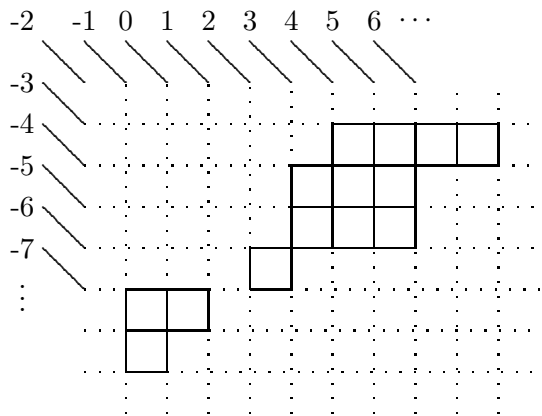
\matrix{
\beginpicture
\setcoordinatesystem units <0.5cm,0.5cm>         
\setplotarea x from -1 to 10, y from -1 to 7    
\linethickness=0.5pt                          
\put{$\vdots$}  at -2.5 2
\put{-7} at -2.5 3
\put{-6} at -2.5 4
\put{-5} at -2.5 5
\put{-4}  at -2.5 6
\put{-3}  at -2.5 7
\put{-2}  at -2.5 8.5
\put{-1}  at -1 8.5
\put{0} at 0 8.5
\put{1} at 1 8.5
\put{2}  at 2 8.5
\put{3}  at 3 8.5
\put{4}  at 4 8.5
\put{5}  at 5 8.5
\put{6}  at 6 8.5
\put{$\cdots$}  at 7 8.5
\putrule from 5 6 to 9 6          %
\putrule from 4 5 to 9 5          
\putrule from 4 4 to 7 4          %
\putrule from 3 3 to 7 3          %
\putrule from 3 2 to 4 2          %
\putrule from 0 2 to 2 2          %
\putrule from 0 1 to 2 1          %
\putrule from 0 0 to 1 0          %
\putrule from 0 0 to 0 2        %
\putrule from 1 0 to 1 2        %
\putrule from 2 1 to 2 2        %
\putrule from 3 2 to 3 3        
\putrule from 4 2 to 4 5        %
\putrule from 5 3 to 5 6        %
\putrule from 6 3 to 6 6        %
\putrule from 7 3 to 7 6        %
\putrule from 8 5 to 8 6        %
\putrule from 9 5 to 9 6        %
\plot -2 3  -1 2 /
\plot -2 4  -1 3 /
\plot -2 5  -1 4 /
\plot -2 6  -1 5 /
\plot -2 7  -1 6 /
\plot -2 8  -1 7 /
\plot -1 8  0 7 /
\plot 0 8  1 7 /
\plot 1 8  2 7 /
\plot 2 8  3 7 /
\plot 3 8  4 7 /
\plot 4 8  5 7 /
\plot 5 8  6 7 /
\plot 6 8  7 7 /
\setdots
\putrule from 9 6 to 9 7
\putrule from 8 6 to 8 7
\putrule from 7 6 to 7 7
\putrule from 6 6 to 6 7
\putrule from 5 6 to 5 7
\putrule from 4 5 to 4 7
\putrule from 3 3 to 3 7
\putrule from 2 2 to 2 7
\putrule from 1 2 to 1 7
\putrule from 0 2 to 0 7
\putrule from 9 -1 to 9 5
\putrule from 8 -1 to 8 5
\putrule from 7 -1 to 7 3
\putrule from 6 -1 to 6 3
\putrule from 5 -1 to 5 3
\putrule from 4 -1 to 4 2
\putrule from 3 -1 to 3 2
\putrule from 2 -1 to 2 1
\putrule from 1 -1 to 1 0
\putrule from 0 -1 to 0 0
\putrule from -1 6 to 5 6
\putrule from  9 6 to 10 6
\putrule from -1 5 to 4 5
\putrule from  9 5 to 10 5
\putrule from -1 4 to 4 4
\putrule from  7 4 to 10 4
\putrule from -1 3 to 3 3
\putrule from  7 3 to 10 3
\putrule from -1 2 to 0 2
\putrule from  2 2 to 3 2
\putrule from  4 2 to 10 2
\putrule from -1 1 to 0 1
\putrule from  2 1 to 10 1
\putrule from -1 0 to 0 0
\putrule from  1 0 to 10 0
\endpicture
\cr
}
$$
The {\it content} of a box $b$ is 
$$c(b) = \hbox{diagonal number of box $b$}.$$
Identify the sequence
$$
\gamma=(c({\rm box}_1),c({\rm box}_2),\ldots,c({\rm box}_n))
\qquad\hbox{with}\qquad
\gamma=\sum_{i=1}^n c({\rm box}_i)\varepsilon_i \in \RR^n.
\formula$$
The pair $(\gamma,\lambda/\mu)$ is a {\it placed skew shape}.
It follows from the definitions in (5.1) that
$$
\eqalign{
Z(\gamma) &= \{ \varepsilon_j-\varepsilon_i \ |\ 
\hbox{$j>i$ and 
${\rm box}_j$ and ${\rm box}_i$ are in the same diagonal} \}, 
\quad\hbox{and} \cr
P(\gamma) &= \{ \varepsilon_j-\varepsilon_i \ |\ 
\hbox{$j>i$ and 
${\rm box}_j$ and ${\rm box}_i$ are in adjacent diagonals} \}. \cr
}
$$
Define 
$$
J=\left\{\varepsilon_j-\varepsilon_i\ \  \Bigg|\ \ 
\matrix{
j>i \hfill \cr 
\hbox{${\rm box}_j$ and ${\rm box}_i$ are in adjacent diagonals} \hfill \cr
\hbox{${\rm box}_j$ is northwest of ${\rm box}_i$ } \hfill\cr
} 
\right\},
\formula
$$
where ${\it northwest}$ means strictly north and weakly west.

\subsection{Example.}

The following diagrams illustrate standard tableaux and
the numbering of boxes in a skew shape $\lambda/\mu$.
$$
\matrix{
\beginpicture
\setcoordinatesystem units <0.5cm,0.5cm>         
\setplotarea x from 0 to 4, y from 0 to 3    
\linethickness=0.5pt                          
\put{1} at 0.5 0.5
\put{2} at 0.5 1.5
\put{3} at 1.5 1.5
\put{4}  at 3.5 2.5
\put{5} at 4.5 3.5
\put{6} at 4.5 4.5
\put{7}  at 5.5 3.5
\put{8}  at 5.5 4.5
\put{10}  at 5.5 5.5
\put{9}  at 6.5 3.5
\put{11}  at 6.5 4.5
\put{12}  at 6.5 5.5
\put{13}  at 7.5 5.5
\put{14}  at 8.5 5.5
\putrule from 5 6 to 9 6          %
\putrule from 4 5 to 9 5          
\putrule from 4 4 to 7 4          %
\putrule from 3 3 to 7 3          %
\putrule from 3 2 to 4 2          %
\putrule from 0 2 to 2 2          %
\putrule from 0 1 to 2 1          %
\putrule from 0 0 to 1 0          %
\putrule from 0 0 to 0 2        %
\putrule from 1 0 to 1 2        %
\putrule from 2 1 to 2 2        %
\putrule from 3 2 to 3 3        
\putrule from 4 2 to 4 5        %
\putrule from 5 3 to 5 6        %
\putrule from 6 3 to 6 6        %
\putrule from 7 3 to 7 6        %
\putrule from 8 5 to 8 6        %
\putrule from 9 5 to 9 6        %
\endpicture
&\qquad\qquad
&
\beginpicture
\setcoordinatesystem units <0.5cm,0.5cm>         
\setplotarea x from 0 to 4, y from 0 to 3    
\linethickness=0.5pt                          
\put{11} at 0.5 0.5
\put{6} at 0.5 1.5
\put{8} at 1.5 1.5
\put{2}  at 3.5 2.5
\put{7} at 4.5 3.5
\put{1} at 4.5 4.5
\put{13}  at 5.5 3.5
\put{5}  at 5.5 4.5
\put{3}  at 5.5 5.5
\put{14}  at 6.5 3.5
\put{10}  at 6.5 4.5
\put{4}  at 6.5 5.5
\put{9}  at 7.5 5.5
\put{12}  at 8.5 5.5
\putrule from 5 6 to 9 6          %
\putrule from 4 5 to 9 5          
\putrule from 4 4 to 7 4          %
\putrule from 3 3 to 7 3          %
\putrule from 3 2 to 4 2          %
\putrule from 0 2 to 2 2          %
\putrule from 0 1 to 2 1          %
\putrule from 0 0 to 1 0          %
\putrule from 0 0 to 0 2        %
\putrule from 1 0 to 1 2        %
\putrule from 2 1 to 2 2        %
\putrule from 3 2 to 3 3        
\putrule from 4 2 to 4 5        %
\putrule from 5 3 to 5 6        %
\putrule from 6 3 to 6 6        %
\putrule from 7 3 to 7 6        %
\putrule from 8 5 to 8 6        %
\putrule from 9 5 to 9 6        %
\endpicture
\cr
\hbox{$\lambda/\mu$ with boxes numbered}
&&\hbox{A standard tableau $p$ of shape $\lambda/\mu$} \cr
}
$$
The word of the standard tableau $p$ is the permutation
$w_p=(11,6,8,2,7,1,13,5,14,3,10,4,9,12)$ (in one-line 
notation).

The following picture shows the contents of the boxes in the
placed skew shape $(\gamma,\lambda/\mu)$ with
$\gamma = (-7,-6,-5,-2,0,1,1,2,2,3,3,4,5,6)$.
$$
\matrix{
\beginpicture
\setcoordinatesystem units <0.5cm,0.5cm>         
\setplotarea x from 0 to 4, y from 0 to 3    
\linethickness=0.5pt                          
\put{-7} at 0.5 0.5
\put{-6} at 0.5 1.5
\put{-5} at 1.5 1.5
\put{-2}  at 3.5 2.5
\put{0} at 4.5 3.5
\put{1} at 4.5 4.5
\put{1}  at 5.5 3.5
\put{2}  at 5.5 4.5
\put{3}  at 5.5 5.5
\put{2}  at 6.5 3.5
\put{3}  at 6.5 4.5
\put{4}  at 6.5 5.5
\put{5}  at 7.5 5.5
\put{6}  at 8.5 5.5
\putrule from 5 6 to 9 6          %
\putrule from 4 5 to 9 5          
\putrule from 4 4 to 7 4          %
\putrule from 3 3 to 7 3          %
\putrule from 3 2 to 4 2          %
\putrule from 0 2 to 2 2          %
\putrule from 0 1 to 2 1          %
\putrule from 0 0 to 1 0          %
\putrule from 0 0 to 0 2        %
\putrule from 1 0 to 1 2        %
\putrule from 2 1 to 2 2        %
\putrule from 3 2 to 3 3        
\putrule from 4 2 to 4 5        %
\putrule from 5 3 to 5 6        %
\putrule from 6 3 to 6 6        %
\putrule from 7 3 to 7 6        %
\putrule from 8 5 to 8 6        %
\putrule from 9 5 to 9 6        %
\endpicture
\cr
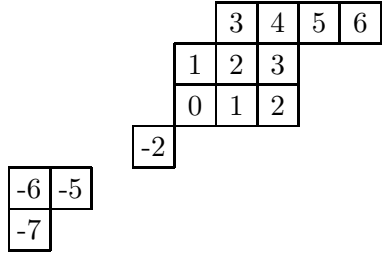
\hbox{Contents of the boxes of $(\gamma,\lambda/\mu)$} \cr
}
$$
In this case
$J=\{
\varepsilon_2-\varepsilon_1,
\varepsilon_6-\varepsilon_5,
\varepsilon_8-\varepsilon_7,
\varepsilon_{10}-\varepsilon_8,
\varepsilon_{10}-\varepsilon_9,
\varepsilon_{11}-\varepsilon_9,
\varepsilon_{12}-\varepsilon_{11}\}.$
\endpf

\thm  Let $(\gamma,\lambda/\mu)$ be a placed skew shape and
let $J$ be as defined in (5.7).  Let
$\cF^{\lambda/\mu}$ be the set of standard tableaux of shape $\lambda/\mu$
and let $\cF^{(\gamma,J)}$ be the set defined in (5.1).
Then the map 
$$
\matrix{
\cF^{\lambda/\mu} &\mapleftright{1-1} & \cF^{(\gamma,J)} \cr
p & \longleftrightarrow & w_p, \cr
}
$$
where $w_p$ is as defined in (5.4),
is a bijection.
\pf
If $w=(w(1)\cdots w(n))$ is a permutation in $S_n$ then 
$$R(w)=\{\varepsilon_j-\varepsilon_i\ |\ 
\hbox{$j>i$ such that $w(j)<w(i)$} \ \}.$$  
The theorem is a consequence of the following chain of equivalences:

\medskip\noindent
The filling $p$ is a standard tableau if and only if for all $1\le i<j\le n$
\item{(a)} $p({\rm box}_i)<p({\rm box}_j)$ 
if ${\rm box}_i$ and ${\rm box}_j$ are on the same diagonal,
\item{(b)} $p({\rm box}_i)<p({\rm box}_j)$ 
if ${\rm box}_j$ is immediately to the right of ${\rm box}_i$, and
\item{(c)} $p({\rm box}_i)>p({\rm box}_j)$ 
if ${\rm box}_j$ is  immediately above ${\rm box}_i$.
\smallskip\noindent
These conditions hold if and only if
\item{(a)} $\varepsilon_j-\varepsilon_i\not\in R(w_p)$ if
$\varepsilon_j-\varepsilon_i\in Z(\gamma)$,
\item{(b)} $\varepsilon_j-\varepsilon_i\not\in R(w_p)$ if
$\varepsilon_j-\varepsilon_i\in P(\gamma)\setminus J$,
\item{(c)} $\varepsilon_j-\varepsilon_i\in R(w_p)$ if
$\varepsilon_j-\varepsilon_i\in J$,
\smallskip\noindent
which hold if and only if
$$\hbox{(a)\enspace $\alpha\not\in R(w_p)$ if $\alpha\in Z(\gamma)$,}
\quad
\hbox{(b)\enspace $\alpha\not\in R(w_p)$ if $\alpha\in
P(\gamma)\setminus J$,}
\quad\hbox{and}\quad
\hbox{(c)\enspace $\alpha\in R(w_p)$ if $\alpha\in J$.}
$$
Finally, these are equivalent to the conditions
$R(w_p)\cap Z(\gamma)=\emptyset$ and 
$R(w_p)\cap P(\gamma)=J$.
\endpf

\subsection{Placed configurations.}

We have described how one can identify placed
skew shapes $(\gamma,\lambda/\mu)$ with certain pairs $(\gamma, J)$.
One can extend this conversion to associate placed configurations of boxes
to more general pairs $(\gamma,J)$.  The resulting configurations
are not always skew shapes.

Let $(\gamma,J)$ be a pair such that $\gamma=(\gamma_1,\ldots,\gamma_n)$
is a dominant integral weight and $J\subseteq P(\gamma)$. 
(The sequence $\gamma$ is a dominant integral weight
if $\gamma_1\leq\cdots\leq\gamma_n$ and $\gamma_i\in\ZZ$ for all $i$.)
If $J$ satisfies the condition
$$\hbox{{\sl If $\beta\in J$, $\alpha\in Z(\gamma)$, and $\beta-\alpha\in R^+$
then $\beta-\alpha\in J$}}
$$
then $(\gamma,J)$ will determine a placed configuration of boxes
(see Theorem 4.6).
As in the placed skew shape case, think of the boxes as being placed on graph
paper where the  boxes on a given diagonal all have the same content.  
(The boxes on each diagonal are allowed to slide along the diagonal
as long as they don't pass through the corner of a box on an adjacent diagonal.)
The sequence $\gamma$ describes how many boxes are on each diagonal and the set
$J$ determines how the boxes on adjacent diagonals are placed relative to each
other.  We want
$$\gamma=\sum_{i=1}^n c({\rm box}_i)\varepsilon_i,$$
and
\itemitem{(a)} If $\varepsilon_j-\varepsilon_i\in J$ then ${\rm box}_j$ is
northwest of ${\rm box}_i$, and
\itemitem{(b)} If $\varepsilon_j-\varepsilon_i\in P(\gamma)\backslash J$ then
${\rm box}_j$ is southeast of ${\rm box}_i$,
\smallskip\noindent
where the boxes are numbered along diagonals in the same way as 
for skew shapes,
{\it southeast} means weakly south and strictly east, 
and {\it northwest} means strictly north and weakly west. 

If we view the pair $(\gamma,J)$ as a placed configuration of boxes then
the {\it standard tableaux} are fillings $p$ of the $n$ boxes in the
configuration with $1,2,\ldots, n$ such that for all $i<j$
\item{(a)} $p({\rm box}_i)<p({\rm box}_j)$ 
if ${\rm box}_i$ and ${\rm box}_j$ are on the same diagonal,
\item{(b)} $p({\rm box}_i)<p({\rm box}_j)$ if ${\rm box}_i$ and 
${\rm box}_j$ are on adjacent diagonals and ${\rm box}_j$ is southeast 
of ${\rm box}_i$, and
\item{(c)} $p({\rm box}_i)>p({\rm box}_j)$ if ${\rm box}_i$ and 
${\rm box}_j$ are on adjacent diagonals and ${\rm box}_j$ is northwest 
of ${\rm box}_i$.
\smallskip\noindent
As in Theorem 5.6 the permutation in $\cF^{(\gamma,J)}$ which corresponds
to the standard tableau $p$ is
$w_p=(p({\rm box}_1),\ldots,p({\rm box}_n))$.
The following example illustrates the conversion.

\bigskip\noindent
{\sl Example.}
Suppose $\gamma=(-1,-1,-1,0,0,0,1,1,1,2,2,2)$ and
$$\eqalign{
J &=\left\{ \varepsilon_4-\varepsilon_1, \varepsilon_4-\varepsilon_2,
\varepsilon_4-\varepsilon_3, \varepsilon_5-\varepsilon_2,
\varepsilon_5-\varepsilon_3, \varepsilon_7-\varepsilon_5,
\varepsilon_7-\varepsilon_6, \varepsilon_8-\varepsilon_6,
\varepsilon_{10}-\varepsilon_9, \varepsilon_{10}-\varepsilon_8,\right. \cr
&\phantom{J = } \left.\varepsilon_{10}-\varepsilon_7,
\varepsilon_{11}-\varepsilon_9,
\varepsilon_{11}-\varepsilon_8, \varepsilon_{11}-\varepsilon_7,
\varepsilon_{12}-\varepsilon_9 \right\}. \cr
}
$$
The placed configuration of boxes corresponding to $(\gamma,J)$
is as given below.
$$
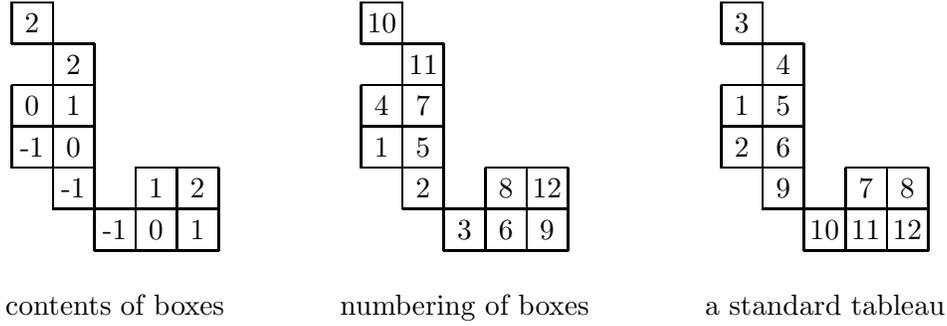
\matrix{
\beginpicture
\setcoordinatesystem units <0.5cm,0.5cm> point at 6 0       
\setplotarea x from 0 to 5, y from -1 to 6    
\linethickness=0.5pt                          
\put{-1}  at 0.5 2.5    %
\put{-1}  at 1.5 1.5    %
\put{-1}  at 2.5 0.5    %
\put{0}  at 0.5 3.5    %
\put{0}  at 1.5 2.5    
\put{0}  at 3.5 0.5    %
\put{1}  at 1.5 3.5    %
\put{1}  at 3.5 1.5    %
\put{1}  at 4.5 0.5    %
\put{2}  at 0.5 5.5    %
\put{2}  at 1.5 4.5    %
\put{2}  at 4.5 1.5    %
\putrule from 0 6 to 1 6          %
\putrule from 0 5 to 2 5          %
\putrule from 0 4 to 2 4          %
\putrule from 0 3 to 2 3          %
\putrule from 3 2 to 5 2          %
\putrule from 0 2 to 2 2          
\putrule from 1 1 to 5 1          %
\putrule from 2 0 to 5 0          %
\putrule from 0 2 to 0 4        %
\putrule from 0 5 to 0 6        %
\putrule from 1 1 to 1 6        %
\putrule from 2 0 to 2 5        
\putrule from 3 0 to 3 2        %
\putrule from 4 0 to 4 2        %
\putrule from 5 0 to 5 2        %
\endpicture
&\qquad
&\beginpicture
\setcoordinatesystem units <0.5cm,0.5cm> point at -1 0       
\setplotarea x from 0 to 5, y from -1 to 6    
\linethickness=0.5pt                          
\put{1}  at 0.5 2.5    %
\put{2}  at 1.5 1.5    %
\put{3}  at 2.5 0.5    %
\put{4}  at 0.5 3.5    %
\put{5}  at 1.5 2.5    
\put{6}  at 3.5 0.5    %
\put{7}  at 1.5 3.5    %
\put{8}  at 3.5 1.5    %
\put{9}  at 4.5 0.5    %
\put{10} at 0.5 5.5    %
\put{11} at 1.5 4.5    %
\put{12} at 4.5 1.5    %
\putrule from 0 6 to 1 6          %
\putrule from 0 5 to 2 5          %
\putrule from 0 4 to 2 4          %
\putrule from 0 3 to 2 3          %
\putrule from 3 2 to 5 2          %
\putrule from 0 2 to 2 2          
\putrule from 1 1 to 5 1          %
\putrule from 2 0 to 5 0          %
\putrule from 0 2 to 0 4        %
\putrule from 0 5 to 0 6        %
\putrule from 1 1 to 1 6        %
\putrule from 2 0 to 2 5        
\putrule from 3 0 to 3 2        %
\putrule from 4 0 to 4 2        %
\putrule from 5 0 to 5 2        %
\endpicture 
&\qquad
&\beginpicture
\setcoordinatesystem units <0.5cm,0.5cm> point at -1 0       
\setplotarea x from 0 to 5, y from -1 to 6    
\linethickness=0.5pt                          
\put{2}  at 0.5 2.5    %
\put{9}  at 1.5 1.5    %
\put{10}  at 2.5 0.5    %
\put{1}  at 0.5 3.5    %
\put{6}  at 1.5 2.5    
\put{11}  at 3.5 0.5    %
\put{5}  at 1.5 3.5    %
\put{7}  at 3.5 1.5    %
\put{12}  at 4.5 0.5    %
\put{3} at 0.5 5.5    %
\put{4} at 1.5 4.5    %
\put{8} at 4.5 1.5    %
\putrule from 0 6 to 1 6          %
\putrule from 0 5 to 2 5          %
\putrule from 0 4 to 2 4          %
\putrule from 0 3 to 2 3          %
\putrule from 3 2 to 5 2          %
\putrule from 0 2 to 2 2          
\putrule from 1 1 to 5 1          %
\putrule from 2 0 to 5 0          %
\putrule from 0 2 to 0 4        %
\putrule from 0 5 to 0 6        %
\putrule from 1 1 to 1 6        %
\putrule from 2 0 to 2 5        
\putrule from 3 0 to 3 2        %
\putrule from 4 0 to 4 2        %
\putrule from 5 0 to 5 2        %
\endpicture 
\cr
\hbox{contents of boxes}
&&\hbox{numbering of boxes}
&&\hbox{a standard tableau} \cr
}
$$
\endexample

\subsection{Books of placed configurations.}

The general case, when
$\gamma=(\gamma_1,\ldots,\gamma_n)$ is an arbitrary element of
$\RR^n$ and $J\subseteq P(\gamma)$,  is handled as follows.
First group the entries of $\gamma$ according to their $\ZZ$-coset in $\RR$.
Each group of entries in $\gamma$ can be arranged to form a sequence
$$
\beta+C_\beta = \beta+(z_1,\ldots,z_k) = (\beta+z_1,\ldots,\beta+z_k),
\quad\hbox{
where $0\leq\beta<1$, $z_i\in\ZZ$ and $z_1\leq \cdots
\leq z_k$.}
$$
Fix some ordering of these groups and let 
$$\vec{\gamma}=(\beta_1+C_{\beta_1},\ldots,\beta_r+C_{\beta_r})$$
be the rearrangement of the sequence $\gamma$ with the groups listed in
order.  Since $\vec \gamma$ and $\gamma$ are in the same
orbit it is sufficient to analyze $\vec \gamma$ ($\gamma$ corresponds
to the central character of the corresponding affine Hecke algebra
representations and thus any convenient element of the orbit
is appropriate, see (2.3)).

The decomposition of $\vec \gamma$ into
groups induces decompositions
$$
Z(\vec{\gamma})=\bigcup_{\beta_i}Z_{\beta_i}, \quad\quad
P(\vec{\gamma})=\bigcup_{\beta_i}P_{\beta_i}, \quad\hbox{
and, if $J\subseteq P(\vec{\gamma})$, then }\quad
J=\bigcup_{\beta_i}J_{\beta_i},$$
where $J_{\beta_i}=J\cap P_{\beta_i}$.
Each pair $(C_\beta,J_{\beta})$ is a placed shape of the type considered
in the previous subsection and we may
identify $(\vec{\gamma},J)$ with the {\it book of placed shapes}
$\left( (C_{\beta_1},J_{\beta_1}),\ldots,(C_{\beta_r},J_{\beta_r}) \right)$.
We think of this as a {\it book} with {\it pages} numbered by the values
$\beta_1,\ldots,\beta_r$ and with the placed configuration determined by
$(C_{\beta_i},J_{\beta_i})$ on page $\beta_i$.
In this form the {\it standard tableaux} of shape $(\vec \gamma,J)$ are
fillings of the $n$ boxes in the book with the numbers $1,\ldots,n$ such
that the filling on each page satisfies the conditions for a standard
tableau in (5.10).
 
\bigskip\noindent
{\sl Example.}   
If
$\gamma=(1/2,1/2,1,1,1,3/2,-2,-2,-1/2,-1,-1,-1,-1/2,1/2,0,0,0)$ then one
possibility for $\vec{\gamma}$ is
$$
\vec{\gamma} = (-2,-2,-1,-1,-1,0,0,0,1,1,1,-1/2,-1/2,1/2,1/2,1/2,3/2).$$
In this case $\beta_1=0$, $\beta_2=1/2$,
$$\beta_1+C_{\beta_1}=(-2,-2,-1,-1,-1,0,0,0,1,1,1)\quad\hbox{and}\quad
\beta_2+C_{\beta_2}=(-1/2,-1/2,1/2,1/2,1/2,3/2).$$
If $J=J_{\beta_1}\cup J_{\beta_2}$ where 
$J_{\beta_2} = \{ \varepsilon_{14}-\varepsilon_{13}, 
\varepsilon_{17}-\varepsilon_{16} \}$ and 
$$
J_{\beta_1} = \{ \varepsilon_3-\varepsilon_2, \varepsilon_4-\varepsilon_2,
\varepsilon_5-\varepsilon_2, \varepsilon_6-\varepsilon_3,
\varepsilon_6-\varepsilon_4, \varepsilon_6-\varepsilon_5,
\varepsilon_9-\varepsilon_7, \varepsilon_9-\varepsilon_8,
\varepsilon_{10}-\varepsilon_7, \varepsilon_{10}-\varepsilon_8 \} 
$$
then the book of shapes is
$$
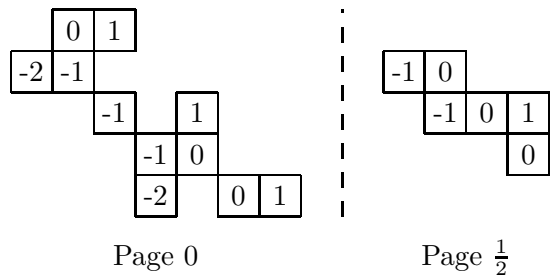
\beginpicture
\setcoordinatesystem units <0.5cm,0.5cm>         
\setplotarea x from 0 to 13, y from -1 to 5   
\linethickness=0.5pt                          
\multiput{-1} at 1.5 3.5 *2 1 -1 /
\put{0}  at 1.5 4.5    %
\put{0}  at 4.5 1.5    %
\put{0}  at 5.5 0.5    
\put{-2} at 0.5 3.5    %
\put{-2} at 3.5 0.5    %
\put{1}  at 2.5 4.5    %
\put{1}  at 4.5 2.5    %
\put{1}  at 6.5 0.5    %
\putrule from 3 0 to 4 0          %
\putrule from 5 0 to 7 0          
\putrule from 3 1 to 7 1          
\putrule from 2 2 to 5 2          %
\putrule from 0 3 to 3 3          %
\putrule from 4 3 to 5 3          %
\putrule from 0 4 to 3 4          %
\putrule from 1 5 to 3 5          %
\putrule from 0 3 to 0 4        %
\putrule from 1 3 to 1 5        
\putrule from 2 2 to 2 5        
\putrule from 3 0 to 3 3        %
\putrule from 3 4 to 3 5        %
\putrule from 4 0 to 4 3        %
\putrule from 5 0 to 5 3        %
\putrule from 6 0 to 6 1        %
\putrule from 7 0 to 7 1        %
\put{0}  at 10.5 3.5    %
\put{0}  at 11.5 2.5    %
\put{0}  at 12.5 1.5    
\put{-1} at 9.5  3.5    %
\put{-1} at 10.5 2.5    %
\put{1}  at 12.5 2.5    %
\putrule from 12 1 to 13 1          %
\putrule from 10 2 to 13 2          
\putrule from 9  3 to 13 3          
\putrule from 9  4 to 11 4          %
\putrule from 9  3 to 9  4        %
\putrule from 10 2 to 10 4        
\putrule from 11 2 to 11 4        
\putrule from 12 1 to 12 3        %
\putrule from 13 1 to 13 3        %
\setdashes
\putrule from 8 0 to 8 5           
\put{Page 0}   at 3.5 -1.0      %
\put{Page ${1\over 2}$} at 11 -1.0       
\endpicture
$$
where the numbers in the boxes are the contents of the boxes.
The filling 
$$
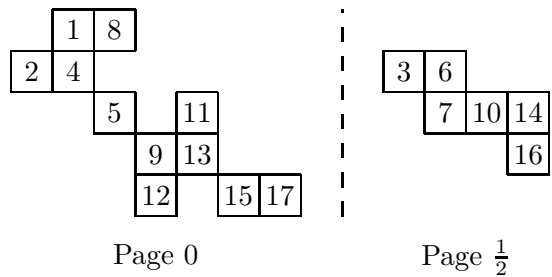
\beginpicture
\setcoordinatesystem units <0.5cm,0.5cm>         
\setplotarea x from 0 to 13, y from -1 to 5   
\linethickness=0.5pt                          
\put{4}  at 1.5 3.5    %
\put{5}  at 2.5 2.5    %
\put{9}  at 3.5 1.5    %
\put{2}  at 0.5 3.5    %
\put{12} at 3.5 0.5    %
\put{1}  at 1.5 4.5    %
\put{13} at 4.5 1.5    %
\put{15} at 5.5 0.5    
\put{8}  at 2.5 4.5    %
\put{11} at 4.5 2.5    %
\put{17} at 6.5 0.5    %
\putrule from 3 0 to 4 0          %
\putrule from 5 0 to 7 0          
\putrule from 3 1 to 7 1          
\putrule from 2 2 to 5 2          %
\putrule from 0 3 to 3 3          %
\putrule from 4 3 to 5 3          %
\putrule from 0 4 to 3 4          %
\putrule from 1 5 to 3 5          %
\putrule from 0 3 to 0 4        %
\putrule from 1 3 to 1 5        
\putrule from 2 2 to 2 5        
\putrule from 3 0 to 3 3        %
\putrule from 3 4 to 3 5        %
\putrule from 4 0 to 4 3        %
\putrule from 5 0 to 5 3        %
\putrule from 6 0 to 6 1        %
\putrule from 7 0 to 7 1        %
\put{6}  at 10.5 3.5    %
\put{10} at 11.5 2.5    %
\put{16} at 12.5 1.5    
\put{3}  at 9.5  3.5    %
\put{7}  at 10.5 2.5    %
\put{14} at 12.5 2.5    %
\putrule from 12 1 to 13 1          %
\putrule from 10 2 to 13 2          
\putrule from 9  3 to 13 3          
\putrule from 9  4 to 11 4          %
\putrule from 9  3 to 9  4        %
\putrule from 10 2 to 10 4        
\putrule from 11 2 to 11 4        
\putrule from 12 1 to 12 3        %
\putrule from 13 1 to 13 3        %
\setdashes
\putrule from 8 0 to 8 5           
\put{Page 0}   at 3.5 -1.0      %
\put{Page ${1\over2}$} at 11 -1.0       
\endpicture
$$
is a standard tableau of shape $(\vec{\gamma},J)$.  This filling corresponds to
the permutation 
$$w=(2,12,4,5,9,1,13,15,8,11,17,3,7,6,10,16,14)
\qquad\hbox{in $\cF^{(\vec\gamma,J)}\subseteq S_{16}$.}
\qquad\qquad\hbox{\qed}$$


\section 6.  Skew shapes, ribbons, conjugation, etc. in Type A

In this section we shall explain how the definitions in 
Section 5.1 correspond to classical notions in 
Young tableaux theory.
As in the previous section let $R$ be the root system of Type
$A_{n-1}$ as given in (5.2).
For clarity, we shall state all of the results in this section
for placed shapes $(\gamma,J)$ such that $\gamma$ is dominant
and integral, i.e. $\gamma=(\gamma_1,\ldots,\gamma_n)$
with $\gamma_1\le \cdots\le \gamma_n$ and $\gamma_i\in \ZZ$.  
This assumption is purely for notational clarity.

\subsection{Axial distance.}

Let $(\gamma,J)$ be a local region such that $\gamma$ is dominant
and integral.  Let $w_p\in \cF^{(\gamma,J)}$ and let $p$ be the corresponding
standard tableau as defined by the map in Theorem 5.9.  
Then it follows from the definitions of $\gamma$ and $w_p$ in 
(5.6) and (5.4) that
$$\langle w\gamma,\varepsilon_i\rangle
=\langle \gamma, w_p^{-1}\varepsilon_i\rangle
=c({\rm box}_{w_p^{-1}(i)})=c(p(i)),
\formula$$
where $p(i)$ is the box of $p$ containing the entry $i$.

In classical standard tableau theory the {\it axial distance} between
two boxes in a standard tableau is defined as follows.
Let $\lambda$ be a partition and let $p$ be a standard tableau of 
shape $\lambda$.  Let $1\le i,j\le n$ and let $p(i)$ and $p(j)$
be the boxes which are filled with $i$ and $j$ respectively.
Let $(r_i,c_i)$ and $(r_j,c_j)$ be the positions of these boxes,
where the rows and columns of $\lambda$ are numbered in the same way as for 
matrices.  Then the {\it axial distance} from $j$ to $i$ in $p$ is
$$d_{ji}(p) = c_j-c_i+r_i-r_j,$$
(see [Wz]).  Rewriting this in terms of the local region
$(\gamma,J)$ determined by (5.7)
$$d_{ji}(p) = c(p(j))-c(p(i)) 
=\langle w_p\gamma, \varepsilon_j-\varepsilon_i\rangle
=d_{\varepsilon_j-\varepsilon_i}(w),$$
where $w_p\in \cF^{(\gamma,J)}$ is the permutation corresponding to the 
standard tableau $p$ and $d_\alpha(w_p)$ is the
$\alpha$-axial distance
defined in (5.1f).  This shows that the axial distance defined in
(5.1f) is a generalization of the classical notion of axial distance.
These numbers are crucial to the classical construction of 
the seminormal representations of the symmetric group given by
Young (see Remark (3) of Section 5.1).

\subsection{Skew shapes.}

The following proposition shows that, in the case of a root system
of type $A$, the definition of skew local region coincides with 
the classical notion of a skew shape.  

\prop  
Let $(\gamma,J)$ be a local region with $\gamma$ 
dominant and integral. Then the configuration of boxes 
associated to $(\gamma,J)$ is a placed
skew shape if and only if $(\gamma,J)$ is a skew local region.
\endprop
\pf
$\Longleftarrow:$  We shall show that if the placed configuration 
corresponding to the pair $(\gamma,J)$ has any $2\times 2$ blocks of the 
forms
$$
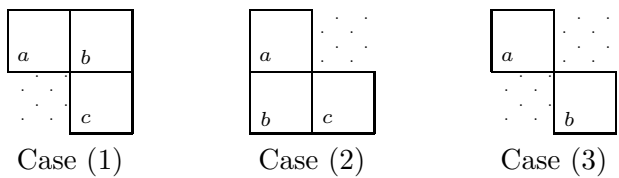
\matrix{
\beginpicture
\setcoordinatesystem units <0.75cm,0.75cm>         
\setplotarea x from 0 to 2, y from 0 to 2    
\linethickness=0.5pt                          
\put{$\scriptstyle{a}$} at 0.25 1.25
\put{$\scriptstyle{b}$}  at 1.25 1.25
\put{$\scriptstyle{c}$} at 1.25 0.25
\putrule from 0 2 to 2 2          %
\putrule from 0 1 to 2 1          
\putrule from 1 0 to 2 0          %
\putrule from 0 1 to 0 2        %
\putrule from 1 0 to 1 2        %
\putrule from 2 0 to 2 2        
\vshade 0 0 1    1 0 1  /           
\endpicture
&\qquad
&\beginpicture
\setcoordinatesystem units <0.75cm,0.75cm>         
\setplotarea x from 0 to 2, y from 0 to 2    
\linethickness=0.5pt                          
\put{$\scriptstyle{a}$} at 0.25 1.25
\put{$\scriptstyle{b}$}  at 0.25 0.25
\put{$\scriptstyle{c}$} at 1.25 0.25
\putrule from 0 2 to 1 2          %
\putrule from 0 1 to 2 1          
\putrule from 0 0 to 2 0          %
\putrule from 0 0 to 0 2        %
\putrule from 1 0 to 1 2        %
\putrule from 2 0 to 2 1        
\vshade 1 1 2    2 1 2  /           
\endpicture
&\qquad
&\beginpicture
\setcoordinatesystem units <0.75cm,0.75cm>         
\setplotarea x from 0 to 2, y from 0 to 2    
\linethickness=0.5pt                          
\put{$\scriptstyle{a}$} at 0.25 1.25
\put{$\scriptstyle{b}$}  at 1.25 0.25
\putrule from 0 2 to 1 2          %
\putrule from 0 1 to 2 1          
\putrule from 1 0 to 2 0          %
\putrule from 0 1 to 0 2        %
\putrule from 1 0 to 1 2        %
\putrule from 2 0 to 2 1        
\vshade 0 0 1    1 0 1  /           
\vshade 1 1 2    2 1 2  /           
\endpicture
\cr
\hbox{Case (1)} 
&&\hbox{Case (2)}
&&\hbox{Case (3)} \cr
}
$$
then there exists a $w\in {\cal F}^{(\gamma,J)}$ such that
$w\gamma$ violates one of the two conditions in (5.1c).
This will show that if $(\gamma,J)$ is a skew local region then
the corresponding placed configuration of boxes must be 
a placed skew shape.
In the pictures above the shaded regions indicate the
absence of a box and, for reference, we have labeled the 
boxes with $a,b,c$.

\medskip\noindent
{\it Case} (1):  Create a standard
tableau $p$ such that the $2\times 2$ block is filled with
$$
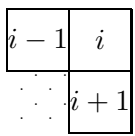
\beginpicture
\setcoordinatesystem units <0.75cm,0.75cm>         
\setplotarea x from 0 to 2, y from 0 to 2    
\linethickness=0.5pt                          
\put{$i-1$} at 0.5 1.5
\put{$i$}  at 1.5 1.5
\put{$i+1$} at 1.5 0.5
\putrule from 0 2 to 2 2          %
\putrule from 0 1 to 2 1          
\putrule from 1 0 to 2 0          %
\putrule from 0 1 to 0 2        %
\putrule from 1 0 to 1 2        %
\putrule from 2 0 to 2 2        
\vshade 0 0 1    1 0 1  /           
\endpicture
$$
by filling the region of the configuration strictly north and weakly 
west of box c in row reading order (sequentially left to right across
the rows starting at the top), putting the next entry in box c, 
and filling the remainder of the
configuration in column reading order (sequentially down the columns
beginning at the leftmost available column).  
Let $w=w_p$ be the permutation in $\cF^{(\gamma,J)}$
which corresponds to the standard tableau $p$.  Let $p(i)$ denote the
box containing $i$ in $p$. Then, using the identity (6.2),
$$\langle w\gamma,\alpha_i+\alpha_{i+1}\rangle
=\langle w\gamma,\varepsilon_{i+1}-\varepsilon_{i-1}\rangle
=c(p(i+1))-c(p(i-1)) = 0,$$
since the boxes $p(i+1)$ and $p(i-1)$ are on the same diagonal.
However, 
$$\eqalign{
\langle w\gamma,\alpha_i\rangle
&= \langle w\gamma, \varepsilon_i-\varepsilon_{i-1}\rangle
=c(p(i))-c(p(i-1))=1, \quad\hbox{and}\cr
\langle w\gamma,\alpha_{i+1}\rangle
&= \langle w\gamma, \varepsilon_{i+1}-\varepsilon_i\rangle
=c(p(i+1))-c(p(i))=-1, \cr
}$$
and so condition (2) in (5.1c) is violated.

\medskip\noindent
{\it Case} (2):   
Create a standard tableau $p$ such that the $2\times 2$ block is
filled with
$$
\beginpicture
\setcoordinatesystem units <0.75cm,0.75cm>         
\setplotarea x from 0 to 2, y from 0 to 2    
\linethickness=0.5pt                          
\put{$i-1$} at 0.5 1.5
\put{$i$}  at 0.5 0.5
\put{$i+1$} at 1.5 0.5
\putrule from 0 2 to 1 2          %
\putrule from 0 1 to 2 1          
\putrule from 0 0 to 2 0          %
\putrule from 0 0 to 0 2        %
\putrule from 1 0 to 1 2        %
\putrule from 2 0 to 2 1        
\vshade 1 1 2    2 1 2  /           
\endpicture
$$
by filling the region weakly north and strictly west of box c in column 
reading order, putting the next entry in box c, and filling the remainder 
of the configuration in row reading order.  Using this standard tableau $p$,
the remainder of the argument is the same as for case (1).

\medskip\noindent
{\it Case} (3): 
Create a standard tableau $p$ such that the $2\times 2$ block is
filled with
$$
\beginpicture
\setcoordinatesystem units <0.75cm,0.75cm>         
\setplotarea x from 0 to 2, y from 0 to 2    
\linethickness=0.5pt                          
\put{$i-1$} at 0.5 1.5
\put{$i$}  at 1.5 0.5
\putrule from 0 2 to 1 2          %
\putrule from 0 1 to 2 1          
\putrule from 1 0 to 2 0          %
\putrule from 0 1 to 0 2        %
\putrule from 1 0 to 1 2        %
\putrule from 2 0 to 2 1        
\vshade 0 0 1    1 0 1  /           
\vshade 1 1 2    2 1 2  /           
\endpicture
$$
by filling the region strictly north and strictly west of box b in 
column reading order, putting the next entry in box b, and 
filling the remainder of the
configuration in row reading order.  Let $w=w_p$ be the permutation
in $\cF^{(\gamma,J)}$ corresponding to $p$ and let $p(i)$ denote the
box containing $i$ in $p$. Then 
$$\langle w\gamma,\alpha_i\rangle
=\langle w\gamma,\varepsilon_i-\varepsilon_{i-1}\rangle
=c(p(i))-c(p(i-1)) = 0,$$
since $t(i)$ and $t(i-1)$ are on the same diagonal.
Hence, condition (1) in (5.1c) is violated.

$\Longrightarrow:$  
Let $\gamma\in \ZZ^n$ and $\lambda/\mu$ describe a placed skew
shape (a skew shape placed on infinite graph paper).
Let $(\gamma,J)$ be the corresponding local region as 
defined in (5.7).  We will show that every $w\gamma$ is calibratable
for every $w\in {\cal F}^{(\gamma,J)}$.

Let $w\in \cF^{(\gamma,J)}$ and let $p$ be the corresponding
standard tableau of shape $\lambda/\mu$.  Consider a 
$2\times 2$ block of boxes of $p$.
If these boxes are filled with
$$
\beginpicture
\setcoordinatesystem units <0.75cm,0.75cm>         
\setplotarea x from 0 to 4, y from 0 to 2    
\linethickness=0.5pt                          
\put{$i$} at 0.5 1.5
\put{$j$}  at 1.5 1.5
\put{$k$}  at 0.5 0.5
\put{$\ell$} at 1.5 0.5
\putrule from 0 2 to 2 2          %
\putrule from 0 1 to 2 1          
\putrule from 0 0 to 2 0          %
\putrule from 0 0 to 0 2        %
\putrule from 1 0 to 1 2        %
\putrule from 2 0 to 2 2        
\endpicture
$$
then either $i<j<k<\ell$ or $i<k<j<\ell$.  In both cases we have
$i<\ell-1$ and it follows that $\ell-1$ and $\ell$ are not on the same diagonal.
Thus
$$\langle w\gamma,\alpha_\ell\rangle=c(p(\ell))-c(p(\ell-1))\ne 0,$$
and so $w\gamma$ stasfies condition (a) in the definition of
calibratable.

The same argument shows that one can never get a standard tableau in which
$\ell$ and $\ell-2$ occur in adjacent boxes of the same diagonal and thus
it follows that $w\gamma$ satisfies condition (b)
in the definition of calibratable.  
Thus $(\gamma,J)$ is a skew local region.
\endpf

\subsection{Ribbon Shapes.}

Classically, a {\it border strip} (or {\it ribbon}) is a skew
shape which contains at most one box in each diagonal.  Although the 
convention,
[Mac, I \S 1 p. 5], is to assume that border strips are connected skew
shapes we shall {\it not} assume this.

Recall from (5.1b) that a placed shape $(\gamma,J)$ is a
placed {\it ribbon} shape if $\gamma$ is regular, 
i.e. $\langle \gamma,\alpha\rangle \ne 0$ for all $\alpha\in R$.

\prop 
Let $(\gamma,J)$ be a placed ribbon shape such that $\gamma$
is dominant and integral. 
Then the configuration of boxes coresponding
to $(\gamma,J)$ is a placed border strip.
\pf
Let $(\gamma,J)$ be a placed ribbon shape with $\gamma$ dominant
and regular.  Since $\gamma=(\gamma_1,\ldots,\gamma_n)$ is regular,
$\gamma_i\ne \gamma_j$ for all $i\ne j$.  In terms of the placed 
configuration $\gamma_i=c({\rm box}_i)$ is the diagonal that ${\rm box}_i$ 
is on.  Thus the configuration of boxes corresponding to $(\gamma,J)$ contains
at most one box in each diagonal.
\endpf

\noindent
{\sl Example.}
If $\gamma=(-6,-5,-4,0,1,3,4,5,6,7)$ and
$J=\{\varepsilon_2-\varepsilon_1,
\varepsilon_5-\varepsilon_4,
\varepsilon_7-\varepsilon_6,
\varepsilon_9-\varepsilon_8,
\varepsilon_{10}-\varepsilon_9\}$
then the placed configuration of boxes corresponding to $(\gamma,J)$
is the placed border strip
$$
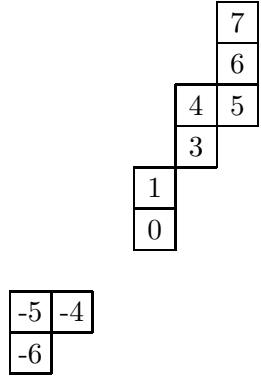
\beginpicture
\setcoordinatesystem units <0.5cm,0.5cm> point at 6 0       
\setplotarea x from 0 to 6, y from -1 to 9    
\linethickness=0.5pt                          
\put{-6}  at 0.5 0.5    %
\put{-5}  at 0.5 1.5    %
\put{-4}  at 1.5 1.5    %
\put{0}  at 3.5 3.5    %
\put{1}  at 3.5 4.5    
\put{3}  at 4.5 5.5    %
\put{4}  at 4.5 6.5    %
\put{5}  at 5.5 6.5    %
\put{6}  at 5.5 7.5    %
\put{7}  at 5.5 8.5    %
\putrule from 0 0 to 1 0          %
\putrule from 0 1 to 2 1          %
\putrule from 0 2 to 2 2          %
\putrule from 3 3 to 4 3          %
\putrule from 3 4 to 4 4          %
\putrule from 3 5 to 5 5          %
\putrule from 4 6 to 6 6          %
\putrule from 4 7 to 6 7          %
\putrule from 5 8 to 6 8          %
\putrule from 5 9 to 6 9          %
\putrule from 0 0 to 0 2        %
\putrule from 1 0 to 1 2        %
\putrule from 2 1 to 2 2        %
\putrule from 3 3 to 3 5        
\putrule from 4 3 to 4 7        %
\putrule from 5 5 to 5 9        %
\putrule from 6 6 to 6 9        %
\endpicture
$$
where the boxes are labeled with their contents.
\endpf

\subsection{Conjugation of Shapes.}

Let $(\gamma,J)$ be a placed shape with $\gamma$ dominant and integral
(i.e. $\gamma=(\gamma_1,\ldots,\gamma_n)$ with 
$\gamma_1\leq\cdots\leq\gamma_n$ and $\gamma_i\in\ZZ$)
and view $(\gamma,J)$ as a placed configuration of boxes.  
In terms of placed configurations,
conjugation of shapes is equivalent to transposing
the placed configuration across the diagonal of boxes of content $0$.
The following example illustrates this.

\bigskip\noindent
{\sl Example.}
Suppose $\gamma=(-1,-1,-1,0,0,1,1)$ and $J=(\varepsilon_4-\varepsilon_2,
\varepsilon_4-\varepsilon_3, \varepsilon_6-\varepsilon_5,
\varepsilon_7-\varepsilon_5)$. Then the placed configuration of boxes
corresponding to  $(\gamma,J)$ is 
$$
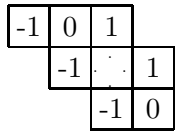
\beginpicture
\setcoordinatesystem units <0.5cm,0.5cm>         
\setplotarea x from 0 to 4, y from 0 to 3    
\linethickness=0.5pt                          
\put{-1} at 0.5 2.5
\put{0}  at 1.5 2.5
\put{-1} at 1.5 1.5
\put{-1} at 2.5 0.5
\put{1}  at 2.5 2.5
\put{1}  at 3.5 1.5
\put{0}  at 3.5 0.5
\putrule from 0 3 to 3 3          %
\putrule from 0 2 to 4 2          
\putrule from 1 1 to 4 1          %
\putrule from 2 0 to 4 0          %
\putrule from 0 2 to 0 3        %
\putrule from 1 1 to 1 3        %
\putrule from 2 0 to 2 3        
\putrule from 3 0 to 3 3        %
\putrule from 4 0 to 4 2        %
\vshade 2 1 2   3 1 2 /           
\endpicture
$$
in which the shaded box is not a box in the configuration.

The minimal length representative of the coset $w_0W_\gamma$ is
the permutation
$$u=\left(
\matrix{
1&2&3&4&5&6&7\cr
5&6&7&3&4&1&2 \cr
}
\right).$$ 
We have $-u\gamma=-w_0\gamma=(-1,-1,0,0,1,1,1)$ and 
$$
\eqalign{
-u(P(\gamma)\setminus J) &= 
-u\left\{ \varepsilon_4-\varepsilon_1,
\varepsilon_5-\varepsilon_1, \varepsilon_5-\varepsilon_2,
\varepsilon_5-\varepsilon_3, \varepsilon_6-\varepsilon_4,
\varepsilon_7-\varepsilon_4 \right\} \cr
&= -\left\{ \varepsilon_3-\varepsilon_5,
\varepsilon_4-\varepsilon_5, \varepsilon_4-\varepsilon_6,
\varepsilon_4-\varepsilon_7, \varepsilon_1-\varepsilon_3,
\varepsilon_2-\varepsilon_3 \right\} \cr
&= \left\{ \varepsilon_5-\varepsilon_3,
\varepsilon_5-\varepsilon_4, \varepsilon_6-\varepsilon_4,
\varepsilon_7-\varepsilon_4, \varepsilon_3-\varepsilon_1,
\varepsilon_3-\varepsilon_2 \right\}. \cr
}
$$
Thus the configuration of boxes corresponding to the placed shape 
$(\gamma,J)'$ is 
$$
\beginpicture
\setcoordinatesystem units <0.5cm,0.5cm>         
\setplotarea x from 0 to 3, y from 0 to 4    
\linethickness=0.5pt                          
\put{0}  at 2.5 0.5 
\put{1}  at 2.5 1.5 
\put{-1} at 1.5 0.5 
\put{1}  at 1.5 2.5 
\put{-1} at 0.5 1.5 
\put{0}  at 0.5 2.5 
\put{1}  at 0.5 3.5 
\putrule from 3 0 to 3 2          %
\putrule from 2 0 to 2 3          
\putrule from 1 0 to 1 4          %
\putrule from 0 1 to 0 4          %
\putrule from 1 0 to 3 0        %
\putrule from 0 1 to 3 1        %
\putrule from 0 2 to 3 2        
\putrule from 0 3 to 2 3        %
\putrule from 0 4 to 1 4        %
\vshade 1 1 2   2 1 2 /           
\endpicture
$$
\endexample

\subsection{Row reading and column reading tableaux.}

Let $(\gamma,J)$ be a placed shape such that $\gamma$ is 
dominant and integral and consider the placed configuration of 
boxes corresponding to $(\gamma,J)$.  
The {\it minimal box} of the configuration is the box such that
\smallskip
\itemitem{($m_1$)} there is no box immediately above,
\smallskip
\itemitem{($m_2$)} there is no box immediately to the left, 
\smallskip
\itemitem{($m_3$)} there is no box northwest in the same diagonal, and
\smallskip
\itemitem{($m_4$)} it has the minimal content of the boxes satisfying 
($m_1$), ($m_2$) and ($m_3$).
\smallskip\noindent
There is at most one box in each diagonal satisfying ($m_1$), ($m_2$),
and ($m_3$).  Thus, ($m_4$) guarantees that the minimal box is unique.
It is clear that the minimal box of the configuration
always exists.

The {\it column reading} tableaux of shape $(\gamma,J)$ is
the filling $p_{\rm min}$ which is created inductively by 
\smallskip
\item{(a)} filling the minimal box of the configuration with $1$, and
\smallskip
\item{(b)} if $1,2,\ldots, i$ have been filled in then fill the minimal 
box of the configuration formed by the unfilled boxes with $i+1$.
\smallskip\noindent
The {\it row reading tableau} of shape $(\gamma,J)$ is the
standard tableau $p_{\rm max}$ whose conjugate 
$(p_{\rm max})'$ is the column reading tableaux
for the shape $(\gamma,J)'$ (the conjugate shape to $(\gamma,J)$).

Recall the definitions of the weak Bruhat order and
closed subsets of roots given after equation (4.5).

\thm  Let $(\gamma,J)$ be a placed shape such that $\gamma$ is 
dominant and integral
(i.e. $\gamma=(\gamma_1,\ldots,\gamma_n)$ with 
$\gamma_1\leq\cdots\leq\gamma_n$ and $\gamma_i\in\ZZ$).
Let $p_{\rm min}$ and $p_{\rm max}$
be the column reading and row reading tableaux of shape $(\gamma,J)$,
respectively, and let $w_{\rm min}$ and $w_{\rm max}$ be the corresponding
permutations in $\cF^{(\gamma,J)}$.
Then 
$$R(w_{\rm min})=\overline{J},
\qquad
R(w_{\rm max})= \overline{(P(\gamma)\setminus J)\cup Z(\gamma)}^c,
\qquad\hbox{and}\qquad
\cF^{(\gamma,J)} = 
[w_{min},w_{max}],$$
where $K^c$ denotes the complement of $K$ in $R^+$ and $[w_{min},w_{max}]$ denotes
the interval between $w_{\rm min}$ and $\,w_{\rm max}$ in the weak Bruhat order.
\pf
(a)  Consider the configuration of boxes corresponding to $(\gamma,J)$.
If $k>i$ then either
$c({\rm box}_k)>c({\rm box}_i)$, or
${\rm box}_k$ is in the same diagonal and southeast of 
${\rm box}_i$.
Thus when we create $p_{\rm min}$ we have that
$$\hbox{If $k>i$ then ${\rm box}_k$ gets filled before
${\rm box}_i$ if and only if ${\rm box}_k$ is {\sl northwest} of
${\rm box}_i$,}
$$
where the {\sl northwest} is in a very strong sense:
There is a sequence of boxes
$${\rm box}_i={\rm box}_{i_0},\enspace {\rm box}_{i_1},\enspace\ldots, 
\enspace{\rm box}_{i_r}={\rm box}_k$$
such that ${\rm box}_{i_m}$ is either directly above ${\rm box}_{i_{m-1}}$
or in the same diagonal and directly northwest of ${\rm box}_{i_{m-1}}$.
In other words,
$$
\hbox{ If $k>i$ then $p_{\rm min}({\rm box}_k)<p_{\rm min}({\rm box}_i)$
\quad$\Longleftrightarrow$\quad ${\rm box}_k$ is {\sl northwest} of ${\rm
box}_i$.}
$$
So, from the formula for $w_p$ in (5.4) we get
$$
\hbox{If $k>i$ then $w_{\rm min}(k)<w_{\rm min}(i)$
$\quad\Longleftrightarrow\quad \varepsilon_k-\varepsilon_i\in \overline{J}$, }
$$
where $w_{\rm min}$ is the permutation in $\cF^{(\gamma,J)}$
which corresponds to the filling $t_{\rm min}$  
and $\overline{J}$ is the closure of $J$ in $R$.
It follows that 
$$R(w_{\rm min})=\overline{J}.$$

(b) There are at least two ways to prove that
$R(w_{\rm max})= \overline{(P(\gamma)\setminus J)\cup Z(\gamma)}^c$.
One can mimic the proof of part (a) by defining the maximal
box of a configuration and a corresponding filling.  
Alternatively one can use the definition of conjugation 
and the fact that $R(w_0w)=R(w)^c$.
The permutation $w_{\rm min}$ is the unique minimal element
of $\cF^{(\gamma,J)}$ and the conjugate of $w_{\rm max}$ is the
unique minimal element of $\cF^{(\gamma,J)'}$.
We shall leave the details to the reader.

(c)  An element $w\in W$ 
is an element of $\cF^{(\gamma,J)}$ if and only if
$R(w)\cap P(\gamma)=J$ and $R(w)\cap Z(\gamma)=\emptyset$.
Thus $\cF^{(\gamma,J)}$ consists of those permutations
$w\in W$ such that 
$$\overline{J}\subseteq R(w) \subseteq 
\overline{(P(\gamma)\setminus J)\cup Z(\gamma)}^c.$$
Since the weak Bruhat order is the ordering determined
by inclusions of $R(w)$, it follows that $\cF^{(\gamma,J)}$ is
the interval between $w_{\rm min}$ and $w_{\rm max}$.
\endpf

\noindent
{\sl Example.}
Suppose $\gamma=(-1,-1,-1,0,0,1,1)$ and $J=\{\varepsilon_4-\varepsilon_2,
\varepsilon_4-\varepsilon_3, \varepsilon_6-\varepsilon_5,
\varepsilon_7-\varepsilon_5\}$. The minimal and maximal elements in 
$\cF^{(\gamma,J)}$ are the 
permutations
$$
w_{\rm min}=\pmatrix{
1&2&3&4&5&6&7\cr
1&3&4&2&7&5&6\cr}
\qquad\hbox{and}\qquad
w_{\rm max}=\pmatrix{
1&2&3&4&5&6&7\cr
1&5&6&2&7&3&4\cr}.
$$
The permutations correspond to the standard tableaux
$$
\beginpicture
\setcoordinatesystem units <0.5cm,0.5cm>         
\setplotarea x from 0 to 4, y from 0 to 3    
\linethickness=0.5pt                          
\put{1} at 0.5 2.5
\put{2} at 1.5 2.5
\put{3} at 1.5 1.5
\put{4} at 2.5 0.5
\put{5} at 2.5 2.5
\put{6} at 3.5 1.5
\put{7} at 3.5 0.5
\putrule from 0 3 to 3 3          %
\putrule from 0 2 to 4 2          
\putrule from 1 1 to 4 1          %
\putrule from 2 0 to 4 0          %
\putrule from 0 2 to 0 3        %
\putrule from 1 1 to 1 3        %
\putrule from 2 0 to 2 3        
\putrule from 3 0 to 3 3        %
\putrule from 4 0 to 4 2        %
\vshade 2 1 2   3 1 2 /           
\endpicture
\qquad\qquad\hbox{and}\qquad\qquad
\beginpicture
\setcoordinatesystem units <0.5cm,0.5cm>         
\setplotarea x from 0 to 4, y from 0 to 3    
\linethickness=0.5pt                          
\put{1} at 0.5 2.5
\put{2} at 1.5 2.5
\put{5} at 1.5 1.5
\put{6} at 2.5 0.5
\put{3} at 2.5 2.5
\put{4} at 3.5 1.5
\put{7} at 3.5 0.5
\putrule from 0 3 to 3 3          %
\putrule from 0 2 to 4 2          
\putrule from 1 1 to 4 1          %
\putrule from 2 0 to 4 0          %
\putrule from 0 2 to 0 3        %
\putrule from 1 1 to 1 3        %
\putrule from 2 0 to 2 3        
\putrule from 3 0 to 3 3        %
\putrule from 4 0 to 4 2        %
\vshade 2 1 2   3 1 2 /           
\endpicture
.
$$
\endexample

\section 7. The type A, root of unity case

This section describes the sets ${\cal F}^{(t,J)}$ in the
case of the root system of (5.2) when $q^2=e^{2\pi i/\ell}$,
a primitive $\ell^{\rm th}$ root of unity, $\ell>2$.

Let $t\in T$.  Identify $t$ with a sequence
$$t=(t_1,\ldots, t_n)\in \CC^n,
\qquad\hbox{where}\qquad t(X^{\varepsilon_i})=t_i.$$
For the purposes of representation theory (see Theorem 3.6)
$t$ indexes a central character (see (2.3)) and so $t$ can safely
be replaced by any element of its $W$-orbit.  In this case $W$ is 
the symmetric group, $S_n$, acting by permuting the sequence
$t=(t_1,\ldots,t_n)$.  

The cyclic group $\langle q^2\rangle$ of order $\ell$ generated by $q^2$,
acts on $\CC^*$.  Fix a choice of a set $\{\xi\}$ of 
coset representatives of the $\langle q^2\rangle$ cosets in $\CC^*$.
Replace $t$ with the sequence obtained by rearranging
its entries to group entries in the same $\langle q^2\rangle$-orbit, 
so that
$$t=(\xi_1t^{(1)},\cdots, \xi_kt^{(k)}),
\qquad\hbox{where}\quad
\hbox{$\xi_1,\ldots, \xi_k$ are distinct representatives
of the cosets in $\CC^*/\langle q^2\rangle$},$$
and each $t^{(j)}$ is a sequence of the form
$$t^{(j)}=(q^{2\gamma_1},\ldots,q^{2\gamma_r}),
\qquad\hbox{with}\quad
\gamma_1,\ldots,\gamma_r\in \{0,1,\ldots, \ell-1\}
\quad\hbox{and}\quad
\gamma_1\le \cdots\le \gamma_r.$$
As in (5.11)
this decomposition of $t$ into groups induces decompositions
$$Z(t) = \bigcup_{j=1}^k Z_{\xi_j}(t)
\qquad\hbox{and}\qquad
P(t) = \bigcup_{j=1}^k P_{\xi_j}(t),$$
and it is sufficient to analyze the case when $t$ consists
of only one group, i.e.\ all the entries of $t$ are in the
same $\langle q^2\rangle$ coset.

Now assume that 
$$t=(q^{2\gamma_1},\ldots, q^{2\gamma_n}),
\qquad\hbox{with}\qquad
\gamma_1\le \cdots\le\gamma_n, \ \ \gamma_i\in \{0,\ldots,\ell-1\}.$$
Consider a page of graph paper with diagonals labeled 
by $\ldots, 0,1,\ldots,\ell-1,0,1,\ldots,\ell-1,0,1,\ldots$
from southwest to northeast.  For each local region
$(t,J)$, $J\subseteq P(t)$, we will construct an $\ell$-periodic
configuration of boxes for which the $\ell$-periodic standard
tableaux defined below will be in bijection
with the elements of ${\cal F}^{(t,J)}$.
For each $1\le i\le n$, 
the configuration will have a box numbered $i$, ${\rm box}_i$,
on each diagonal which is labeled $\gamma_i$.
There are an infinite number of such diagonals containing a
box numbered $i$, since
the diagonals are labeled in an $\ell$-periodic fashion,
but each strip of consecutive diagonals labeled $0,1,\ldots, \ell-1$
will contain $n$ boxes.
The {\it content} of a box $b$ (see [Mac, I \S 1 Ex.\ 3]) is
$$c(b) = \hbox{(the diagonal number of the box $b$)}.$$
Then 
$$\eqalign{
Z(t) 
&= \{ \varepsilon_j-\varepsilon_i\ |\ i<j,\ \gamma_i=\gamma_j\}  \cr
&= \left\{ \varepsilon_j-\varepsilon_i\ |\ 
i<j,\ \hbox{${\rm box}_i$ and ${\rm box}_j$ are in the same
diagonal} \right\}, \cr
}$$
and
$$
\eqalign{
P(t) 
&= \left\{ \varepsilon_j-\varepsilon_i\ \Bigg\vert\ 
\matrix{ \hbox{$i<j$ and $\gamma_j=\gamma_i+1$, or}\hfill \cr
\hbox{$i<j$ and $\gamma_j=\ell-1$ and $\gamma_i=0$}\hfill \cr}
\right\} \cr
&= \{ \varepsilon_j-\varepsilon_i\ \vert\ 
\hbox{$i<j$ and 
${\rm box}_i$ and ${\rm box}_j$ are in adjacent diagonals} \}. \cr
}$$
We will use $J\subseteq P(t)$ to organize the relative
positions of the boxes in adjacent diagonals.  
$$\matrix{
\hbox{If $\varepsilon_j-\varepsilon_i\in J$ and
if $c({\rm box}_j)\ne \ell-1$ or $c({\rm box}_i)\ne 0$} \hfill 
&
&\hbox{place} \hfill
&
&\hbox{${\rm box}_j$ northwest of ${\rm box}_i$,} \hfill
\cr
\hbox{if $\varepsilon_j-\varepsilon_i\not\in J$ and
if $c({\rm box}_j)\ne \ell-1$ or $c({\rm box}_i)\ne 0$} \hfill 
&
&\hbox{place} \hfill
&
&\hbox{${\rm box}_j$ southeast of ${\rm box}_i$,}\hfill
\cr
\hbox{if $\varepsilon_j-\varepsilon_i\in J$ and
$c({\rm box}_j)=\ell-1$ and $c({\rm box}_i)=0$} \hfill 
&
&\hbox{place} \hfill
&
&\hbox{${\rm box}_j$ southeast of ${\rm box}_i$,} \hfill
\cr
\hbox{if $\varepsilon_j-\varepsilon_i\not\in J$ and
$c({\rm box}_j)=\ell-1$ and $c({\rm box}_i)=0$} \hfill
&
&\hbox{place} \hfill
&
&\hbox{${\rm box}_j$ northwest of ${\rm box}_i$.} \hfill
\cr
}$$
Thus, $t$ determines the number of boxes in each diagonal
and $J$ determines the relative positions of the boxes
in adjacent diagonals.  This information completely determines
the $\ell$-periodic configuration of boxes associated
to the pair $(t,J)$.

A {\it $\ell$-periodic standard tableau} is an $\ell$-periodic 
filling $p$ of the boxes with $1,2,\ldots, n$ such that
\smallskip\noindent
\item{(a)}
if $i<j$ and ${\rm box}_i$ and ${\rm box}_j$ are
in the same diagonal then $p(i)<p(j)$, 
\smallskip\noindent
\item{(b)}
if $i<j$ and ${\rm box}_i$ and ${\rm box}_j$ are in 
adjacent diagonals with ${\rm box}_j$ southwest of
${\rm box}_i$ then $p(i)<p(j)$, 
\smallskip\noindent
\item{(c)}
if $i<j$ and ${\rm box}_i$ and ${\rm box}_j$ are in 
adjacent diagonals with ${\rm box}_j$ northeast of
${\rm box}_i$ then $p(i)>p(j)$, 
\smallskip\noindent
where $p(i)$ denotes the entry in ${\rm box}_i$.
An $\ell$-periodic standard tableau $p$ corresponds to a 
permutation in $S_n$ via the correspondence
$$\matrix{
\{\hbox{standard tableaux}\}
&\longleftrightarrow 
&{\cal F}^{(t,J)} \cr
p &\longmapsto 
&\pmatrix{ 1 &2 &\cdots &n \cr p(1) &p(2) &\cdots &p(n) \cr}
\cr
}$$

\smallskip\noindent
{\bf Example.}  Suppose that $q^2=e^{2\pi i/4}$ and 
$$t=(q^0,q^0,q^0,q^0,q^2,q^2,q^2,q^4,q^4,q^6,q^6,q^6,q^6,q^6).$$
Then 
$$\eqalign{
Z(t) &= \{\varepsilon_2-\varepsilon_1,\varepsilon_3-\varepsilon_1,
\varepsilon_4-\varepsilon_1,
\varepsilon_3-\varepsilon_2,\varepsilon_4-\varepsilon_2,
\varepsilon_4-\varepsilon_3, \varepsilon_6-\varepsilon_5,
\varepsilon_7-\varepsilon_5, \ldots\} \quad\hbox{and} \cr
P(t) &= \{
\varepsilon_5-\varepsilon_1,\varepsilon_5-\varepsilon_2,
\varepsilon_5-\varepsilon_3,\varepsilon_5-\varepsilon_4,
\varepsilon_6-\varepsilon_1,\ldots,
\varepsilon_{14}-\varepsilon_9,
\varepsilon_{10}-\varepsilon_1,
\varepsilon_{10}-\varepsilon_2,\ldots,
\varepsilon_{14}-\varepsilon_4\}.
\cr}$$
If 
$$\eqalign{
J
&=\{
\varepsilon_5-\varepsilon_2,
\varepsilon_5-\varepsilon_3,
\varepsilon_5-\varepsilon_4,
\varepsilon_6-\varepsilon_3,
\varepsilon_6-\varepsilon_4,
\varepsilon_8-\varepsilon_5,
\varepsilon_8-\varepsilon_6,
\varepsilon_8-\varepsilon_7,
\varepsilon_9-\varepsilon_7,
\varepsilon_{10}-\varepsilon_9, \cr
&\phantom{J=\ }
\varepsilon_{11}-\varepsilon_9, 
\varepsilon_{12}-\varepsilon_9,
\varepsilon_{12}-\varepsilon_2,
\varepsilon_{12}-\varepsilon_3,
\varepsilon_{12}-\varepsilon_4,
\varepsilon_{13}-\varepsilon_2,
\varepsilon_{13}-\varepsilon_3,
\varepsilon_{13}-\varepsilon_4,
\varepsilon_{14}-\varepsilon_3,
\varepsilon_{14}-\varepsilon_4\} \cr
}$$
then the corresponding $\ell$-periodic
configuration of boxes and a sample $\ell$-periodic standard tableau are
$$
\matrix{
\beginpicture
\setcoordinatesystem units <0.5cm,0.5cm> point at 6 0       
\setplotarea x from 0 to 14, y from -1 to 8    
\linethickness=0.5pt                          
\put{$\cdots$}  at 0.5 4.5     %
\put{$\cdots$}  at 11.5 5     %
\put{{\bf 1}}  at 3.5 5.5     %
\put{{\bf 2}}  at 5.5 3.5     %
\put{{\bf 3}}  at 6.5 2.5     %
\put{{\bf 4}}  at 7.5 1.5     %
\put{{\bf 5}}  at 4.5 5.5     
\put{{\bf 6}}  at 6.5 3.5     %
\put{{\bf 7}}  at 8.5 1.5     %
\put{{\bf 8}}  at 4.5 6.5     %
\put{{\bf 9}}  at 8.5 2.5     %
\put{{\bf 10}}  at 5.5 6.5    %
\put{{\bf 11}}  at 6.5 5.5    %
\put{{\bf 12}}  at 8.5 3.5    %
\put{{\bf 13}}  at 9.5 2.5    %
\put{{\bf 14}}  at 10.5 1.5   %
\put{9}   at 3.5 3.5    %
\put{10}  at 0.5 7.5    %
\put{11}  at 1.5 6.5    %
\put{12}  at 3.5 4.5    %
\put{13}  at 4.5 3.5    %
\put{14}  at 5.5 2.5    %
\put{1}  at 8.5 4.5     %
\put{2}  at 10.5 2.5    %
\put{3}  at 11.5 1.5    %
\put{4}  at 12.5 0.5    %
\put{5}  at 9.5 4.5     
\put{6}  at 11.5 2.5    %
\put{7}  at 13.5 0.5    %
\put{8}  at 9.5 5.5    %
\putrule from 7 1 to 9 1          %
\putrule from 10 1 to 11 1        %
\putrule from 6 2 to 11 2         %
\putrule from 5 3 to 7 3          %
\putrule from 8 3 to 10 3         %
\putrule from 5 4 to 7 4          
\putrule from 8 4 to 9 4          %
\putrule from 3 5 to 5 5          %
\putrule from 6 5 to 7 5          %
\putrule from 3 6 to 7 6          %
\putrule from 4 7 to 6 7          %
\putrule from 3 5 to 3 6        %
\putrule from 4 5 to 4 7        %
\putrule from 5 3 to 5 4        %
\putrule from 5 5 to 5 7        
\putrule from 6 2 to 6 4        %
\putrule from 6 5 to 6 7        %
\putrule from 7 1 to 7 4        %
\putrule from 7 5 to 7 6        %
\putrule from  8 1 to 8 4       %
\putrule from 9 1 to 9 4        %
\putrule from 10 1 to 10 3      %
\putrule from 11 1 to 11 2      %
\setdashes
\putrule from 12 0 to 14 0        %
\putrule from 11 1 to 14 1        %
\putrule from 5 2 to 6 2          %
\putrule from 11 2 to 12 2        %
\putrule from 3 3 to 5 3          %
\putrule from 10 3 to 12 3        
\putrule from 3 4 to 5 4          %
\putrule from 9 4 to 10 4         %
\putrule from 8 5 to 10 5         %
\putrule from 1 6 to 2 6          %
\putrule from 9 6 to 10 6         %
\putrule from 0 7 to 2 7          %
\putrule from 0 8 to 1 8          %
\putrule from 0 7 to 0 8       %
\putrule from 1 6 to 1 8       %
\putrule from 2 6 to 2 7       %
\putrule from 3 3 to 3 5       %
\putrule from 4 3 to 4 5       
\putrule from 5 2 to 5 3       %
\putrule from 8 4 to 8 5       %
\putrule from 9 4 to 9 6       %
\putrule from 10 4 to 10 6     %
\putrule from 11 2 to 11 3     %
\putrule from 12 0 to 12 3     %
\putrule from 13 0 to 13 1     %
\putrule from 14 0 to 14 1     %
\endpicture
&
&
\beginpicture
\setcoordinatesystem units <0.5cm,0.5cm> point at 6 0       
\setplotarea x from 0 to 14, y from -1 to 8    
\linethickness=0.5pt                          
\put{$\cdots$}  at 0.5 4.5     %
\put{$\cdots$}  at 11.5 5     %
\put{{\bf 0}}  at 3.5 5.5     %
\put{{\bf 0}}  at 5.5 3.5     %
\put{{\bf 0}}  at 6.5 2.5     %
\put{{\bf 0}}  at 7.5 1.5     %
\put{{\bf 1}}  at 4.5 5.5     
\put{{\bf 1}}  at 6.5 3.5     %
\put{{\bf 1}}  at 8.5 1.5     %
\put{{\bf 2}}  at 4.5 6.5     %
\put{{\bf 2}}  at 8.5 2.5     %
\put{{\bf 3}}  at 5.5 6.5    %
\put{{\bf 3}}  at 6.5 5.5    %
\put{{\bf 3}}  at 8.5 3.5    %
\put{{\bf 3}}  at 9.5 2.5    %
\put{{\bf 3}}  at 10.5 1.5   %
\put{2}   at 3.5 3.5    %
\put{3}  at 0.5 7.5    %
\put{3}  at 1.5 6.5    %
\put{3}  at 3.5 4.5    %
\put{3}  at 4.5 3.5    %
\put{3}  at 5.5 2.5    %
\put{0}  at 8.5 4.5     %
\put{0}  at 10.5 2.5    %
\put{0}  at 11.5 1.5    %
\put{0}  at 12.5 0.5    %
\put{1}  at 9.5 4.5     
\put{1}  at 11.5 2.5    %
\put{1}  at 13.5 0.5    %
\put{2}  at 9.5 5.5    %
\putrule from 7 1 to 9 1          %
\putrule from 10 1 to 11 1        %
\putrule from 6 2 to 11 2         %
\putrule from 5 3 to 7 3          %
\putrule from 8 3 to 10 3         %
\putrule from 5 4 to 7 4          
\putrule from 8 4 to 9 4          %
\putrule from 3 5 to 5 5          %
\putrule from 6 5 to 7 5          %
\putrule from 3 6 to 7 6          %
\putrule from 4 7 to 6 7          %
\putrule from 3 5 to 3 6        %
\putrule from 4 5 to 4 7        %
\putrule from 5 3 to 5 4        %
\putrule from 5 5 to 5 7        
\putrule from 6 2 to 6 4        %
\putrule from 6 5 to 6 7        %
\putrule from 7 1 to 7 4        %
\putrule from 7 5 to 7 6        %
\putrule from  8 1 to 8 4       %
\putrule from 9 1 to 9 4        %
\putrule from 10 1 to 10 3      %
\putrule from 11 1 to 11 2      %
\setdashes
\putrule from 12 0 to 14 0        %
\putrule from 11 1 to 14 1        %
\putrule from 5 2 to 6 2          %
\putrule from 11 2 to 12 2        %
\putrule from 3 3 to 5 3          %
\putrule from 10 3 to 12 3        
\putrule from 3 4 to 5 4          %
\putrule from 9 4 to 10 4         %
\putrule from 8 5 to 10 5         %
\putrule from 1 6 to 2 6          %
\putrule from 9 6 to 10 6         %
\putrule from 0 7 to 2 7          %
\putrule from 0 8 to 1 8          %
\putrule from 0 7 to 0 8       %
\putrule from 1 6 to 1 8       %
\putrule from 2 6 to 2 7       %
\putrule from 3 3 to 3 5       %
\putrule from 4 3 to 4 5       
\putrule from 5 2 to 5 3       %
\putrule from 8 4 to 8 5       %
\putrule from 9 4 to 9 6       %
\putrule from 10 4 to 10 6     %
\putrule from 11 2 to 11 3     %
\putrule from 12 0 to 12 3     %
\putrule from 13 0 to 13 1     %
\putrule from 14 0 to 14 1     %
\endpicture
\cr
\hbox{numbering of boxes}
&&\hbox{contents of boxes} \cr
}
$$
\medskip\noindent
$$
\matrix{
\beginpicture
\setcoordinatesystem units <0.5cm,0.5cm> point at 6 0       
\setplotarea x from 0 to 14, y from -1 to 8    
\linethickness=0.5pt                          
\put{$\cdots$}  at 0.5 4.5     %
\put{$\cdots$}  at 11.5 5     %
\put{{\bf 1}}  at 3.5 5.5     %
\put{{\bf 9}}  at 5.5 3.5     %
\put{{\bf 13}}  at 6.5 2.5     %
\put{{\bf 14}}  at 7.5 1.5     %
\put{{\bf 6}}  at 4.5 5.5     
\put{{\bf 12}}  at 6.5 3.5     %
\put{{\bf 4}}  at 8.5 1.5     %
\put{{\bf 5}}  at 4.5 6.5     %
\put{{\bf 3}}  at 8.5 2.5     %
\put{{\bf 8}}  at 5.5 6.5    %
\put{{\bf 11}}  at 6.5 5.5    %
\put{{\bf 2}}  at 8.5 3.5    %
\put{{\bf 7}}  at 9.5 2.5    %
\put{{\bf 10}}  at 10.5 1.5   %
\put{3}   at 3.5 3.5    %
\put{8}  at 0.5 7.5    %
\put{11}  at 1.5 6.5    %
\put{2}  at 3.5 4.5    %
\put{7}  at 4.5 3.5    %
\put{10}  at 5.5 2.5    %
\put{1}  at 8.5 4.5     %
\put{9}  at 10.5 2.5    %
\put{13}  at 11.5 1.5    %
\put{14}  at 12.5 0.5    %
\put{6}  at 9.5 4.5     
\put{12}  at 11.5 2.5    %
\put{4}  at 13.5 0.5    %
\put{5}  at 9.5 5.5    %
\putrule from 7 1 to 9 1          %
\putrule from 10 1 to 11 1        %
\putrule from 6 2 to 11 2         %
\putrule from 5 3 to 7 3          %
\putrule from 8 3 to 10 3         %
\putrule from 5 4 to 7 4          
\putrule from 8 4 to 9 4          %
\putrule from 3 5 to 5 5          %
\putrule from 6 5 to 7 5          %
\putrule from 3 6 to 7 6          %
\putrule from 4 7 to 6 7          %
\putrule from 3 5 to 3 6        %
\putrule from 4 5 to 4 7        %
\putrule from 5 3 to 5 4        %
\putrule from 5 5 to 5 7        
\putrule from 6 2 to 6 4        %
\putrule from 6 5 to 6 7        %
\putrule from 7 1 to 7 4        %
\putrule from 7 5 to 7 6        %
\putrule from  8 1 to 8 4       %
\putrule from 9 1 to 9 4        %
\putrule from 10 1 to 10 3      %
\putrule from 11 1 to 11 2      %
\setdashes
\putrule from 12 0 to 14 0        %
\putrule from 11 1 to 14 1        %
\putrule from 5 2 to 6 2          %
\putrule from 11 2 to 12 2        %
\putrule from 3 3 to 5 3          %
\putrule from 10 3 to 12 3        
\putrule from 3 4 to 5 4          %
\putrule from 9 4 to 10 4         %
\putrule from 8 5 to 10 5         %
\putrule from 1 6 to 2 6          %
\putrule from 9 6 to 10 6         %
\putrule from 0 7 to 2 7          %
\putrule from 0 8 to 1 8          %
\putrule from 0 7 to 0 8       %
\putrule from 1 6 to 1 8       %
\putrule from 2 6 to 2 7       %
\putrule from 3 3 to 3 5       %
\putrule from 4 3 to 4 5       
\putrule from 5 2 to 5 3       %
\putrule from 8 4 to 8 5       %
\putrule from 9 4 to 9 6       %
\putrule from 10 4 to 10 6     %
\putrule from 11 2 to 11 3     %
\putrule from 12 0 to 12 3     %
\putrule from 13 0 to 13 1     %
\putrule from 14 0 to 14 1     %
\endpicture
\cr
\hbox{a standard tableau $p$} \cr
}$$

\section 8. Standard tableaux for type $C$ in terms of boxes

\subsection{The root system.}

Let $\{\varepsilon_1,\ldots,\varepsilon_n\}$ be an orthonormal basis  of
$\fh_{\RR}^*=\RR^n$ and view elements
$\gamma=\sum_i \gamma_i\varepsilon_i$ of $\RR^n$ as sequences
$$\gamma=(\gamma_{-n},\ldots,\gamma_{-1};\gamma_1,\ldots,\gamma_n),
\qquad\hbox{such that}\qquad \gamma_{-i}=-\gamma_i.
\formula$$
The root system of type $C_n$ is given by the sets
$$
R = \left\{\pm2\varepsilon_i,\pm(\varepsilon_j\pm\varepsilon_i)
\ |\  1\leq i,j\leq n\right\}
\quad\hbox{and}\quad
R^+=\left\{2\varepsilon_i,\varepsilon_j\pm\varepsilon_i
\ |\  1\leq i<j\leq n\right\}.
\formula$$
The simple roots are given by $\alpha_1=2\varepsilon_1$,
$\alpha_i=\varepsilon_i-\varepsilon_{i-1}$, $2\le i\le n$.
The Weyl group $W=WC_n$ is the {\it hyperoctahedral group}
of permutations of $-n,\ldots,-1,1,\ldots,n$ such that
$w(-i)=-w(i)$.  This groups acts on the $\varepsilon_i$ by the
rule $w\varepsilon_i=\varepsilon_{w(i)}$, with the convention
that $\varepsilon_{-i}=-\varepsilon_i$.

For this type C case there is a nice trick.
View the root system as 
$$\eqalign{
R &= \left\{\pm(\varepsilon_j\pm\varepsilon_i)
\ |\  i<j,\  i,j\in \{\pm1,\ldots,\pm n\} \right\}
\quad\hbox{and}\cr
R^+&=\left\{\varepsilon_j-\varepsilon_i
\ |\  i<j,\ i,j\in \{\pm1,\ldots,\pm n\} \right\}. \cr
}
\formula$$
with the convention that $\varepsilon_{-i}=-\varepsilon_i$.
In this notation $\varepsilon_i-\varepsilon_{-i}=2\varepsilon_i$,
and $\varepsilon_{-i}-\varepsilon_{-j}=\varepsilon_j-\varepsilon_i$.
This way the type C root system ``looks like'' a 
type A root system and many computations can be done in the
same way as in type A.

\subsection{Rearranging $\gamma$.}

We analyze the structure of the sets ${\cal F}^{(\gamma,J)}$
as considered in (4.3).  This corresponds to when
the $q$ in the affine Hecke algebra is not a root of unity.
The analysis in this case is analogous to the method that was used
in (5.11) to create books of placed configurations in the type A case. 

Let $\gamma\in \RR^n$.  Apply an element of the Weyl group to $\gamma$
to ``arrange'' the entries of $\gamma$ so that, 
for each $i\in \{1, \ldots, n\}$, 
$$\hbox{$\gamma_i\in [z+{1\over2}, z]$,
\quad for some $z\in \ZZ$}.
\qquad\hbox{Then $\gamma_{-i}=-\gamma_i\in 
[z', z'+{1\over 2}]$,\quad for some $z'\in \ZZ$.}
$$
As in the type A case, the sets $Z(\gamma)$ and $P(\gamma)$
can be partitioned according to the $\ZZ$ cosets of the
elements of $\gamma$ and it is sufficient to consider
each $\ZZ$-coset separately and then assemble the results
in ``books of pages''.  There are three cases to consider:
\smallskip\noindent
Case ($\beta$): the $\ZZ$-coset $\beta+\ZZ$, $\beta\in ({1\over 2},1)$.
Then 
$$\gamma=(-\beta-z_n\le\cdots\le  -\beta-z_2\le -\beta-z_1\ ;\ 
\beta+z_1\le \beta+z_2 \le\cdots\le  \beta+z_n),
\qquad z_i\in \ZZ,$$
Case ($1\over2$): The $\ZZ$-coset ${1\over 2}+\ZZ$.
Then
$$\gamma=(-\hbox{$1\over2$}-z_n\le\cdots\le  -\hbox{$1\over2$}-z_2
\le -\hbox{$1\over2$}-z_1\ ;\ 
\hbox{$1\over2$}+z_1\le \hbox{$1\over2$}+z_2 
\le \cdots\le \hbox{$1\over2$}+z_n),
\qquad z_i\in \ZZ_{\ge 0},$$
Case ($0$): The $\ZZ$-coset $\ZZ$.
Then
$$\gamma=(-z_n\le \cdots\le -z_2\le -z_1\ ;\ 
z_1\le z_2 \le \cdots\le z_n),
\qquad z_i\in \ZZ_{\ge 0}.$$
It is notationally convenient to let $z_{-i} = -z_i$.

\subsection{Boxes and standard tableaux.}

Let us assume that the entries of $\gamma$ all lie in a single
$\ZZ$-coset and decribe the resulting standard tableaux.
The general case is obtained by creating books of pages
of standard tableaux where the pages correspond to the
different $\ZZ$-cosets of entries in $\gamma$.

The placed configuration of boxes is determined as follows.

\subsection{Case $\beta$, $\beta\in ({1\over2},1)$:}  

Assume that $\gamma\in \fh_{\RR}^*$ is of the form
$$\gamma=(-\beta-z_n\le\cdots\le  -\beta-z_2\le -\beta-z_1\ ;\ 
\beta+z_1\le \beta+z_2 \le\cdots\le  \beta+z_n),
\qquad z_i\in \ZZ.$$
Place boxes on two pages of infinite
graph paper.  These pages are numbered $\beta$ and
$-\beta$ and each page has the diagonals numbered consecutively with 
the elements of $\ZZ$, from bottom left to top right.  View these
two pages, page $\beta$ and page $-\beta$, as ``linked''.
For each $1\le i\le n$ place ${\rm box}_i$ on diagonal $z_i$ of 
page $\beta$ and ${\rm box}_{-i}$ on diagonal $-z_i$ of page $-\beta$.  
The boxes on each diagonal are arranged in increasing order from
top left to bottom right.  The placement
of boxes on page $-\beta$ is a $180^{\circ}$ rotation of the
placement of the boxes on page $\beta$.

Using the notation for the root system of type $C_n$ in (8.4)
$$
\eqalign{
P(\gamma) &= \left\{
\varepsilon_j-\varepsilon_i \ |\ 
\hbox{$j>i$ and ${\rm box}_i$ and ${\rm box}_j$ are in adjacent diagonals} 
\right\}
\quad\hbox{and} \cr
Z(\gamma) &= \left\{
\varepsilon_j-\varepsilon_i \ |\  
\hbox{$j>i$ and ${\rm box}_i$ and ${\rm box}_j$ are in the same diagonal}
\right\}. \cr
}
$$
Note that 
$\varepsilon_{-i}-\varepsilon_{-j}\in Z(\gamma)$
if and only if $\varepsilon_j-\varepsilon_i\in Z_\beta(\gamma)$,
and similarly
$\varepsilon_{-i}-\varepsilon_{-j}\in P_\beta(\gamma)$
if and only if $\varepsilon_j-\varepsilon_i\in P_\beta(\gamma)$.
If $J\subseteq P(\gamma)$ arrange the boxes on adjacent diagonals 
according to the rules
\smallskip
\itemitem{(1)}
if $\varepsilon_j-\varepsilon_i\in J$ place 
${\rm box}_j$ northwest of ${\rm box}_i$, and
\smallskip
\itemitem{(2)}
if $\varepsilon_j-\varepsilon_i\in P(\gamma)\backslash J$ place 
${\rm box}_j$ southeast of ${\rm box}_i$.
\smallskip\noindent
A {\it standard tableau} is a negative rotationally
symmetric filling $p$ of the $2n$ boxes with 
$-n,\ldots,-1,1,\ldots,n$ such that 
\item{(a)} $p({\rm box}_i)<p({\rm box}_j)$
if $j>i$ and 
${\rm box}_j$ and ${\rm box}_i$ are in the same diagonal,
\item{(b)} 
$p({\rm box}_i)>p({\rm box}_j)$ if $j>i$, 
${\rm box}_i$ and ${\rm box}_j$ are in adjacent diagonals
and ${\rm box}_j$ is northwest of ${\rm box}_i$,
\item{(c)} 
$p({\rm box}_i)<p({\rm box}_j)$
if $j>i$, 
${\rm box}_i$ and ${\rm box}_j$ are in adjacent diagonals
and ${\rm box}_j$ is southeast of ${\rm box}_i$.
\smallskip\noindent
The negative rotational symmetry means that the filling of the boxes
on page $-\beta$ is the same as the filling
on page $\beta$ except rotated by $180^{\circ}$ and
with all entries in the boxes multiplied by $-1$.

\bigskip\noindent
{\sl Example.}
Suppose $\beta\in ({1\over2},1)$, and  
$$\eqalign{
\gamma&=(-\beta;\beta)
+(-2,-2,-2,-1,-1,-1,0,0,0,1,1,1;-1,-1,-1,0,0,0,1,1,1,2,2,2) \cr
&=(-\beta-2,-\beta-2,-\beta-2,-\beta-1,-\beta-1,-\beta-1,
-\beta,-\beta,-\beta, -\beta+1,-\beta+1,-\beta+1;\cr
&\qquad
\beta-1,\beta-1,\beta-1,\beta,\beta,\beta,\beta+1,\beta+1,\beta+1,
\beta+2,\beta+2,\beta+2)\cr
}$$ 
and
$$\eqalign{
J &=\left\{ 
\varepsilon_4-\varepsilon_1, \varepsilon_{-1}-\varepsilon_{-4}, 
\varepsilon_4-\varepsilon_2, \varepsilon_{-2}-\varepsilon_{-1},
\varepsilon_4-\varepsilon_3, \varepsilon_{-3}-\varepsilon_{-4},
\varepsilon_5-\varepsilon_2, \varepsilon_{-2}-\varepsilon_{-5},
\right. \cr
&\phantom{J = } \left.
\varepsilon_5-\varepsilon_3, \varepsilon_{-3}-\varepsilon_{-5},
\varepsilon_7-\varepsilon_5, \varepsilon_{-5}-\varepsilon_{-7},
\varepsilon_7-\varepsilon_6, \varepsilon_{-6}-\varepsilon_{-7},
\varepsilon_8-\varepsilon_6, \varepsilon_{-6}-\varepsilon_{-8},
\right. \cr
&\phantom{J = } \left.
\varepsilon_{10}-\varepsilon_9, \varepsilon_{-9}-\varepsilon_{-10},
\varepsilon_{10}-\varepsilon_8, \varepsilon_{-8}-\varepsilon_{-10},
\varepsilon_{10}-\varepsilon_7, \varepsilon_{-7}-\varepsilon_{-10},
\varepsilon_{11}-\varepsilon_9, \varepsilon_{-9}-\varepsilon_{-11},
\right. \cr
&\phantom{J = } \left.
\varepsilon_{11}-\varepsilon_8, \varepsilon_{-8}-\varepsilon_{-11},
\varepsilon_{11}-\varepsilon_7, \varepsilon_{-7}-\varepsilon_{-11},
\varepsilon_{12}-\varepsilon_9, \varepsilon_{-9}-\varepsilon_{-12},
\right\}. \cr
}
$$
The placed configuration of boxes corresponding to 
$(\gamma,J)$ is 
$$
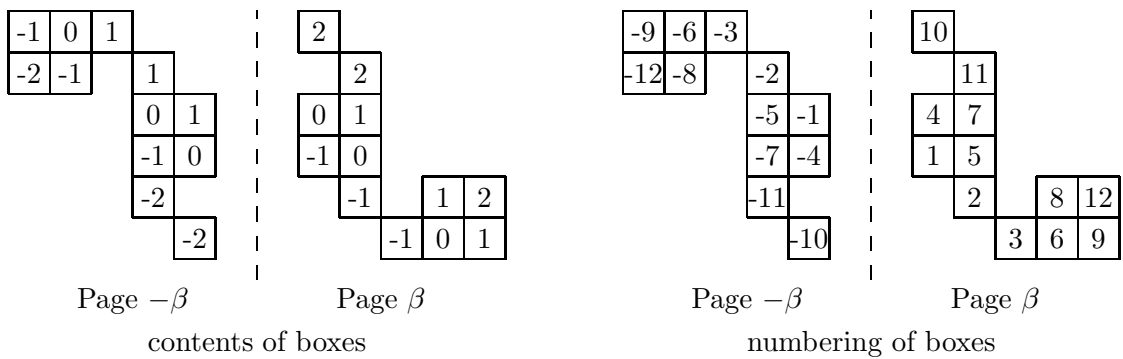
\matrix{
\beginpicture
\setcoordinatesystem units <0.5cm,0.5cm> point at 6 0       
\setplotarea x from -6 to 6, y from -1 to 6.5    
\linethickness=0.5pt                          
\put{1}  at -1.5 3.5    %
\put{1}  at -2.5 4.5    %
\put{1}  at -3.5 5.5    %
\put{0}  at -1.5 2.5    %
\put{0}  at -2.5 3.5    
\put{0}  at -4.5 5.5    %
\put{-1}  at -2.5 2.5    %
\put{-1}  at -4.5 4.5    %
\put{-1}  at -5.5 5.5    %
\put{-2}  at -1.5 0.5    %
\put{-2}  at -2.5 1.5    %
\put{-2}  at -5.5 4.5    %
\putrule from -1 0 to -2 0          %
\putrule from -1 1 to -3 1          %
\putrule from -1 2 to -3 2          %
\putrule from -1 3 to -3 3          %
\putrule from -4 4 to -6 4          %
\putrule from -1 4 to -3 4          
\putrule from -2 5 to -6 5          %
\putrule from -3 6 to -6 6          %
\putrule from -1 4 to -1 2        %
\putrule from -1 1 to -1 0        %
\putrule from -2 5 to -2 0        %
\putrule from -3 6 to -3 1        
\putrule from -4 6 to -4 4        %
\putrule from -5 6 to -5 4        %
\putrule from -6 6 to -6 4        %
\put{-1}  at 1.5 2.5    %
\put{-1}  at 2.5 1.5    %
\put{-1}  at 3.5 0.5    %
\put{0}  at 1.5 3.5    %
\put{0}  at 2.5 2.5    
\put{0}  at 4.5 0.5    %
\put{1}  at 2.5 3.5    %
\put{1}  at 4.5 1.5    %
\put{1}  at 5.5 0.5    %
\put{2}  at 1.5 5.5    %
\put{2}  at 2.5 4.5    %
\put{2}  at 5.5 1.5    %
\putrule from 1 6 to 2 6          %
\putrule from 1 5 to 3 5          %
\putrule from 1 4 to 3 4          %
\putrule from 1 3 to 3 3          %
\putrule from 4 2 to 6 2          %
\putrule from 1 2 to 3 2          
\putrule from 2 1 to 6 1          %
\putrule from 3 0 to 6 0          %
\putrule from 1 2 to 1 4        %
\putrule from 1 5 to 1 6        %
\putrule from 2 1 to 2 6        %
\putrule from 3 0 to 3 5        
\putrule from 4 0 to 4 2        %
\putrule from 5 0 to 5 2        %
\putrule from 6 0 to 6 2        %
\setdashes
\putrule from 0 -0.5 to 0 6
\put{Page $-\beta$} at -3 -1
\put{Page $\beta$} at 3 -1
\endpicture
&\qquad 
&\beginpicture
\setcoordinatesystem units <0.5cm,0.5cm> point at 6 0       
\setplotarea x from -6 to 6, y from -1 to 6.5    
\linethickness=0.5pt                          
\put{-1}  at -1.5 3.5    %
\put{-2}  at -2.5 4.5    %
\put{-3}  at -3.5 5.5    %
\put{-4}  at -1.5 2.5    %
\put{-5}  at -2.5 3.5    
\put{-6}  at -4.5 5.5    %
\put{-7}  at -2.5 2.5    %
\put{-8}  at -4.5 4.5    %
\put{-9}  at -5.5 5.5    %
\put{-10}  at -1.5 0.5    %
\put{-11}  at -2.5 1.5    %
\put{-12}  at -5.5 4.5    %
\putrule from -1 0 to -2 0          %
\putrule from -1 1 to -3 1          %
\putrule from -1 2 to -3 2          %
\putrule from -1 3 to -3 3          %
\putrule from -4 4 to -6 4          %
\putrule from -1 4 to -3 4          
\putrule from -2 5 to -6 5          %
\putrule from -3 6 to -6 6          %
\putrule from -1 4 to -1 2        %
\putrule from -1 1 to -1 0        %
\putrule from -2 5 to -2 0        %
\putrule from -3 6 to -3 1        
\putrule from -4 6 to -4 4        %
\putrule from -5 6 to -5 4        %
\putrule from -6 6 to -6 4        %
\put{1}  at 1.5 2.5    %
\put{2}  at 2.5 1.5    %
\put{3}  at 3.5 0.5    %
\put{4}  at 1.5 3.5    %
\put{5}  at 2.5 2.5    
\put{6}  at 4.5 0.5    %
\put{7}  at 2.5 3.5    %
\put{8}  at 4.5 1.5    %
\put{9}  at 5.5 0.5    %
\put{10}  at 1.5 5.5    %
\put{11}  at 2.5 4.5    %
\put{12}  at 5.5 1.5    %
\putrule from 1 6 to 2 6          %
\putrule from 1 5 to 3 5          %
\putrule from 1 4 to 3 4          %
\putrule from 1 3 to 3 3          %
\putrule from 4 2 to 6 2          %
\putrule from 1 2 to 3 2          
\putrule from 2 1 to 6 1          %
\putrule from 3 0 to 6 0          %
\putrule from 1 2 to 1 4        %
\putrule from 1 5 to 1 6        %
\putrule from 2 1 to 2 6        %
\putrule from 3 0 to 3 5        
\putrule from 4 0 to 4 2        %
\putrule from 5 0 to 5 2        %
\putrule from 6 0 to 6 2        %
\setdashes
\putrule from 0 -0.5 to 0 6
\put{Page $-\beta$} at -3 -1
\put{Page $\beta$} at 3 -1
\endpicture
\cr
\cr
\hbox{contents of boxes}
&&\hbox{numbering of boxes} \cr\cr}
$$
and a sample negative rotationally symmetric standard tableau is
$$
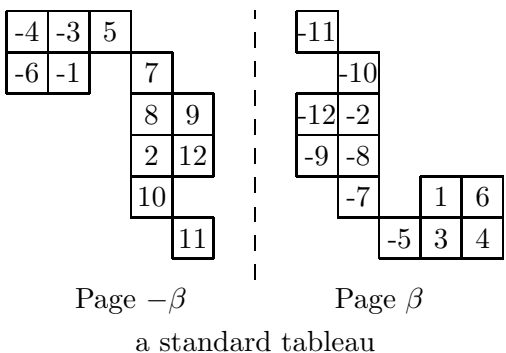
\matrix{
\cr
\beginpicture
\setcoordinatesystem units <0.5cm,0.5cm> point at 6 0       
\setplotarea x from -6 to 6, y from -1 to 6.5    
\linethickness=0.5pt                          
\put{9}  at -1.5 3.5    %
\put{7}  at -2.5 4.5    %
\put{5}  at -3.5 5.5    %
\put{12}  at -1.5 2.5    %
\put{8}  at -2.5 3.5    
\put{-3}  at -4.5 5.5    %
\put{2}  at -2.5 2.5    %
\put{-1}  at -4.5 4.5    %
\put{-4}  at -5.5 5.5    %
\put{11}  at -1.5 0.5    %
\put{10}  at -2.5 1.5    %
\put{-6}  at -5.5 4.5    %
\putrule from -1 0 to -2 0          %
\putrule from -1 1 to -3 1          %
\putrule from -1 2 to -3 2          %
\putrule from -1 3 to -3 3          %
\putrule from -4 4 to -6 4          %
\putrule from -1 4 to -3 4          
\putrule from -2 5 to -6 5          %
\putrule from -3 6 to -6 6          %
\putrule from -1 4 to -1 2        %
\putrule from -1 1 to -1 0        %
\putrule from -2 5 to -2 0        %
\putrule from -3 6 to -3 1        
\putrule from -4 6 to -4 4        %
\putrule from -5 6 to -5 4        %
\putrule from -6 6 to -6 4        %
\put{-9}  at 1.5 2.5    %
\put{-7}  at 2.5 1.5    %
\put{-5}  at 3.5 0.5    %
\put{-12}  at 1.5 3.5    %
\put{-8}  at 2.5 2.5    
\put{3}  at 4.5 0.5    %
\put{-2}  at 2.5 3.5    %
\put{1}  at 4.5 1.5    %
\put{4}  at 5.5 0.5    %
\put{-11}  at 1.5 5.5    %
\put{-10}  at 2.5 4.5    %
\put{6}  at 5.5 1.5    %
\putrule from 1 6 to 2 6          %
\putrule from 1 5 to 3 5          %
\putrule from 1 4 to 3 4          %
\putrule from 1 3 to 3 3          %
\putrule from 4 2 to 6 2          %
\putrule from 1 2 to 3 2          
\putrule from 2 1 to 6 1          %
\putrule from 3 0 to 6 0          %
\putrule from 1 2 to 1 4        %
\putrule from 1 5 to 1 6        %
\putrule from 2 1 to 2 6        %
\putrule from 3 0 to 3 5        
\putrule from 4 0 to 4 2        %
\putrule from 5 0 to 5 2        %
\putrule from 6 0 to 6 2        %
\setdashes
\putrule from 0 -0.5 to 0 6
\put{Page $-\beta$} at -3 -1
\put{Page $\beta$} at 3 -1
\endpicture
\cr
\cr
\hbox{a standard tableau} \cr
}$$
\endexample

\subsection{Case ${1\over2}$:}  

Assume that $\gamma\in \fh_{\RR}^*$ is of the form
$$\gamma=(-\hbox{$1\over2$}-z_n\le\cdots\le  -\hbox{$1\over2$}-z_2
\le -\hbox{$1\over2$}-z_1\ ;\ 
\hbox{$1\over2$}+z_1\le \hbox{$1\over2$}+z_2 
\le \cdots\le \hbox{$1\over2$}+z_n),
\qquad z_i\in \ZZ_{\ge 0},$$
Place boxes on a page of infinite graph paper
which has its diagonals numbered consecutively with the elements 
of ${1\over2}+\ZZ$, from bottom left to top right.  
This page has page number $1\over 2$.
For each $i\in \{ \pm1,\ldots,\pm n\}$
place ${\rm box}_i$ on diagonal ${1\over2}+z_i$ and
${\rm box}_{-i}$ on diagonal $-{1\over2}-z_i$.
The boxes on each diagonal are arranged in increasing order from
top left to bottom right and the 
placement of boxes is negative rotationally symmetric
in the sense that a $180^{\circ}$ rotation takes ${\rm box}_i$
to ${\rm box}_{-i}$.

Using the root system notation in (8.4),
$$\eqalign{
P(\gamma) &= \{
\varepsilon_j-\varepsilon_i \ |\ 
\hbox{ $j>i$ and 
${\rm box}_i$ and ${\rm box}_j$ are in adjacent diagonals}\}
\quad\hbox{and}\cr
Z(\gamma) &= \{
\varepsilon_j-\varepsilon_i \ |\ \hbox{
$j>i$ and ${\rm box}_i$ and ${\rm box}_j$ are in the same diagonal} 
\}. 
\cr} $$
Note that it is the formulation of the root system of type $C_n$ in (8.4)
which makes the description of $P(\gamma)$ and $Z(\gamma)$ nice
in this case.
If $J\subseteq P(\gamma)$ arrange the boxes on adjacent diagonals 
according to the rules
\smallskip
\itemitem{(1)}
if $\varepsilon_j-\varepsilon_i\in J$ place 
${\rm box}_j$ northwest of ${\rm box}_i$, and
\smallskip
\itemitem{(2)}
if $\varepsilon_j-\varepsilon_i\in P(\gamma)\backslash J$ place 
${\rm box}_j$ southeast of ${\rm box}_i$.
\smallskip\noindent
A {\it standard tableau} is a negative rotationally
symmetric filling $p$ of the $2n$ boxes with 
$-n,\ldots,-1,1,\ldots,n$ such that 
\item{(a)} $p({\rm box}_i)<p({\rm box}_j)$
if $j>i$ and 
${\rm box}_j$ and ${\rm box}_i$ are in the same diagonal,
\item{(b)} 
$p({\rm box}_i)>p({\rm box}_j)$ if $j>i$, 
${\rm box}_i$ and ${\rm box}_j$ are in adjacent diagonals
and ${\rm box}_j$ is northwest of ${\rm box}_i$,
\item{(c)} 
$p({\rm box}_i)<p({\rm box}_j)$
if $j>i$, 
${\rm box}_i$ and ${\rm box}_j$ are in adjacent diagonals
and ${\rm box}_j$ is southeast of ${\rm box}_i$.
\smallskip\noindent
The negative rotational symmetry means that the filling of the boxes
is the same if all entries in the boxes multiplied by $-1$ and 
it is rotated by $180^{\circ}$.

\bigskip\noindent
{\sl Example.}  
Suppose $\gamma=
\left(-{7\over 2},-{5\over2},-{5\over2}, 
-{3\over2},-{3\over2},-{3\over2},-{3\over2},
-{1\over2},-{1\over2},-{1\over2},-{1\over2};
{1\over2},{1\over2},{1\over2},{1\over2},
{3\over2},{3\over2},{3\over2},{3\over2},
{5\over2},{5\over2},{7\over2}\right)$ 
and
$$\eqalign{
J 
&=\left\{ 
\varepsilon_{11}-\varepsilon_{10}, \varepsilon_{-10}-\varepsilon_{-11}, 
\varepsilon_{10}-\varepsilon_8, \varepsilon_{-8}-\varepsilon_{-10},
\varepsilon_9-\varepsilon_7, \varepsilon_{-7}-\varepsilon_{-9}, 
\varepsilon_9-\varepsilon_8, \varepsilon_{-8}-\varepsilon_{-9},
\right.\cr
&\phantom{J = } \left.
\varepsilon_7-\varepsilon_3, \varepsilon_{-3}-\varepsilon_{-7}, 
\varepsilon_7-\varepsilon_4, \varepsilon_{-4}-\varepsilon_{-7},
\varepsilon_6-\varepsilon_2, \varepsilon_{-2}-\varepsilon_{-6}, 
\varepsilon_6-\varepsilon_3, \varepsilon_{-3}-\varepsilon_{-6},
\right.\cr
&\phantom{J = } \left.
\varepsilon_6-\varepsilon_4, \varepsilon_{-4}-\varepsilon_{-6}, 
\varepsilon_5-\varepsilon_4, \varepsilon_{-4}-\varepsilon_{-5},
\varepsilon_5-\varepsilon_3, \varepsilon_{-3}-\varepsilon_{-5},
\varepsilon_5-\varepsilon_2, \varepsilon_{-2}-\varepsilon_{-5},
\right. \cr
&\phantom{J = } \left.
\varepsilon_2-\varepsilon_{-1}, 
\varepsilon_3-\varepsilon_{-1},
\varepsilon_4-\varepsilon_{-1}, 
\varepsilon_1-\varepsilon_{-1} 
\right\}. \cr
&=\left\{ 
\varepsilon_{11}-\varepsilon_{10}, 
\varepsilon_{10}-\varepsilon_8, 
\varepsilon_9-\varepsilon_7,  
\varepsilon_9-\varepsilon_8, 
\varepsilon_7-\varepsilon_3, 
\varepsilon_7-\varepsilon_4,
\varepsilon_6-\varepsilon_2, 
\varepsilon_6-\varepsilon_3,
\varepsilon_6-\varepsilon_4, 
\varepsilon_5-\varepsilon_4,
\right. \cr
&\phantom{J = } \left.
\varepsilon_5-\varepsilon_3,
\varepsilon_5-\varepsilon_2,
\varepsilon_2+\varepsilon_1, 
\varepsilon_3+\varepsilon_1,
\varepsilon_4+\varepsilon_1, 
2\varepsilon_1 
\right\}. \cr
}
$$
The placed configuration of boxes corresponding to 
$(\gamma,J)$
is as given below.
$$
\matrix{
\beginpicture
\setcoordinatesystem units <0.5cm,0.5cm> point at 6 0       
\setplotarea x from 0 to 8, y from -1 to 7    
\linethickness=0.5pt                          
\put{-${7\over2}$}  at 0.5 3.5    %
\put{-${5\over2}$}  at 0.5 4.5    %
\put{-${5\over2}$}  at 2.5 2.5    %
\put{-${3\over2}$}  at 0.5 5.5    %
\put{-${3\over2}$}  at 2.5 3.5    %
\put{-${3\over2}$}  at 3.5 2.5    %
\put{-${3\over2}$}  at 4.5 1.5    %
\put{-${1\over2}$}  at 1.5 5.5    %
\put{-${1\over2}$}  at 2.5 4.5    %
\put{-${1\over2}$}  at 3.5 3.5    %
\put{-${1\over2}$}  at 6.5 0.5    %
\put{$1\over2$}  at 1.5 6.5    %
\put{$1\over2$}  at 4.5 3.5    %
\put{$1\over2$}  at 5.5 2.5    %
\put{$1\over2$}  at 6.5 1.5    %
\put{$3\over2$}  at 3.5 5.5    %
\put{$3\over2$}  at 4.5 4.5    %
\put{$3\over2$}  at 5.5 3.5    %
\put{$3\over2$}  at 7.5 1.5    %
\put{$5\over2$}  at 5.5 4.5    %
\put{$5\over2$}  at 7.5 2.5    %
\put{$7\over2$}  at 7.5 3.5    %
\putrule from 1 7 to 2 7          %
\putrule from 0 6 to 2 6          %
\putrule from 3 6 to 4 6          %
\putrule from 0 5 to 6 5          %
\putrule from 0 4 to 1 4          %
\putrule from 2 4 to 6 4          %
\putrule from 7 4 to 8 4          %
\putrule from 0 3 to 1 3          %
\putrule from 2 3 to 6 3          %
\putrule from 7 3 to 8 3          %
\putrule from 2 2 to 8 2          
\putrule from 4 1 to 5 1          %
\putrule from 6 1 to 8 1          %
\putrule from 6 0 to 7 0          %
\putrule from 0 3 to 0 6        %
\putrule from 1 3 to 1 7        %
\putrule from 2 2 to 2 7        %
\putrule from 3 2 to 3 6        
\putrule from 4 1 to 4 6        %
\putrule from 5 1 to 5 5        %
\putrule from 6 0 to 6 5        %
\putrule from 7 0 to 7 4        %
\putrule from 8 1 to 8 4        %
\endpicture
&\qquad
&\beginpicture
\setcoordinatesystem units <0.5cm,0.5cm> point at 6 0       
\setplotarea x from 0 to 8, y from -1 to 7    
\linethickness=0.5pt                          
\put{-11}  at 0.5 3.5    %
\put{-10}  at 0.5 4.5    %
\put{-9}  at 2.5 2.5    %
\put{-8}  at 0.5 5.5    %
\put{-7}  at 2.5 3.5    %
\put{-6}  at 3.5 2.5    %
\put{-5}  at 4.5 1.5    %
\put{-4}  at 1.5 5.5    %
\put{-3}  at 2.5 4.5    %
\put{-2}  at 3.5 3.5    %
\put{-1}  at 6.5 0.5    %
\put{1}  at 1.5 6.5    %
\put{2}  at 4.5 3.5    %
\put{3}  at 5.5 2.5    %
\put{4}  at 6.5 1.5    %
\put{5}  at 3.5 5.5    %
\put{6}  at 4.5 4.5    %
\put{7}  at 5.5 3.5    %
\put{8}  at 7.5 1.5    %
\put{9}  at 5.5 4.5    %
\put{10}  at 7.5 2.5    %
\put{11}  at 7.5 3.5    %
\putrule from 1 7 to 2 7          %
\putrule from 0 6 to 2 6          %
\putrule from 3 6 to 4 6          %
\putrule from 0 5 to 6 5          %
\putrule from 0 4 to 1 4          %
\putrule from 2 4 to 6 4          %
\putrule from 7 4 to 8 4          %
\putrule from 0 3 to 1 3          %
\putrule from 2 3 to 6 3          %
\putrule from 7 3 to 8 3          %
\putrule from 2 2 to 8 2          
\putrule from 4 1 to 5 1          %
\putrule from 6 1 to 8 1          %
\putrule from 6 0 to 7 0          %
\putrule from 0 3 to 0 6        %
\putrule from 1 3 to 1 7        %
\putrule from 2 2 to 2 7        %
\putrule from 3 2 to 3 6        
\putrule from 4 1 to 4 6        %
\putrule from 5 1 to 5 5        %
\putrule from 6 0 to 6 5        %
\putrule from 7 0 to 7 4        %
\putrule from 8 1 to 8 4        %
\endpicture
&\qquad
&\beginpicture
\setcoordinatesystem units <0.5cm,0.5cm> point at 6 0       
\setplotarea x from 0 to 8, y from -1 to 7    
\linethickness=0.5pt                          
\put{-3}  at 0.5 3.5    %
\put{-5}  at 0.5 4.5    %
\put{8}  at 2.5 2.5    %
\put{-7}  at 0.5 5.5    %
\put{-2}  at 2.5 3.5    %
\put{9}  at 3.5 2.5    %
\put{10}  at 4.5 1.5    %
\put{-6}  at 1.5 5.5    %
\put{-4}  at 2.5 4.5    %
\put{-1}  at 3.5 3.5    %
\put{11}  at 6.5 0.5    %
\put{-11}  at 1.5 6.5    %
\put{1}  at 4.5 3.5    %
\put{4}  at 5.5 2.5    %
\put{6}  at 6.5 1.5    %
\put{-10}  at 3.5 5.5    %
\put{-9}  at 4.5 4.5    %
\put{2}  at 5.5 3.5    %
\put{7}  at 7.5 1.5    %
\put{-8}  at 5.5 4.5    %
\put{5}  at 7.5 2.5    %
\put{3}  at 7.5 3.5    %
\putrule from 1 7 to 2 7          %
\putrule from 0 6 to 2 6          %
\putrule from 3 6 to 4 6          %
\putrule from 0 5 to 6 5          %
\putrule from 0 4 to 1 4          %
\putrule from 2 4 to 6 4          %
\putrule from 7 4 to 8 4          %
\putrule from 0 3 to 1 3          %
\putrule from 2 3 to 6 3          %
\putrule from 7 3 to 8 3          %
\putrule from 2 2 to 8 2          
\putrule from 4 1 to 5 1          %
\putrule from 6 1 to 8 1          %
\putrule from 6 0 to 7 0          %
\putrule from 0 3 to 0 6        %
\putrule from 1 3 to 1 7        %
\putrule from 2 2 to 2 7        %
\putrule from 3 2 to 3 6        
\putrule from 4 1 to 4 6        %
\putrule from 5 1 to 5 5        %
\putrule from 6 0 to 6 5        %
\putrule from 7 0 to 7 4        %
\putrule from 8 1 to 8 4        %
\endpicture
\cr
\hbox{Page ${1\over2}$}
&&\hbox{Page ${1\over2}$}
&&\hbox{Page ${1\over2}$}\cr
\cr
\hbox{contents of boxes}
&&\hbox{numbering of boxes}
&&\hbox{a standard tableau} \cr
}
$$
\endexample

\subsection{Case $0$:}  

Assume that $\gamma\in \fh_{\RR}^*$ is of the form
$$\gamma=(-z_n\le \cdots\le -z_2\le -z_1\ ;\ 
z_1\le z_2 \le \cdots\le z_n),
\qquad z_i\in \ZZ_{\ge 0}.$$
Place boxes on a page of infinite graph paper
which has its diagonals numbered consecutively with the elements 
of $\ZZ$, from bottom left to top right.  
This page has page number $0$.
For each $i\in \{ \pm1,\ldots,\pm n\}$
place ${\rm box}_i$ on diagonal $z_i$ and
${\rm box}_{-i}$ on diagonal $-z_i$.
The boxes on each diagonal are arranged in increasing order from
top left to bottom right and the 
placement of boxes is negative rotationally symmetric
in the sense that a $180^{\circ}$ rotation takes ${\rm box}_i$
to ${\rm box}_{-i}$.

Using the root system notation in (8.4),
$$\eqalign{
P(\gamma) &= \{
\varepsilon_j-\varepsilon_i \ |\ 
\hbox{ $j>i$ and 
${\rm box}_i$ and ${\rm box}_j$ are in adjacent diagonals}\}
\quad\hbox{and}\cr
Z(\gamma) &= \{
\varepsilon_j-\varepsilon_i \ |\ \hbox{
$j>i$ and ${\rm box}_i$ and ${\rm box}_j$ are in the same diagonal} 
\}. 
\cr} $$
If $J\subseteq P(\gamma)$ arrange the boxes on adjacent diagonals 
according to the rules
\smallskip
\itemitem{(1)}
if $\varepsilon_j-\varepsilon_i\in J$ place 
${\rm box}_j$ northwest of ${\rm box}_i$, and
\smallskip
\itemitem{(2)}
if $\varepsilon_j-\varepsilon_i\in P(\gamma)\backslash J$ place 
${\rm box}_j$ southeast of ${\rm box}_i$.
\smallskip\noindent
A {\it standard tableau} is a negative rotationally
symmetric filling $p$ of the $2n$ boxes with 
$-n,\ldots,-1,1,\ldots,n$ such that 
\item{(a)} $p({\rm box}_i)<p({\rm box}_j)$
if $j>i$ and 
${\rm box}_j$ and ${\rm box}_i$ are in the same diagonal,
\item{(b)} 
$p({\rm box}_i)>p({\rm box}_j)$ if $j>i$, 
${\rm box}_i$ and ${\rm box}_j$ are in adjacent diagonals
and ${\rm box}_j$ is northwest of ${\rm box}_i$,
\item{(c)} 
$p({\rm box}_i)<p({\rm box}_j)$
if $j>i$, 
${\rm box}_i$ and ${\rm box}_j$ are in adjacent diagonals
and ${\rm box}_j$ is southeast of ${\rm box}_i$.
\smallskip\noindent
The negative rotational symmetry means that the filling of the boxes
is the same if all entries in the boxes multiplied by $-1$ and 
it is rotated by $180^{\circ}$.

\bigskip\noindent
{\sl Example.}  
Suppose $\gamma=(-2,-1,-1,-1,0,0,0;0,0,0,1,1,1,2)$ and
$$\eqalign{
J 
&=\left\{ 
\varepsilon_4-\varepsilon_1, \varepsilon_{-1}-\varepsilon_{-4}, 
\varepsilon_4-\varepsilon_2, \varepsilon_{-2}-\varepsilon_{-4},
\varepsilon_4-\varepsilon_3, \varepsilon_{-3}-\varepsilon_{-2}, 
\varepsilon_5-\varepsilon_1, \varepsilon_{-1}-\varepsilon_{-5},
\right. \cr
&\phantom{J = } \left.
\varepsilon_5-\varepsilon_2, \varepsilon_{-2}-\varepsilon_{-5}, 
\varepsilon_5-\varepsilon_3, \varepsilon_{-3}-\varepsilon_{-5},
\varepsilon_6-\varepsilon_1, \varepsilon_{-1}-\varepsilon_{-6}, 
\varepsilon_6-\varepsilon_2, \varepsilon_{-2}-\varepsilon_{-6},
\right. \cr
&\phantom{J = } \left.
\varepsilon_6-\varepsilon_3, \varepsilon_{-3}-\varepsilon_{-6}, 
\varepsilon_7-\varepsilon_6, \varepsilon_{-6}-\varepsilon_{-7},
\varepsilon_6-\varepsilon_{-1}, \varepsilon_1-\varepsilon_{-6},
\varepsilon_5-\varepsilon_{-1}, \varepsilon_1-\varepsilon_{-5},
\right. \cr
&\phantom{J = } \left.
\varepsilon_4-\varepsilon_{-1}, \varepsilon_1-\varepsilon_{-4},
\varepsilon_5-\varepsilon_{-2}, \varepsilon_2-\varepsilon_{-5}, 
\varepsilon_4-\varepsilon_{-2}, \varepsilon_2-\varepsilon_{-4} 
\right\}. \cr
&=\left\{ \varepsilon_4-\varepsilon_1, \varepsilon_4-\varepsilon_2,
\varepsilon_4-\varepsilon_3, \varepsilon_5-\varepsilon_1,
\varepsilon_5-\varepsilon_2, \varepsilon_5-\varepsilon_3,
\varepsilon_6-\varepsilon_1, \varepsilon_6-\varepsilon_2,
\varepsilon_6-\varepsilon_3, \varepsilon_7-\varepsilon_6,\right. \cr
&\phantom{J = } \left.
\varepsilon_6+\varepsilon_1,
\varepsilon_5+\varepsilon_1,
\varepsilon_5+\varepsilon_2, \varepsilon_4+\varepsilon_1,
\varepsilon_4+\varepsilon_2 \right\}. \cr
}
$$
The placed configuration of boxes corresponding to $(\gamma,J)$
is as given below.
$$
\matrix{
\beginpicture
\setcoordinatesystem units <0.5cm,0.5cm> point at 6 0       
\setplotarea x from 0 to 8, y from -1 to 7    
\linethickness=0.5pt                          
\put{0}  at 0.5 7.5    %
\put{0}  at 2.5 5.5    %
\put{0}  at 3.5 4.5    %
\put{-1}  at 4.5 2.5    %
\put{-1}  at 5.5 1.5    %
\put{-1}  at 6.5 0.5    %
\put{-2}  at 4.5 1.5    %
\put{0}  at 4.5 3.5    
\put{0}  at 5.5 2.5    %
\put{0}  at 7.5 0.5    %
\put{1}  at 1.5 7.5    %
\put{1}  at 2.5 6.5    %
\put{1}  at 3.5 5.5    %
\put{2}  at 3.5 6.5    %
\putrule from 0 8 to 2 8          %
\putrule from 0 7 to 4 7          %
\putrule from 2 6 to 4 6          %
\putrule from 2 5 to 4 5          %
\putrule from 3 4 to 5 4          %
\putrule from 4 3 to 6 3          %
\putrule from 4 2 to 6 2          
\putrule from 4 1 to 8 1          %
\putrule from 6 0 to 8 0          %
\putrule from 0 7 to 0 8        %
\putrule from 1 7 to 1 8        %
\putrule from 2 5 to 2 8        %
\putrule from 3 4 to 3 7        %
\putrule from 4 1 to 4 7        
\putrule from 5 1 to 5 4        %
\putrule from 6 0 to 6 3        %
\putrule from 7 0 to 7 1        %
\putrule from 8 0 to 8 1        %
\endpicture
&\qquad
&\beginpicture
\setcoordinatesystem units <0.5cm,0.5cm> point at 6 0       
\setplotarea x from 0 to 8, y from -1 to 7    
\linethickness=0.5pt                          
\put{-3}  at 0.5 7.5    %
\put{-2}  at 2.5 5.5    %
\put{-1}  at 3.5 4.5    %
\put{-6}  at 4.5 2.5    %
\put{-5}  at 5.5 1.5    %
\put{-4}  at 6.5 0.5    %
\put{-7}  at 4.5 1.5    %
\put{1}  at 4.5 3.5    
\put{2}  at 5.5 2.5    %
\put{3}  at 7.5 0.5    %
\put{4}  at 1.5 7.5    %
\put{5}  at 2.5 6.5    %
\put{6}  at 3.5 5.5    %
\put{7}  at 3.5 6.5    %
\putrule from 0 8 to 2 8          %
\putrule from 0 7 to 4 7          %
\putrule from 2 6 to 4 6          %
\putrule from 2 5 to 4 5          %
\putrule from 3 4 to 5 4          %
\putrule from 4 3 to 6 3          %
\putrule from 4 2 to 6 2          
\putrule from 4 1 to 8 1          %
\putrule from 6 0 to 8 0          %
\putrule from 0 7 to 0 8        %
\putrule from 1 7 to 1 8        %
\putrule from 2 5 to 2 8        %
\putrule from 3 4 to 3 7        %
\putrule from 4 1 to 4 7        
\putrule from 5 1 to 5 4        %
\putrule from 6 0 to 6 3        %
\putrule from 7 0 to 7 1        %
\putrule from 8 0 to 8 1        %
\endpicture
&\qquad
&\beginpicture
\setcoordinatesystem units <0.5cm,0.5cm> point at 6 0       
\setplotarea x from 0 to 8, y from -1 to 7    
\linethickness=0.5pt                          
\put{-7}  at 0.5 7.5    %
\put{-3}  at 2.5 5.5    %
\put{-1}  at 3.5 4.5    %
\put{2}  at 4.5 2.5    %
\put{5}  at 5.5 1.5    %
\put{6}  at 6.5 0.5    %
\put{4}  at 4.5 1.5    %
\put{1}  at 4.5 3.5    
\put{3}  at 5.5 2.5    %
\put{7}  at 7.5 0.5    %
\put{-6}  at 1.5 7.5    %
\put{-5}  at 2.5 6.5    %
\put{-2}  at 3.5 5.5    %
\put{-4}  at 3.5 6.5    %
\putrule from 0 8 to 2 8          %
\putrule from 0 7 to 4 7          %
\putrule from 2 6 to 4 6          %
\putrule from 2 5 to 4 5          %
\putrule from 3 4 to 5 4          %
\putrule from 4 3 to 6 3          %
\putrule from 4 2 to 6 2          
\putrule from 4 1 to 8 1          %
\putrule from 6 0 to 8 0          %
\putrule from 0 7 to 0 8        %
\putrule from 1 7 to 1 8        %
\putrule from 2 5 to 2 8        %
\putrule from 3 4 to 3 7        %
\putrule from 4 1 to 4 7        
\putrule from 5 1 to 5 4        %
\putrule from 6 0 to 6 3        %
\putrule from 7 0 to 7 1        %
\putrule from 8 0 to 8 1        %
\endpicture
\cr
\hbox{Page $0$}
&&\hbox{Page $0$}
&&\hbox{Page $0$}\cr
\cr
\hbox{contents of boxes}
&&\hbox{numbering of boxes}
&&\hbox{a standard tableau} \cr
}
$$
\endexample

\subsection{}

A posteriori the analysis of the three cases, Case $\beta$, 
Case ${1\over2}$, and Case $0$, it becomes evident that the
trick of using the formulation of the root system of type $C_n$
in (8.4) provides a completely uniform
description of the configurations of boxes
and standard tableaux corresponding to type $C_n$ local
regions.  All three cases give
negative rotationally invariant tableaux.
We could not ask for nature to work out more perfectly.

\bigskip
\bigskip
\centerline{\smallcaps References}
\bigskip

\medskip
\item{[AK]} {\smallcaps S.\ Ariki and K.\ Koike}, {\it A Hecke algebra of
$(\ZZ/r\ZZ)\wr S_n$ and construction of its irreducible representations},
Adv.\ in Math.\ {\bf 106} (1994), 216--243.

\medskip
\item{[AL]} {\smallcaps C. Athanasiadis and S. Linusson}, 
{\it A simple bijection for the regions of the Shi arrangement of hyperplanes},
Discrete Math., to appear.

\medskip
\item{[Au]} {\smallcaps M. Aubert}, {\it Dualit\'e dans le groupe de
Grothendieck de la categorie des repr\'esentations lisses de longueur
finie d'un groupe r\'eductif $p$-adique},  Trans. Amer. Soc. {\bf 347}
(1995), 2179--2189.

\medskip
\item{[Bj]} {\smallcaps A.\ Bj\"orner}, {\it Orderings of Coxeter groups},
{\sl Combinatorics and Algebra (Boulder, Colo. 1983)}, Contemp. Math.
{\bf 34}, Amer. Math. Soc., Providence 1984, 175--195.

\medskip
\item{[BW]} {\smallcaps A.\ Bj\"orner and M.\ Wachs},
{\it Generalized quotients in Coxeter groups}, 
Trans.\ Amer.\ Math.\ Soc.\ {\bf 308} no.\ 1 (1988), 1--37.

\medskip
\item{[Bou]} {\smallcaps N.\ Bourbaki}, 
{\it Groupes et alg\`ebres de Lie, Chapitres 4,5 et 6},
Elements de Math\'ematique, Hermann, Paris 1968.

\medskip
\item{[CSM]} {\smallcaps R.\ Carter, G.\ Segal and I.G.\ Macdonald}
{\sl Lectures on Lie groups and Lie algebras}, 
London Mathematical Society Student Texts {\bf 32},
Cambridge University Press, Cambridge, 1995. 

\medskip
\item{[CG]} {\smallcaps N.\ Chriss and V.\ Ginzburg}, 
{\sl Representation Theory and Complex Geometry}, Birkh\"auser, 1997.

\medskip
\item{[GR]} {\smallcaps A.\ Garsia and C.\ Reutenauer},
{\it A decomposition of Solomon's descent algebra}, Adv.\ Math.\ {\bf 77} (1989), 
no.\ 2, 189--262. 

\medskip
\item{[H]} {\smallcaps P.N.\ Hoefsmit}, 
{\it Representations of Hecke algebras of finite groups with $BN$-pairs 
of classical type}, 
Ph.D.\ Thesis, University of British Columbia, 1974.

\medskip
\item{[IM]} {\smallcaps N.\ Iwahori and H.\ Matsumoto}, 
{\it On some Bruhat decomposition and the structure of the Hecke rings of 
$p$-adic Chevalley groups}, 
Publ.\ Math.\ IHES {\bf 40} (1972), 81--116.

\medskip\medskip
\item{[Kt1]}  {\smallcaps S.\ Kato}, 
{\it Irreducibility of principal series representations for Hecke algebras
of affine type}, J.\ Fac.\ Sci.\ Univ.\ Tokyo sec.\ 1A {\bf 28} (1981),
929--943.

\medskip
\item{[Kt2]}  {\smallcaps S.\ Kato}, {\it Duality for representations of
a Hecke algebra}, Proc.\ Amer.\ Math.\ Soc.\ {\bf 119} (1993), 941--946.

\medskip
\item{[KL]} {\smallcaps D.\ Kazhdan and G.\ Lusztig}, 
{\it Proof of the Deligne-Langlands conjecture for Hecke algebras}, 
Invent.\ Math.\ {\bf 87} (1987), 153--215.

\medskip
\item{[Ks]} {\smallcaps B.\ Kostant},
{\it On Dale Peterson's $2^{\rm rank}$ abelian ideal theorem, symmetric spaces,
discrete series, and Euler number multiplets of representations},
to appear.

\medskip
\item{[KOP]} {\smallcaps C.\ Krattenthaler, L.\ Orsina and P.\ Papi},
{\it Enumeration of ad-nilpotent ${\frak b}$-ideals for simple Lie algebras}
to appear in Adv.\ in Appl.\ Math.

\medskip
\item{[KR]} {\smallcaps C.\ Kriloff and A.\ Ram},
{\it Representations of graded Hecke algebras},
Representation Theory {\bf 6} (2002), 31--69, 
{\tt http://www.ams.org/ert/home-2002.html}.

\medskip
\item{[KZ]} {\smallcaps H.\ Knight and A.\ Zelevinsky},
{\it Representations of quivers of type A and the multisegment
duality}, Adv.\ Math.\ {\bf 117} (1996), no. 2, 273--293.

\medskip
\item{[LTV]} {\smallcaps B.\ Leclerc, J.-Y.\ Thibon and E.\ Vasserot}
{\it Zelevinsky's involution at roots of unity},
J.\ reine angew.\ Math.\ {\bf 513} (1999), 33--51.

\medskip
\item{[Li1]} {\smallcaps P.\ Littelmann}, {\it Paths and root operators in
representation theory}, Ann.\ of Math.\ (2) {\bf 142} (1995), 499--525.

\medskip
\item{[Li2]} {\smallcaps P.\ Littelmann}, {\it The path model for
representations of symmetrizable Kac-Moody algebras}, Proceedings of the
International Congress of Mathematicians, Vol.\ 1 (Z\"urich, 1994), 298--308,
Birkh\"auser, Basel, 1995. 

\medskip
\item{[Lo]} {\smallcaps J.\ Losonczy},
{\it Standard Young tableaux in the Weyl group setting},
J.\ Algebra {\bf 220} (1999), 255--260.

\medskip
\item{[Lu]} {\smallcaps G.\ Lusztig},
{\it Singularities, character formulas, and a $q$-analog of weight
multiplicities}, Analysis and topology on singular spaces, II, III 
(Luminy, 1981),
 Ast\'erisque {\bf 101-102}, Soc. Math. France, Paris, 1983, 208--229.

\medskip
\item{[Mac]} {\smallcaps I.G. Macdonald}, 
Symmetric functions and Hall polynomials, Second edition, 
Oxford Mathematical Monographs, Oxford Univ. Press, New York, 1995.


\medskip
\item{[MW]} {\smallcaps C.\ Moeglin and J.-L.\ Waldspurger},
{\it L'involution de Zelevinski}, J.\ Reine Angew.\ Math.\ {\bf 372} (1986),
136--177. 

\medskip
\item{[NT]} {\smallcaps M.\ Nazarov and V.\ Tarasov},
{\it Representations of Yangians with Gelfand-Zetlin bases},
J.\ Reine Angew.\ Math.\ {\bf 496} (1998),
181--212. 

\medskip
\item{[R1]} {\smallcaps A.\ Ram}, 
{\it Irreducible representations of rank two affine Hecke algebras}, 
preprint 1998.

\medskip
\item{[R2]} {\smallcaps A.\ Ram}, 
{\it Calibrated representations of affine Hecke algebras},
preprint 1998.

\medskip
\item{[R3]} {\smallcaps A.\ Ram}, 
{\it Standard Young tableaux for finite root systems},
preprint 1998.

\medskip
\item{[RR]} {\smallcaps A.\ Ram and J. Ramagge}, 
{\it  Affine Hecke algebras, cyclotomic Hecke algebras
and Clifford theory}, preprint 1999.

\medskip\noindent
\item{[Re]} {\smallcaps C.\ Reutenauer},
{\sl Free Lie algebras},
London Mathematical Society Monographs, New Series {\bf 7}, 
Oxford Science Publications, Oxford University Press, New York, 1993. 

\medskip
\item{[Rg]} {\smallcaps J. Rogawski}, {\it On modules over the Hecke algebra
of a $p$-adic group}, Invent. Math. {\bf 79} (1985), 443--465.

\medskip
\item{[Sh1]} {\smallcaps J.-Y. Shi}, {\it The number of $\oplus$-sign types},
Quart. J. Math. Oxford Ser. (2) {\bf 48} (1997), 93--105.

\medskip
\item{[Sh2]} {\smallcaps J.-Y. Shi}, {\it  Left cells in affine Weyl
groups}, T\^ohoku Math. J. (2) {\bf 46} (1994), no. 1, 105--124. 

\medskip
\item{[Sh3]} {\smallcaps J.-Y. Shi}, 
{\it  Sign types corresponding to an affine Weyl group}, 
J. London Math. Soc. (2) {\bf 35} (1987), 56--74.

\medskip
\item{[So]} {\smallcaps L.\ Solomon},
{\it A Mackey formula in the group ring of a Coxeter group},
J.\ Algebra {\bf 41} (1976), no.\ 2, 255--264. 

\medskip
\item{[ST]} {\smallcaps L. Solomon and H. Terao}, 
{\it The double Coxeter arrangement}, preprint 1997.
 
\medskip
\item{[St1]} {\smallcaps R. Stanley}, 
{\it Hyperplane arrangements, interval orders, and trees}, 
Proc. Nat. Acad. Sci. U.S.A. {\bf 93} (1996), 2620--2625.

\medskip
\item{[St2]} {\smallcaps R. Stanley}, 
{\it Hyperplane arrangements, parking functions and tree inversions}, 
in {\it Mathematical Essays in Honor of Gian-Carlo Rota}, 
Birkh\"auser, Boston, 1998.

\medskip
\item{[Sb1]} {\smallcaps R.\ Steinberg}, 
{\sl Lectures on Chevalley groups}, Notes prepared by John Faulkner and Robert Wilson,
Yale University, New Haven, Conn., 1968.

\medskip
\item{[Sb2]} {\smallcaps R.\ Steinberg}, 
{\it Endomorphisms of linear algebraic groups}, 
Memoirs of the American Mathematical Society, {\bf 80} (1968) 1--108.

\medskip
\item{[Sb3]} {\smallcaps R.\ Steinberg},
{\it On a theorem of Pittie}, Topology {\bf 14} (1975), 173--177.

\medskip
\item{[Su]} {\smallcaps S.\ Sundaram},
{\it Tableaux in the representation theory of the classical Lie groups}
in {\sl Invariant theory and tableaux} (Minneapolis, MN, 1988),
191--225, IMA Vol.\ Math.\ Appl., {\bf 19}, Springer, New York, 1990. 

\medskip
\item{[Wa]} {\smallcaps N.\ Wallach},
{\sl Real reductive groups I}, 
Pure and Applied Mathematics {\bf 132},
Academic Press, Inc., Boston, MA, 1988. 

\medskip
\item{[Wz]} {\smallcaps H.\ Wenzl},
{\it Hecke algebras of type $A_n$ and subfactors},
Invent.\ Math.\ {\bf 92} (1988), 349--383.

\medskip
\item{[Xi]} {\smallcaps N. Xi}, 
{\it  Representations of affine Hecke algebras}, 
Lect.\ Notes in Math.\ {\bf 1587}, Springer-Verlag, Berlin, 1994.

\medskip
\item{[Yg1]} {\smallcaps A.\ Young}, 
{\it The collected papers of Alfred Young 1873-1940},
Mathematical Expositions {\bf 21},
University of Toronto Press, 1977.

\medskip
\item{[Yg2]} {\smallcaps A.\ Young}, 
{\it On quantitative substitutional analysis
(sixth paper)}, Proc. London. Math. Soc. (2) {\bf 34} (1931), 196--230.

\medskip
\item{[Ze]} {\smallcaps A. Zelevinsky}, 
{\it Induced Representations of $\gp$-adic groups II: 
On irreducible representations of $GL(n)$}, 
Ann.\ Scient.\  \'Ec.\ Norm.\  Sup. $4^{\rm e}$ s\`erie, 
{\bf 13} (1980), 165--210.

\vfill\eject
\end

\subsection {Irreducibility of principal series modules.}

Our next goal is to prove Theorem ??? which determines exactly
when the principal series module $M(t)$ is irreducible.  Let
$t\in T$ and let $M(t) = \tilde H\otimes_{\CC[X]} \CC v_t$ be the
corresponding principal series module for $\tilde H$.  The
{\it spherical vector} in $M(t)$ is
$${\bf 1}_t = \sum_{w\in W} q^{\ell(w)}T_w v_t.$$
Up to multiplication by constants this is the unique
vector in $M(t)$ such that $T_w{\bf 1}_t={\bf 1}_t$ for
all $w\in W$.  The following proposition uses the result
of Steinberg [???] to give a simple proof of [Ka] Proposition
1.20 and Lemma 2.3.

\prop Let $t\in T$. 
\smallskip\noindent
\item{(a)}  If $W_t=\{1\}$ and $v_{wt}$,
$w\in W$ is the basis of $M(t)$ defined in (???) then
$${\bf 1}_t = \sum_{z\in W} t(c_z),
\qquad\hbox{where}\qquad
c_z = \prod_{\alpha\in R(w_0z)} {q-q^{-1}X^\alpha\over 1-X^\alpha}.$$
\item{(b)} The spherical vector ${\bf 1}_t$ generates
$M(t)$ if and only if 
$\prod_{\alpha\in R^+} \left({???\over ???}\right) \ne 0$.
\item{(c)} The principal series module $M(t)$ is irreducible
if and only if ${\bf 1}_{wt}$ generates $M(wt)$ for all $w\in W$. 
\pf
Suppose that $\xi_z\in \CC$ are constants such that 
$${\bf 1}_t = \left( \sum_{w\in W} T_w\right)v_t
=\sum_{z\in W} \xi_z v_{zt}.$$
We shall prove that the $\xi_z$ are given by the formula in the
statement.  Since $T_i(\sum q^{\ell(w)}T_w) = q(\sum q^{\ell(w)}T_w)$,
$$\eqalign{
q{\bf 1}_t &= T_i{\bf 1}_t
=\left( \tau_i+{q-q^{-1}\over 1-X^{-\alpha_i}}\right)
\sum_{z\in W} \xi_zv_{zt} \cr
&=\left( \tau_i+{q-q^{-1}\over 1-X^{-\alpha_i}}\right)
\sum_{s_iz>z} \xi_zv_{zt}+\xi_{s_iz}v_{s_izt} \cr
&=\sum_{s_iz>z}
\xi_zv_{s_izt}+\xi_z{q-q^{-1}\over 1-t(X^{-z^{-1}\alpha_i})}v_{zt}
+\xi_{s_izt}\tau_i^2v_{zt}
+\xi_{s_iz}{q-q^{-1}\over 1-t(X^{z^{-1}\alpha_i})}v_{s_izt}.
\cr}$$
Comparing coefficients of $v_{s_izt}$ on each side of this expression
gives
$$q\xi_{s_iz} = \xi_z+\xi_{s_iz}{q-q^{-1}\over 1-t(X^{z^{-1}\alpha_i})},
\qquad\hbox{and so}\qquad
{\xi_z\over \xi_{s_iz}} 
= t\left({q^{-1}-qX^{z^{-1}\alpha_i}\over 1-X^{z^{-1}\alpha_i}}
\right),
\qquad\hbox{if $s_iz>z$.}$$
Using this formula inductively gives
$$\xi_w=\xi_{s_{i_1}\cdots s_{i_p}}
=t\left({1-X^{s_{i_p}\cdots s_{i_2}\alpha_{i_1}}\over
q^{-1}-qX^{s_{i_p}\cdots s_{i_2}\alpha_{i_1}} }\right)
\cdots
t\left({1-X^{\alpha_{i_p}}\over q^{-1}-qX^{\alpha_{i_p}} } \right)\xi_1
=t\left(\prod_{\alpha\in R(w)} 
{1-X^\alpha\over q^{-1}-qX^\alpha } \right)\xi_1.$$ 
Since the transition matrix between the basis $\{T_wv_t\}$ and the 
basis $\{v_{wt}\}$ is upper unitriangular with respect to the
Bruhat order, $\xi_{w_0}=1$.  Thus the last equation implies that
$$\eqalign{
\xi_1 &= t\left(\prod_{\alpha\in R^+} 
{q^{-1}-qX^\alpha\over 1-X^\alpha}\right)
\qquad\hbox{and}\qquad \cr
\xi_w &= t\left(\prod_{\alpha\in R(w)} 
{1-X^\alpha\over q^{-1}-qX^\alpha } \right)\xi_1
=t\left(\prod_{\alpha\in R(w_0w)} 
{q^{-1}-qX^\alpha\over 1-X^\alpha}\right). \cr
}
$$
(b)
By expanding $v_{zt} = \tau_zv_t = \tau_{i_1}\cdots \tau_{i_p}v_t$
for a reduced word $s_{i_1}\cdots s_{i_p}=z$ it follows that there exist
rational functions $m_{uz}$ such that
$$v_{zt} = \sum_{u\in W} t(m_{uz})T_uv_t,$$
for all generic $t\in T$.  Furthermore the matrix $M=(m_{uz})_{u,z\in W}$
with these rational functions as entries is upper unitriangular.

Let $\lambda_w$, $w\in W$, be as in ??? and define polynomials
$q_{uy}\in \CC[X]$, $u,y\in W$, by
$$X^{\lambda_y}\left(\sum_{w\in W} q^{\ell(w)}T_w\right)
=\sum_{u\in W} T_u q_{uy},
\qquad\qquad y\in W,$$
where these equations are equalities in $\tilde H$.
Then,
$$X^{\lambda_y}{\bf 1}_t
= X^{\lambda_y}\left( \sum_{w\in W} q^{\ell(w)}T_w\right)
=\sum_{u\in W} t(q_{uy})(T_u\otimes v_t),
$$
and part(a) implies that if $t$ is generic then
$$X^{\lambda_y}{\bf 1}_t
=X^{\lambda_y}\sum_{z\in W} t(c_z)v_{zt}
= \sum _{z\in W} t(c_zX^{z^{-1}\lambda_y}) v_{zt}
=\sum_{z,u\in W} t(c_zX^{z^{-1}\lambda_y}m_{uz})(T_u\otimes v_t).$$
Since these two expressions are equal for all generic
$t\in T$ it follows that 
$$q_{uy} = \sum_{z\in W} m_{uz} c_z (X^{z^{-1}\lambda_y}),
\qquad u,y\in W,
\formula$$
as rational functions (in fact both sides are in $\CC[X]$).

Since $T_w$, $w\in W$, and $p\in Z(\tilde H)=\CC[X]^W$
act on ${\bf 1}_t$ by constants, 
the $\tilde H$-module $M(t)$ is generated by ${\bf 1}_t$ if and only 
if there exist constants $p_{yw}\in \CC$ such that
$$T_w\otimes v_t = \sum_{y\in W} p_{yw}X^{\lambda_y} {\bf 1}_t,
\qquad\hbox{ for each $w\in W$.}$$
If these constant exist then, for each $w\in W$, 
$$T_w\otimes v_t = \sum_{y\in W} p_{yw}X^{\lambda_y} {\bf 1}_t
=\sum_{y,z,u\in W} t(m_{uz}c_zX^{z^{-1}\lambda_y} p_{yw}) 
~T_u\otimes v_t,$$
where, by ????, there is no restriction that $t$ be generic.
If
$$M=(m_{uz})_{u,z\in W},
\qquad
C = {\rm diag}(c_z)_{z\in W},
\qquad
X = (X^{z^{-1}\lambda_y}),
\qquad
P=(p_{yw})_{y,w\in W},$$
then $P=(t(MCX))^{-1}$ and so $P$ exists if and only if
$\det(t(MCX))\ne 0$.  Now $\det(M)=1$,
$$\det(X)=\xi\prod_{\alpha\in R^+} \alpha^{|W|/2},
\qquad\hbox{and}\qquad
\det(C) = \prod_{z\in W} 
\prod_{\alpha\in R(w)} {q^{-1}-qX^{-\alpha}\over 1-X^{-\alpha}}.$$
where $\xi\in \CC$ is nonzero.
So $P$ exists if and only if 
$$t\left(\prod_{\alpha\in R^+} (q^{-1}-qX^\alpha)\right) \ne 0.
$$
(c) $\Longrightarrow$:  If $M(t)$ is irreducible then, by Proposition
2.8(c), $M(wt)$ is irreducible for all $w\in W$.  
Hence $M(wt)$ is generated by ${\bf 1}_{wt}$.

\smallskip\noindent
$\Longleftarrow$:  Suppose that ${\bf 1}_{wt}$ generates $M(wt)$
for all $w\in W$.  Let $E$ be a nonzero
irreducible submodule of $M(t)$
and let $w\in W$ be such that the weight space $E_{wt}$ is nonzero.
Then, by Proposition ???(a), there is a nonzero surjective
$\tilde H$-module homomorphism $\varphi\colon M(wt)\to E$.  Since
${\bf 1}_{wt}$ generates $M(wt)$, $\varphi({\bf 1}_{wt})$ is a nonzero
vector in $E$ such that 
$T_v\varphi({\bf 1}_{wt})=q\varphi(v_{wt})$ 
for all $w\in W$.  Since there is a unique,
up to constant multiples, spherical vector in $M(t)$,
$\varphi({\bf 1}_{wt})$ is a multiple of ${\bf 1}_t$ and 
${\bf 1}_t$ is nonzero.  This implies that $E=M(t)$ since
${\bf 1}_t$ generates $M(t)$.
\endpf

Together the three parts of Proposition ??? prove the
following irreducibility criterion of Kato [Ka, Theorem 2.1].

\thm Let $t\in T$ and define 
$P(t) = \{\alpha>0\ |\ t(X^\alpha) = q^{\pm 2} \}$.
The principal series module 
$$\hbox{$M(t)$ is irreducible if and only if $P(t)=\emptyset$.}$$
\endthm
 
\medskip
\noindent
{\bf Remark.}  Kato actually proves a more general result and thus needs
a further condition for irreducibility.  We have
simplified matters by specifying the weight lattice $P$ in our 
construction of the
affine Hecke algebra.  One can use any $W$-invariant lattice in
$\RR^n$ and Kato works in this more general situation.  When the
one uses the weight lattice $P$, a result of Steinberg [St, 4.2, 5.3] 
says that the stabilizer $W_t$ of a point $t\in T$ under the action of 
$W$ is always a reflection group.
Because of this Kato's criterion takes a simpler form.

\subsection{Decomposition of the principal series module.}

Suppose $M(\rho)$ is the
principal series module corresponding to $\rho$. As an
$H$-module, $M(\rho)$ is isomorphic to the left regular representation. In
particular, $M(\rho)$ has basis $\{T_wv_\rho\}_{w\in W}$. In the Grothendieck ring
we have 
$$
\left[ M(\rho)\right] = \sum_{J\subseteq P(\rho)} \left[ H^{(\rho,J)}\right].$$
This identity is a module theoretic version of (?) for the case when
$\gamma=\rho$. 

In type $A$ the identity (??) has appeared in the form
$$
\CC\left[ S_n\right]\cong \bigoplus_{\theta} S^{\theta}
$$ 
where the sum is over all ribbons $\theta$ of length $n$ and
$S^{\theta}$ denote the corresponding skew shape representations of
$S_n$. 

Finally, the identity has also appeared in the symmetric function literature as 
$$
h_1^n = \sum_{\theta} s_{\theta}
$$
where the sum is over all ribbons $\theta$ of length $n$, $h_1=x_1+\cdots+x_n$
is the first complete symmetric function and
$s_{\theta}$ are skew Schur functions [Mac I \S5 Ex.~21b].

\subsection Row reading talbeaux

Let
Suppose $j>i$ and 
$\varepsilon_j-\varepsilon_i\in\overline{J}$. Then 
$\varepsilon_j-\varepsilon_i=
(\varepsilon_j-\varepsilon_{j_1})+(\varepsilon_{j_1}-\varepsilon_{j_2})
+\cdots+(\varepsilon_{j_r}-\varepsilon_i)$ where each of the
summands is in $J$.  This means that ${\rm box}_j$ is northwest of
${\rm box}_{j_1}$ is northwest of ${\rm box}_{j_2}$ is northwest of
$\ldots$ is notrthwest of ${\rm box}_i$.  Thus the 
position of ${\rm box}_j$ relative to ${\rm box}_i$ is somewhere in the shaded
region of 
$$
\beginpicture
\setcoordinatesystem units <0.75cm,0.75cm>         
\setplotarea x from 0 to 4, y from 0 to 5    
\linethickness=0.5pt                          
\put{\vdots} at 0 4.5 
\put{\vdots} at 4 4.5 
\put{${\scriptstyle{\rm box}_i}$} at 3.5 0.5 
\putrule from 1 3 to 2 3          
\putrule from 2 2 to 3 2          %
\putrule from 3 1 to 4 1          %
\putrule from 3 0 to 4 0          %
\putrule from 0 4 to 1 4        %
\putrule from 1 3 to 1 4        %
\putrule from 2 2 to 2 3        
\putrule from 3 0 to 3 2        %
\putrule from 4 0 to 4 4        %
\vshade 0 4 5   1 4 5 / 
\vshade 1 3 5   2 3 5 /
\vshade 2 2 5   3 2 5 /
\vshade 3 1 5   4 1 5 /           
\endpicture
$$
and
$\varepsilon_j-\varepsilon_i\in\overline{J}^c$ if ${\rm box}_j$ is in the
following shaded region.
$$
\beginpicture
\setcoordinatesystem units <0.75cm,0.75cm>         
\setplotarea x from 0 to 5, y from 0 to 6    
\linethickness=0.5pt                          
\put{\vdots} at 1 5.5 
\put{\vdots} at 3 0.5 
\put{${\scriptstyle{\rm box}_i}$} at 0.5 3.5 
\putrule from 0 4 to 1 4          %
\putrule from 0 3 to 1 3          
\putrule from 1 2 to 2 2          %
\putrule from 2 1 to 3 1          %
\putrule from 0 3 to 0 4        %
\putrule from 1 2 to 1 5        %
\putrule from 2 1 to 2 2        
\vshade 1 2 6   2 2 6  /
\vshade 2 1 6   3 1 6  /
\vshade 3 0 6   4 0 6 /           
\endpicture
$$

For the purposes of contradiction assume that
$j>i$, $l>k$, with 
$\varepsilon_j-\varepsilon_i,\varepsilon_l-\varepsilon_k\in\overline{J}^c$ and
$(\varepsilon_j-\varepsilon_i)+(\varepsilon_l-\varepsilon_k)
\in R^+\setminus\overline{J}^c$. 
In particular, either $j=k$ or $i=l$. 
Suppose $i=l$. Hence
$\varepsilon_j-\varepsilon_i,\varepsilon_i-\varepsilon_k\in\overline{J}^c$ and
$\varepsilon_j-\varepsilon_k\not\in\overline{J}^c$. 
In other words,
$\varepsilon_j-\varepsilon_i,\varepsilon_i-\varepsilon_k\in\overline{J}^c$ and
$\varepsilon_j-\varepsilon_k\in\overline{J}$. 
By considering their relative
positions in the shape we find this to be impossible. 
Similarly we cannot have $j=k$. 
Thus $\overline{J}^c$ is closed and $\cF^{(\gamma,J)}$ has a unique minimal
element.